\documentclass[12pt]{article}

\newlength{\spse}
\usepackage{bbold}
\usepackage{subfigure}
\usepackage{graphics}
\usepackage{amsmath}
\usepackage{epstopdf}
\usepackage{epsfig}
\usepackage{amsfonts}
\usepackage{verbatim}
\usepackage{multirow}
\newtheorem{Thm}{Theorem}

\newtheorem{Prop}[Thm]{Proposition}

\def\Real{\mathbb R}

\newcommand{\beq}{\begin{equation}}
\newcommand{\eeq}{\end{equation}}
\newcommand{\beqn}{\begin{equation*}}
\newcommand{\eeqn}{\end{equation*}}
\newcommand{\beqa}{\begin{eqnarray}}
\newcommand{\eeqa}{\end{eqnarray}}
\newcommand{\beqan}{\begin{eqnarray*}}
\newcommand{\eeqan}{\end{eqnarray*}}

\def\R{\mathbb{R}}

\title
{Study of  conservation and recurrence of Runge-Kutta discontinuous  Galerkin schemes for Vlasov-Poisson systems}

\author{  Yingda Cheng
\thanks{Department of Mathematics, Michigan State University,
East Lansing, MI 48824 U.S.A.
 {\tt ycheng@math.msu.edu}}
  \and
  Irene M. Gamba
\thanks{Department of Mathematics and ICES, University of Texas at Austin,
Austin, TX 78712 U.S.A.
 {\tt gamba@math.utexas.edu}}
 \and
 Philip J. Morrison
\thanks{Department of Physics and Institute for Fusion Studies, University of Texas at
Austin, Austin, TX 78712 U.S.A.
 {\tt morrison@physics.utexas.edu}}}

\date{\today}
\begin{document}

\maketitle

\begin{abstract}
In this paper we consider Runge-Kutta discontinuous  Galerkin (RKDG) schemes for  Vlasov-Poisson systems that  model collisionless plasmas.  One-dimensional systems are emphasized.  The RKDG method, originally devised to solve conservation laws,  is seen to have  excellent conservation properties,  be readily designed for arbitrary order of accuracy, and  capable of being used with a positivity-preserving limiter  that  guarantees  positivity of the distribution functions.  The RKDG solver for the Vlasov equation is the main focus,  while the electric field is obtained through the classical representation by Green's function for the Poisson equation.  A rigorous study of  recurrence of the DG methods  is presented by Fourier analysis, and the impact of different polynomial spaces and the positivity-preserving limiters on the quality of the solutions is ascertained.  Several  benchmark test problems, such as Landau damping, the two-stream instability,   and the KEEN (Kinetic Electro
 static Electron  Nonlinear) wave, are given.
\end{abstract}


\section{Introduction}
\label{intro}

The Vlasov-Poisson (VP) system is an  important  equation for  modeling  collisionless plasmas, one that possesses computational difficulties of more complete kinetic theories.  Thus, it serves as an important test bed for algorithm development.   The VP system describes the evolution of $f=f(x,v,t)$,  the probability distribution function ($pdf$) for  finding an electron (at position $x$ with velocity $v$ at time $t$) with a uniform background of fixed ions under a self-consistent electrostatic field. In particular, the non-dimensionalized VP system (with time scaled by the inverse plasma frequency $\omega_p^{-1}$ and length scaled by the Debye length $\lambda_D$) is given by
\beqa
\label{eq: vp}
\partial_t f + v \cdot \nabla_x f - E \cdot \nabla_v f=0 && \Omega \times (0, T] \nonumber \\
-\Delta_x \Phi=1- \int_{\Real^n} f \, dv && \Omega_x \times (0, T] \\
E= - \nabla_x \Phi  && \Omega_x \times (0, T] \,.  \nonumber
\eeqa
Here the domain $\Omega= \Omega_x \times \Real^n$, where $\Omega_x$ can be either a finite domain or $\Real^n$. The boundary conditions for the above systems are summarized as follows:  $f\rightarrow0$ as $|x|\rightarrow \infty$ or $|v|\rightarrow \infty$. If $\Omega_x$ is finite, then we can  impose either inflow boundary conditions with $f=f^{in}$ on $\Gamma_I=\{(x,v)| v \cdot \nu_x < 0\}$, where $\nu_x$ is the outward normal vector,  or more simply impose  periodic boundary conditions. For simplicity of discussion, in this paper, we will always assume periodicity in  $x$.  Also, we add that when the VP system is applied to plasmas the total charge neutrality  condition, $\int_{\Omega_x} \!\left(\int_{\Real^n} f \, dv -1 \right)dx=0$,  is imposed.

The following physical quantities associated with this system are related to its conservation properties:
\beqa
\label{eq:vp-conserve}
\textrm{charge density} & &\rho(x,t)= \int_{\Real^n} f(x,v,t) \,dv\,,  \nonumber\\
\textrm{momentum density} & &j(x,t)= \int_{\Real^n} v f(x,v,t) \,dv \,,  \qquad \\
\textrm{kinetic energy density} & &\xi_k(x,t)= \frac{1}{2} \int_{\Real^n} |v|^2 f(x,v,t) \,dv\,,  \nonumber\\
\textrm{electrostatic energy density} & &\xi_e(x,t)= \frac{1}{2} |E(x,t)|^2\,.\nonumber
\eeqa
Indeed, it is well-known that the VP system conserves the total electron charge $\int_{\Omega_x} \rho(x) \, dx$,    momentum $\int_{\Omega_x} j(x) \, dx$,  and energy $\int_{\Omega_x} (\xi_k(x)+\xi_e(x))\, dx$. Moreover, any functional of the form $\int_{\Omega} G(f)\, dx dv$ is a constant of motion. In particular, this includes the $k$-th order invariant $I_k = \int_{\Omega} f^k \,dx dv$ and the entropy $S= -\int_{\Omega} f ln(f) \, dx dv$.  Sometimes the functional $I_2$ is also called the enstrophy, and all of these invariants are called Casimir invariants (see, e.g., \cite{morrison_98}).

The VP system has been studied extensively for the simulation of collisionless plasmas. Popular numerical approaches include Particle-In-Cell (PIC) methods \cite{Birdsall_book1991, Hockney_book1981},   Lagrangian particle methods \cite{barnes_86,shadwick},   semi-Lagrangian methods  \cite{chengknorr_76, Sonnendrucker_1999}, the
WENO method coupled with Fourier collocation \cite{Guo_landau_2001}, finite volume (flux balance) methods  \cite{Boris_1976,Fijalkow_1999, Filbet_PFC_2001}, Fourier-Fourier spectral methods \cite{Klimas_1987, Klimas_1994}, continuous finite element methods \cite{Zaki_1988_1, Zaki_1988_2}, among many others. In this paper, we will focus on the  discontinuous Galerkin (DG) method  to solve the VP system. The original DG method  was introduced by Reed and Hill \cite{Reed_hill_73} for neutron transport.  Lesaint and Raviart \cite{Lesaint_r_74} performed the first convergence study for the original DG method. Cockburn and Shu in a series of papers
\cite{Cockburn_1991_MMNA_RK, Cockburn_1989_MC_RK_DG,Cockburn_1989_JCP_RK_localDG, Cockburn_1990_MC_RK_DG,Cockburn_1998_JCP_RK} developed the Runge-Kutta DG (RKDG) method for hyperbolic equations. The RKDG methods have  been used to simulate the VP system in plasmas by Heath, Gamba, Morrison and Michler \cite{Heath, Heath_thesis} and for the gravitational infinite homogeneous stellar system by Cheng and Gamba \cite{Cheng_jeans}. Theoretical aspects about stability, accuracy and conservation of those methods are discussed in \cite{Heath_thesis, Heath} and  more recently in \cite{Ayuso2009, Ayuso2010} for energy conserving schemes.
Such methods have excellent conservation properties, can be readily designed for arbitrary order of accuracy, and have the  potential for implementation  on unstructured meshes with $hp-$adaptivity. To ensure the positivity of the solution, one can  use a maximum-principle-satisfying limiter that has been recently proposed by Zhang and Shu in  \cite{Zhang_2010_Max} for conservation laws on cartesian meshes, and later extended on triangular meshes \cite{Zhang_JSC_triangle}. This limiter has been used to develop positivity-preserving schemes for compressible Euler \cite{Zhang_JCP_2010_Euler, Zhang_JCP_2011_Euler_source}, shallow water equations \cite{Zhang_Water_2010},  and Vlasov-Boltzmann transport equations \cite{cheng_mathcomp_10}. It has also been employed recently in the framework of semi-Lagrangian DG methods \cite{rossmanith_11, qiu_ppdg_11} for the VP system.

The scope of the present paper is as follows: we  focus on a detailed study of the RKDG scheme for the Vlasov equation from both the numerical and  analytical points of view. Since we are only considering one-dimensional problems, we use the classical representation of the solution by Green's function to compute the  Poisson equation; therefore,  the electric field is explicitly given as a function of the numerical density. This removes all discretization error from the Poisson equation and lets us more accurately investigate our DG solver for the Vlasov equations.
We rigorously study recurrence, which is an important numerical phenomenon that commonly appears with many solvers. We use Fourier analysis and obtain eigenvalues of the amplification matrix, and then  investigate the impact of different polynomial spaces on the quality of the solution by examining  conserved quantities  as well as convergence to  BGK states \cite{BGK} for some choices  of initial states. We consider benchmark test problems such as simulations of Landau damping phenomena for the linearized and nonlinear Vlasov Poisson systems, two-stream instability, and their long time  BGK states  and the formation of  KEEN  waves, both for  the  nonlinear system as well.

The remaining part  of the paper is  organized as follows: in Section \ref{nummeth}, we describe the numerical algorithm and summarize its conservation properties. In Section \ref{recurr}, we study the recurrence phenomena that occurs  for linear Vlasov type transport equations discretized by DG methods with various polynomial orders. Sections \ref{linearVP}  and \ref{NLVP}  are devoted to   discussions of  simulation results for the linearized and nonlinear VP system,  respectively,  for diverse choices of initial data and external drive potentials. We conclude with a few remarks in Section \ref{conclu}.


\section{Numerical methods}
\label{nummeth}

In this section we  first describe the proposed DG numerical algorithm and then discuss some of its basic conservation properties related to the quantities of  \eqref{eq:vp-conserve}.  This is done for both the fully nonlinear VP system of (\ref{eq: vp}) and  the linearized VP  system obtained by linearizing about  the  homogeneous equilibrium $f_{eq}(v)$,  with corresponding vanishing electric equilibrium field.  Periodic boundary conditions in  $x$ are assumed.

Thus, setting   $f(x, v, t)=f_{eq}(v)+\delta f (x, v,t)$  and  expanding  the system to first order approximation,  the perturbation $\delta f$ satisfies the Linear Vlasov-Poisson (LVP) system,
\beqa
\label{eq: linearvp}
\partial_t  f + v \cdot \nabla_x  f  =E \cdot \nabla_v f_{eq}  && \Omega \times (0, T] \nonumber\\
\Delta_x \Phi= \int_{\Real^n}  f \, dv && \Omega_x \times (0, T] \\
E= - \nabla_x \Phi  && \Omega_x \times (0, T] \,,
\nonumber
\eeqa
where  $\delta f$ has been replaced by $f$ to ease the notation.    We find it convenient and efficient to intertwine the  discussion of our algorithms for the  VP and LVP systems.  To avoid confusion in Section \ref{numalgo} we underline the words linear and nonlinear, signaling where discussions specific  to each apply.

\subsection{Numerical algorithm}
\label{numalgo}

For one-dimensional problems, we use a mesh that is a tensor product of grids in the $x$ and $v$ directions, because this simplifies the definitions  of the mesh and polynomial space for the Poisson equation. Specifically, the domain $\Omega$ is  partitioned as follows:
\begin{equation*}
\label{2dcell1} 0=x_{\frac{1}{2}}<x_{\frac{3}{2}}< \ldots
<x_{N_x+\frac{1}{2}}=L , \qquad -V_c=v_{\frac{1}{2}}<v_{\frac{3}{2}}<
\ldots <v_{N_v+\frac{1}{2}}=V_c,
\end{equation*}
where $V_c$ is chosen appropriately large to guarantee $f(x, v, t)=0$ for $|v| \geq V_c$. This is a reasonable assumption, because of the well-posedness of the one-dimensional Vlasov-Poisson system as indicated in \cite{glassey}.   The grid is defined as
\begin{eqnarray*}
\label{2dcell2} && I_{i,j}=[x_{i-\frac{1}{2}},x_{i+\frac{1}{2}}]
\times [v_{j-\frac{1}{2}},v_{j+\frac{1}{2}}] , \nonumber \\
&&J_i=[x_{i-1/2},x_{i+1/2}],
\quad K_j=[v_{j-1/2},v_{j+1/2}]\,, \     \quad i=1,\ldots N_x, \quad j=1,\ldots N_v ,
\end{eqnarray*}
where $x_i=\frac{1}{2}(x_{i-\frac{1}{2}}+x_{i+\frac{1}{2}})$ and $v_j=\frac{1}{2}(v_{j-1/2}+v_{j+1/2})$ are center points of the cells.

We will make use of several approximation spaces. For the $x$-domain, we consider the piecewise polynomial space  of functions $\xi:\Omega_x \rightarrow \R$,
\begin{equation*}
 Z_h^l=\{ \xi: \, \xi|_{J_{i}} \in
P^l(J_{i}), \,\, i=1,\ldots N_x \},
\end{equation*}
where $P^l(J_i)$ is the space of polynomials in one dimension of degree up to $l$.
For the $(x, v)$ space, we consider two approximation spaces of functions $\phi,\varphi:\Omega\rightarrow\R$,
\begin{equation*}
 V_h^l=\{ \phi : \, \phi|_{I_{i,j}} \in
\mathbb{Q}^l(I_{i,j}), \,\, i=1,\ldots N_x, \quad j=1,\ldots N_v\}
\end{equation*} and
\begin{equation*}
 W_h^l=\{ \varphi : \, \varphi|_{I_{i,j}} \in
\mathbb{P}^l(I_{i,j}), \,\, i=1,\ldots N_x, \quad j=1,\ldots N_v\}\,,
\end{equation*}
where $\mathbb{Q}^l(I_{i,j})=P^l(J_i)\otimes P^l(K_j)=span \{ x^{l_1} \,v^{l_2}, \forall \,0 \leq l_1 \leq l, \,0 \leq l_2 \leq l \}$ denotes all polynomials of degree at most $l$
in $x$ and $v$ on $I_{i,j}$, and $\mathbb{P}^l(I_{i,j})=span \{ x^{l_1} \,v^{l_2}, \forall \,0 \leq l_1+l_2 \leq l, l_1 \geq 0, l_2 \geq 0 \}$. These two spaces are  widely considered in the DG literature for multi-dimensional problems. A simple calculation shows that the number of degrees of freedom  of $\mathbb{Q}^l(I_{i,j})$ is $(l+1)^2$.  For  $l \geq 1$ this  is  larger than the number of degrees of freedom  of $\mathbb{P}^l(I_{i,j})$, which is $(l+1)(l+2)/2$.

First we describe the RKDG scheme for the \underline{linear} Vlasov equation. We   seek  $f_h(x,t) \in V_h^l$ (or $W_h^l$), such that
\begin{eqnarray}
\label{lv.scheme}
    \int_{I_{i,j}} (f_h)_t \varphi_h \, dx dv - \int_{I_{i,j}} v f_h (\varphi_h)_x \, dx dv + \int_{K_j} (\widehat{v f_h} \varphi_h^-)_{i+\frac{1}{2}, v} \, dv \nonumber \\
    -\int_{K_j} (\widehat{v f_h} \varphi_h^+)_{i-\frac{1}{2}, v} \, dv =\int_{I_{i,j}} E_h f^{'}_{eq} \varphi_h \, dx dv
\end{eqnarray}
holds for any test function $\varphi_h(x,t) \in V_h^l$ (or $W_h^l$).
Here and below, we use the following notations:  $E_h$ is the discrete electric field,  which is to be computed from  Poisson's equation,  $(\varphi_h)^\pm_{i+1/2, v}=\lim_{\epsilon \rightarrow 0} \varphi_h (x_{i+1/2}\pm \epsilon, v)$,   $(\varphi_h)^\pm_{x, j+1/2}=\lim_{\epsilon \rightarrow 0} \varphi_h (x, v_{j+1/2}\pm \epsilon)$, and $\widehat{v f_h}$ is a  numerical flux. We can assume that in each $K_j$, $v$ has a single sign by properly partitioning the mesh. Then, the  upwind flux is defined as
\begin{eqnarray}
\label{flux1}
  \widehat{v f_h}=  \left\{\begin{array} {l l}
        \displaystyle   v f_h^- & \textrm{if}\quad v \geq 0 \quad \textrm{in} \quad K_j\\
        \displaystyle    v f_h^+ & \textrm{if}\quad v < 0 \quad \textrm{in} \quad K_j
         \end{array}
   \right..
\end{eqnarray}

The scheme for the \underline{nonlinear} Vlasov equation is similar. We seek $f_h(x,t) \in V_h^l$ (or $W_h^l$), such that
\begin{eqnarray}
\label{lv.scheme2}
   && \int_{I_{i,j}} (f_h)_t \varphi_h \, dx dv - \int_{I_{i,j}} v f_h (\varphi_h)_x \, dx dv + \int_{K_j} (\widehat{v f_h} \varphi_h^-)_{i+\frac{1}{2}, v} \, dv
    -\int_{K_j} (\widehat{v f_h} \varphi_h^+)_{i-\frac{1}{2}, v} \, dv \nonumber \\
    &&+ \int_{I_{i,j}} E_h f_h (\varphi_h)_v \, dx dv
    - \int_{J_i} (\widehat{E_h f_h} \varphi_h^-)_{x, j+\frac{1}{2}} \, dx + \int_{J_i} (\widehat{E_h f_h} \varphi_h^+)_{x, j-\frac{1}{2}} \, dx=0
\end{eqnarray}
holds for any test function $\varphi_h(x,t) \in V_h^l$ (or $W_h^l$). The upwind flux for $\widehat{v f_h}$ has been defined in (\ref{flux1}) and the new  flux needed  for the  nonlinear case  is given by
\begin{eqnarray}
\label{flux2}
  \widehat{E_h f_h}=  \left\{\begin{array} {l l}
        \displaystyle   E_h f_h^- & \textrm{if}\quad \int_{J_i} E_h dx \leq 0 \\
        \displaystyle    E_h f_h^+ & \textrm{if}\quad \int_{J_i} E_h dx > 0
         \end{array}
   \right..
\end{eqnarray}

 The above descriptions coupled with a suitable time discretization scheme, such as  the TVD Runge-Kutta method \cite{Shu_1988_JCP_NonOscill},    completes the RKDG methods. For example, suppose the semi-discrete schemes in (\ref{lv.scheme}) and (\ref{lv.scheme2}) are written in the compact form
$$
\int_{I_{i,j}} (f_h)_t \varphi_h \, dx dv=\mathcal{H}_{i,j} (f_h, E_h, \varphi_h)
$$
where for the \underline{linear} Vlasov of (\ref{lv.scheme})
\beqan
\mathcal{H}^{lin}_{i,j} (f_h, E_h, \varphi_h)&=&   \int_{I_{i,j}} v f_h (\varphi_h)_x \,dx dv - \int_{K_j} (\widehat{v f_h} \varphi_h^-)_{i+\frac{1}{2}, v} \, dv
    \nonumber \\
    &+&   \int_{K_j} (\widehat{v f_h} \varphi_h^+)_{i-\frac{1}{2}, v} \,dv +\int_{I_{i,j}} E_h f^{'}_{eq} \varphi_h \,dx dv\,,
\eeqan
while for the \underline{nonlinear} Vlasov of (\ref{lv.scheme2})
\beqan
\mathcal{H}^{nonlin}_{i,j} (f_h, E_h, \varphi_h)&=&  \int_{I_{i,j}} v f_h (\varphi_h)_x \,dx dv - \int_{K_j} (\widehat{v f_h} \varphi_h^-)_{i+\frac{1}{2}, v} \, dv
    +\int_{K_j} (\widehat{v f_h} \varphi_h^+)_{i-\frac{1}{2}, v} \, dv \nonumber \\
  &-&  \int_{I_{i,j}} E_h f_h (\varphi_h)_v \, dx dv
    + \int_{J_i} (\widehat{E_h f_h} \varphi_h^-)_{x, j+\frac{1}{2}} \, dx - \int_{J_i} (\widehat{E_h f_h} \varphi_h^+)_{x, j-\frac{1}{2}} \,dx\,.
\eeqan
The third order TVD-RK method implements  the following procedure for going from  $t^n$ to $t^{n+1}$:
\beqa
\label{rk3}
&&\int_{I_{i,j}} f_h^{(1)} \varphi_h \,dx dv=\int_{I_{i,j}} f_h^n \varphi_h \,dx dv
+ \triangle t\,  \mathcal{H}_{i,j} (f_h^n, E_h^n, \varphi_h)\,,\nonumber \\
&&\int_{I_{i,j}} f_h^{(2)} \varphi_h \,dx dv=\frac{3}{4}\int_{I_{i,j}} f_h^n \varphi_h \, dx dv+\frac{1}{4}\int_{I_{i,j}} f_h^{(1)} \varphi_h \,dx dv+\frac{\triangle t}{4} \mathcal{H}_{i,j} (f_h^{(1)}, E_h^{(1)}, \varphi_h)\,, \\
&&\int_{I_{i,j}} f_h^{n+1} \varphi_h \,dx dv=\frac{1}{3}\int_{I_{i,j}} f_h^n \varphi_h \,dx dv+\frac{2}{3}\int_{I_{i,j}} f_h^{(2)} \varphi_h \,dx dv+\frac{2 \triangle t}{3} \mathcal{H}_{i,j} (f_h^{(2)}, E_h^{(2)}, \varphi_h) \,.\nonumber
\eeqa

Poisson's equation is used to obtain $E_h^n$, $E_h^{(1)}$,  and $E_h^{(2)}$.  Beyond periodicity, we need to enforce some additional conditions to uniquely determine $\Phi$. For example, we can set  one end  of the spatial domain to ground, i.e.\  set $\Phi(0,t)=0$. In the one-dimensional case, then the exact solution can  be   obtained if we enforce $\Phi(0)=\Phi(L)$.  For the \underline{nonlinear}  system we obtain
$$
\Phi_h=\int_0^x \int_0^s \rho_h(z,t) \, dz  ds -\frac{x^2}{2}-C_E  x,
$$
where $C_E= -{L}/{2}+   \int_0^L \int_0^s \rho_h(z,t) \, dz ds/L$, and
\begin{equation}
\label{efield}
E_h=-(\Phi_h)_x=C_E+x-\int_0^x \rho_h(z,t) \, dz\,,
\end{equation}
while for the \underline{linear} system Poisson's equation gives
$$
\Phi_h=\int_0^x \int_0^s \rho_h(z,t) dz \, ds -C_E x,
$$
where $C_E=   \int_0^L \int_0^s \rho_h(z,t) \, dz  ds/L$, and
\begin{equation}
\label{efield2}
E_h=-(\Phi_h)_x=C_E-\int_0^x \rho_h(z,t)\,  dz\,.
\end{equation}

{}From (\ref{efield}) and (\ref{efield2}), we  see that if $f_h \in V_h^l$ (or $W_h^l$), then $\rho_h=\int_{-V_c}^{V_c} f_h \, dv \in Z_h^l$; hence,  $E_h \in Z_h^{l+1} \bigcap C^0$. The above approach  uses  the classical representation of the solution by Green's function and will be referred to as the ``exact" Poisson solver. It is valid only  for the  one-dimensional case. For higher dimensions, a suitable elliptic solver needs to be implemented, such as those discussed in \cite{Heath}. Here we use the exact solver to entirely eliminate discretization error from   Poisson's  equation and, thereby,  spotlight  the performance of the Vlasov solver.

Below we  describe  positivity-preserving limiters,  as summarized in \cite{Zhang_mpp_review}.  We only use such a limiter to enforce the positivity of $f_h$ for the \underline{nonlinear} VP system.  For  each of the forward Euler steps of the Runge-Kutta time discretization, the
following procedures are performed:
\begin{itemize}
\item On each cell $I_{i,j}$, we evaluate $T_{i,j}:=\min_{(x,v) \in
S_{i,j}} f_h(x,v)$, where $S_{i,j}=(S_i^x \otimes \hat{S_j^y}) \bigcup (\hat{S_i^x} \otimes  S_j^y)$, and $S_i^x, S_j^y$ denote the $(l+1)$ Gauss quadrature points, while $\hat{S_i^x}, \hat{S_j^y}$ denote the $(l+1)$ Gauss-Lobatto quadrature points.
\item We compute $\tilde{f_h}(x,v)=\theta (f_h(x,v)-(\overline{f_h})_{i,j})+
(\overline{f_h})_{i,j}$, where $(\overline{f_h})_{i,j}$ is the cell average
of $f_h$ on $I_{i,j}$, and $\theta=\min \{1 ,
|(\overline{f_h})_{i,j}|/|T_{i,j}-(\overline{f_h})_{i,j}|\}$. This limiter has
the effect of maintaining  the cell average, while ``squeezing" the function
to be positive at all points in $S_{i,j}$.
\item Finally, we use $\tilde{f_h}$ instead of $f_h$ to compute the Euler
forward step.
\end{itemize}

This completes the description of the numerical algorithm.

\subsection{Scheme Conservation properties}
\label{cons}
 In the following,  we will briefly review and discuss some of the conservation properties of the above RKDG scheme for the nonlinear VP equations without the positivity-preserving limiter.
Some of those results have been reported in \cite{Heath_thesis} and \cite{Ayuso2009}.

\begin{Prop}(charge conservation)  For both the $V_h^l$ and $W_h^l$ spaces,
\[
\sum_{i,j} \mathcal{H}^{nonlin}_{i,j} (f_h, E_h, 1)=\Theta (f_h, E_h, 1)\,
\]
which implies
\begin{eqnarray}
\sum_{i,j} \int_{I_{i,j}} f_h^{n+1} \,  dx dv&=&\sum_{i,j} \int_{I_{i,j}} f_h^{n} \, dx dv+\frac{2}{3}\triangle t \big(\Theta (f_h^{(2)}, E_h^{(2)}, 1)\nonumber\\
&{\ }&\qquad \quad+ \frac{1}{4}\Theta (f_h^{(1)}, E_h^{(1)}, 1)+\frac{1}{4}\Theta (f_h^n, E_h^n, 1)\big)
\nonumber
\end{eqnarray}
for the fully discrete scheme (\ref{rk3}). Here,
 \[
 \Theta (f_h, E_h, \varphi_h)= \sum_i \int_{J_i} (\widehat{E_h f_h} \varphi_h)_{x, N_v+\frac{1}{2}}\,  dx
 -  \sum_i \int_{J_i} (\widehat{E_h f_h} \varphi_h)_{x, \frac{1}{2}}\, dx
 \]
denotes the contribution from the phase space boundaries located at $v=\pm V_c$, and should be negligible if $V_c$ is chosen large enough.
\end{Prop}

\textbf{Remark}:  Charge conservation (or mass conservation or probability normalization as it is sometimes called)    states that the total charge will be preserved on the discrete level up to approximation errors associated with the phase space boundaries. The proof is straightforward and, therefore, omitted. The same conclusion can be proven for the linearized system. The positivity preserving limiter does not destroy this property because it keeps the cell averages unchanged.

\begin{Prop}(Semi-discrete $L^2$ stability --  enstrophy decay)  For both the $V_h^l$ and $W_h^l$ spaces, $\sum_{i,j} \mathcal{H}^{nonlin}_{i,j} (f_h, E_h, f_h) \leq 0$. Hence, $$\frac{d}{dt} \sum_{i,j} \int_{I_{i,j}} f_h^2 dx dv \leq 0.$$
\end{Prop}
The proof, for an arbitrary field $E_h$, can be found in  \cite{cheng_mathcomp_10}, Theorem~4, which applies directly  here
 by setting the collisional form $Q_\sigma\equiv 0$ in that proof.

\medskip

For the remainder of  this subsection we will assume the DG solution satisfies the velocity boundary conditions $f_h(x, \pm V_c, t)=0$. This is a reasonable assumption when $V_c$ is large enough. In particular, this will guarantee exact charge conservation, which implies that $\int_0^L \rho_h(x, t) dx $ is constant in  time $t$. Therefore, using the definition of $E_h$ in (\ref{efield}), we can obtain  periodicity in  $E_h$, i.e, $E_h(0)=E_h(L)$. Without this assumption  the propositions below  contain multiple boundary terms and the proof  becomes  technical.

\begin{Prop}(Momentum conservation) Assuming  $f_h(x, \pm V_c, t)=0$, for both the $V_h^l$ and $W_h^l$ spaces when $l \geq 1$,  $\sum_{i,j} \mathcal{H}^{nonlin}_{i,j} (f_h, E_h, v)=0$, which implies
$$
\sum_{i,j} \int_{I_{i,j}} f_h^{n+1} v  \, dx dv=\sum_{i,j} \int_{I_{i,j}} f_h^{n} v \, dx dv
$$
for the fully discrete scheme.
\end{Prop}
\emph{Proof.} Choosing  $\varphi_h=v$ in (\ref{lv.scheme2}), we have
\beqan
\sum_{i,j} \mathcal{H}^{nonlin}_{i,j} (f_h, E_h, v)&=&
\sum_{i,j}  \left ( \int_{I_{i,j}} v f_h (v)_x \, dx dv - \int_{K_j} (\widehat{v f_h} v)_{i+\frac{1}{2}, v} \, dv
    +\int_{K_j} (\widehat{v f_h} v)_{i-\frac{1}{2}, v}\,  dv  \right.\nonumber \\
    && \left . \quad - \int_{I_{i,j}} E_h f_h \,  dx dv
    + \int_{J_i} (\widehat{E_h f_h} v)_{x, j+\frac{1}{2}} \, dx - \int_{J_i} (\widehat{E_h f_h} v)_{x, j-\frac{1}{2}}\,  dx \right )\\
    &=&- \sum_{i,j} \int_{I_{i,j}} E_h f_h  \, dx dv = -\sum_i \int_{J_i} E_h \rho_h \, dx\,,
\eeqan
and using the exact Poisson solver together with the periodicity of $E_h$ and  $\Phi_h$ yields the following:
\beqan
  \sum_i \int_{J_i} E_h \rho_h\,  dx &=& \sum_i \int_{J_i} E_h (\rho_h -1) \, dx + \sum_i \int_{J_i} E_h \,  dx\\
 &=& -\sum_i \int_{J_i} E_h (E_h)_x \, dx + \sum_i \int_{J_i} E_h \,  dx\\
 &=&-(E_h^2(L)-E_h^2(0))/2 -\Phi(L)+\Phi(0)=0\,,
 \eeqan
which completes the proof. $\Box$

\medskip

\textbf{Remark}: The above  proof holds  for the linearized system as well. Note, however, it relies on the use of  the exact Poisson solver.  For a full numerical DG  Poisson solver,  such as that developed in  \cite{Heath}  for the discretization Poisson equation, exact momentum conservation remains true,  as was proven in   \cite{Heath_thesis}  by means of the DFUG method developed there for dealing with the discretized Poisson equation.
However, the positivity-preserving limiter we use here  will destroy this property because it was not designed to
conserve the numerical momentum.

\begin{Prop}(Semi-discrete total energy equality)  Assuming  $f_h(x, \pm V_c, t)=0$, for both the $V_h^l$ and $W_h^l$ spaces when $l \geq 2$,
$$
\frac{d}{dt} \left ( \frac{1}{2} \sum_{i,j} \int_{I_{i,j}} f_h v^2 \, dx dv +\frac{1}{2} \sum_i \int_{J_i} E_h^2(x)\,  dx \right )=A(f_h, \Phi_h)= A(f_h-f, \Phi_h-P(\Phi_h))\,,
$$
where the operator $A(w, u):= \sum_{i,j} \int_{I_{i,j}}\big(w u_x v -(w)_t u\big)  \, dx dv,$ and  $P$ is any projection such that $P(\Phi_h) \in Z_h^l$ and $P(\Phi_h)=\Phi_h$ at $x_{i + 1/2}$, for  $i=0, 1, \ldots, N_x$.
\end{Prop}
\emph{Proof.} Choosing  $\varphi_h= {v^2}/{2}$ in (\ref{lv.scheme2}) yields
\beqan
\frac{d}{dt}  \sum_{i,j}  \frac{1}{2} \int_{I_{i,j}} f_h v^2 \, dx dv + \sum_{i,j} \int_{I_{i,j}} E_h f_h v\, dx dv=0
\eeqan
and
\beqan
 \frac{d}{dt}  \sum_i  \frac{1}{2}\int_{J_i} E_h^2(x) \, dx &=& \sum_i  \int_{J_i} E_h (E_h)_t \, dx= -\sum_i  \int_{J_i} (\Phi_h)_x (E_h)_t \, dx\\
&=&\sum_i  \int_{J_i} \Phi_h (E_h)_{xt} \, dx=\sum_i  \int_{J_i} \Phi_h (1-\rho_h)_t\,  dx\\
&=&-\sum_i  \int_{J_i} \Phi_h (\rho_h)_t \, dx=- \sum_{i,j}  \int_{I_{i,j}} \Phi_h (f_h)_t \, dx dv\,,
\eeqan
where in the second line, we have used the periodicity and continuity of $E_h$ and $\Phi_h$. Therefore, we have proven  that
 $$
 \frac{d}{dt} \left ( \frac{1}{2} \sum_{i,j} \int_{I_{i,j}} f_h v^2 dx dv +\frac{1}{2} \sum_i \int_{J_i} E_h^2(x) dx \right )=A(f_h, \Phi_h)\,.
 $$
 On the other hand, upon choosing  $\varphi_h=P(\Phi_h)$ in (\ref{lv.scheme2}) and using the periodicity and continuity of $P(\Phi_h)$, we can verify that $A(f_h, P(\Phi_h))=0$. The exact solution $f$ obviously satisfies $A(f, \Phi_h-P(\Phi_h))=0$ from the continuity and periodicity of $\Phi_h-P(\Phi_h)$, and therefore we are done. $\Box$

 \medskip

The above proof indicates that the variation in the total energy will be related to the error of  the solution, $f_h-f$,  together with the projection error, $\Phi_h-P(\Phi_h)$.  In \cite{Heath_thesis, Heath}, error estimates for DG schemes with NIPG methods for the Poisson equations are provided for multiple dimensions. In \cite{Ayuso2009}, optimal accuracy of order $l+1$ for the semi-discrete scheme with $\mathbb{Q}^l$ spaces has been proven  under certain regularity assumptions.   We   remark that in \cite{Ayuso2009} conservation of the  total numerical energy is proven when the Poisson equation is solved by a local DG method with a special flux. Unfortunately,  no numerical simulations of linear Landau damping or of the nonlinear VP system, such as those done in \cite{Heath} or in Section \ref{numerical} of this present manuscript, have been performed up to this date by the scheme proposed in \cite{Ayuso2009}, so a comparison is not possible.


\section{On recurrence}
\label{recurr}

In this section we study recurrence,  a  numerical phenomenon that is known to occur  in simulations of  Vlasov-like equations.  Its study is important because it provides information about the numerical accuracy of a scheme. Recurrence  was  observed in the  semi-Lagrangian scheme of Cheng and Knorr  \cite{chengknorr_76}, where  a simple argument for its occurrence was provided.  In this section, we  carry out a detailed study of recurrence for  the DG method.

We study  recurrence of  our algorithm applied to the linear advection equation  $f_t+v f_x=0$ on $[0, L=2\pi/k] \times [-V_c, V_c]$,  since it  is tractable and contains the basic recurrence mechanism;  results for the full Vlasov system tend to be quite similar.  The initial condition we consider is  $f_0(x, v)=A \cos (kx) f_{eq} (v)$, and the equilibrium distribution $f_{eq}(v)$ is  taken to be  either the Maxwellian or Lorentzian distribution, viz.
$$
f_M=\frac{1}{\sqrt{2 \pi}} e^{-v^2/2}\qquad {\rm or}\qquad  f_L=\frac{1}{\pi} \frac{1}{v^2+1}\,.
$$
 For the Maxwellian equilibrium, $f_M$, we take $V_c=5$, and for the Lorentzian equilibrium, $f_L$,  we take $V_c=30$.

The exact solution for the advection equation is  $f(t,x,v)=f_0(x-vt, v)$. Hence,   a simple calculation shows $\rho(x,t)= A \cos(kx) e^{-k^2t^2/2}$ for the Maxwellian distribution; similarly,  for the Lorentzian,   $\rho(x,t)= A \cos(kx) e^{-kt}$.  Thus,  we see how the density for each case should decay  to zero. The failure of decay and the occurrence of recurrence noted for the   semi-Lagrangian scheme of \cite{chengknorr_76} stems from the  finite resolution in the velocity space  and, indeed,  the   recurrence time  depends on the mesh size in $v$.

To be specific, we repeat the definition of DG scheme for this equation, which amounts to (\ref{lv.scheme2}) with $E_h$ set to zero: we find $f_h(x,t) \in V_h^k$ (or $W_h^k$) , such that
\begin{eqnarray}
\label{adv.scheme}
    \int_{I_{i,j}} (f_h)_t \varphi_h \, dx dv - \int_{I_{i,j}} v f_h (\varphi_h)_x \, dx dv + \int_{K_j} (\widehat{v f_h} \varphi_h^-)_{i+\frac{1}{2}, v} \, dv
    -\int_{K_j} (\widehat{v f_h} \varphi_h^+)_{i-\frac{1}{2}, v} \, dv =0
\end{eqnarray}
holds for any test function $\varphi_h(x,t) \in V_h^k$ (or $W_h^k$). Again $\widehat{v f_h}$ is  the upwind numerical flux of  (\ref{flux1}). In the analysis below, we always assume time to be continuous, because recurrence is mainly a phenomenon that comes from the spatial and velocity discretization.


\subsection{The case of $l=0$}
\label{elliszero}

For the piecewise constant case, the DG method is equivalent to a simple first order finite volume scheme and we can derive rigorously the behavior for $\rho$. Suppose we define $f_h=f_{ij}$ on cell $I_{ij}$, and assume uniform grids  in both directions. Moreover, we assume $N_v$  to be even for simplicity. With this assumption,  the location of the cell center  is $v_j=(j-\frac{N_v+1}{2}) \triangle v$. Now (\ref{adv.scheme}) simply becomes
\beq
\label{p0scheme}
\begin{array}{ll}
\frac{d f_{ij}}{dt}+ v_j \frac{f_{ij}-f_{i-1,j}}{\triangle x} =0 & \textrm{if} \quad v_j  \geq 0,\\ \\
\frac{d f_{ij}}{dt}+ v_j   \frac{f_{i+1, j}-f_{ij}}{\triangle x}=0 & \textrm{if} \quad  v_j  < 0.
\end{array}
\eeq
The initial condition chosen is clearly equivalent to $f_{ij}(0)= Re\left( A e^{ikx_i} f_{eq}(v_j)\right)$. We  prove that the scheme above gives
\beq
\label{fsol}
f_{ij}(t)= Re\big( A e^{ikx_i+s_j t} f_{eq}(v_j)\big)\,
\eeq
where $s_j$ is given in (\ref{sj}) below.

Upon plugging   (\ref{fsol}) into (\ref{p0scheme}), we have
\beqn
\begin{array}{ll}
s_j f_{ij}+ v_j \frac{1-e^{-ik\triangle x}}{\triangle x} f_{ij} =0 & \textrm{if} \quad v_j  \geq 0\\ \\
s_j f_{ij}+ v_j   \frac{e^{ik\triangle x}-1}{\triangle x} f_{ij}=0 & \textrm{if} \quad  v_j  < 0\,.
\end{array}
\eeqn
Hence,
\beqn
s_j=\left \{
\begin{array}{ll}
- v_j \frac{1-e^{-ik\triangle x}}{\triangle x}= v_j \frac{\cos(k \triangle x)-1}{\triangle x} -v_j \frac{\sin(k \triangle x)}{\triangle x}\,  \mathrm{i}
 & \textrm{if}  \quad v_j  \geq 0\\ \\
- v_j   \frac{e^{ik\triangle x}-1}{\triangle x} = -v_j \frac{\cos(k \triangle x)-1}{\triangle x} -v_j \frac{\sin(k \triangle x)}{\triangle x} \, \mathrm{i}
 & \textrm{if} \quad  v_j  < 0\,,
\end{array}
\right.
\eeqn
which can be summarized as
\beq
\label{sj}
s_j= |v_j| \frac{\cos(k \triangle x)-1}{\triangle x} -v_j \frac{\sin(k \triangle x)}{\triangle x}\, \mathrm{i}\,.
\eeq
 Therefore,  the real part of $s_j$ is always negative, this means the magnitude of $f_{ij}$ will always damp, but because of the $j$-dependence it does so at  different rates for different cells.   Since the density
$$
\rho(x_i)=\sum_j f_{ij} \triangle v= Re\left( \sum_j  A e^{ikx_i+s_j t} f_{eq}(v_j)\right) \triangle v\,,
$$
the density will damp at a rate between $\frac{\triangle v}{2} \frac{\cos(k \triangle x)-1}{\triangle x}$ and $\frac{V_c-\triangle v}{2} \frac{\cos(k \triangle x)-1}{\triangle x}$. Another important fact is that recurring local maxima of the density will have a period $T_R$ that  satisfies $ \frac{\triangle v}{2} \frac{\sin(k \triangle x)}{\triangle x} T_R =\pi$. If we define $k'=\frac{\sin(k \triangle x)}{\triangle x} $, then $T_R= \frac{2 \pi}{k' \triangle v}$.  When $\triangle x \rightarrow 0$, $k' \rightarrow k$, and this coincides with the recurrence time obtained in \cite{chengknorr_76}.

\begin{figure}[htb]
  \begin{center}
    \subfigure[ Maxwellian, $\mathbb{P}^0$]{\includegraphics[width=3in,angle=0]{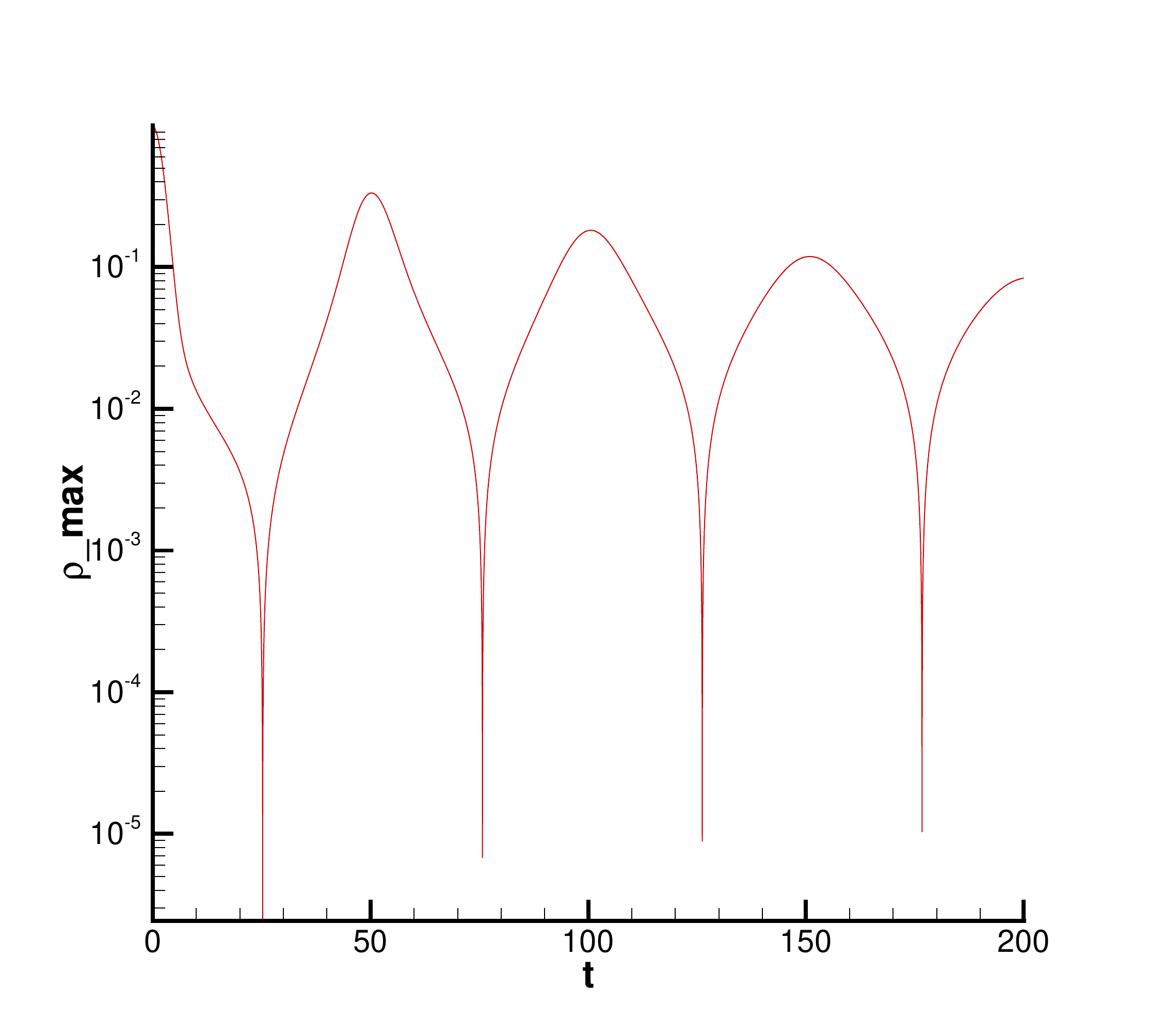}}
    \subfigure[Lorentzian, $\mathbb{P}^0$]{\includegraphics[width=3in,angle=0]{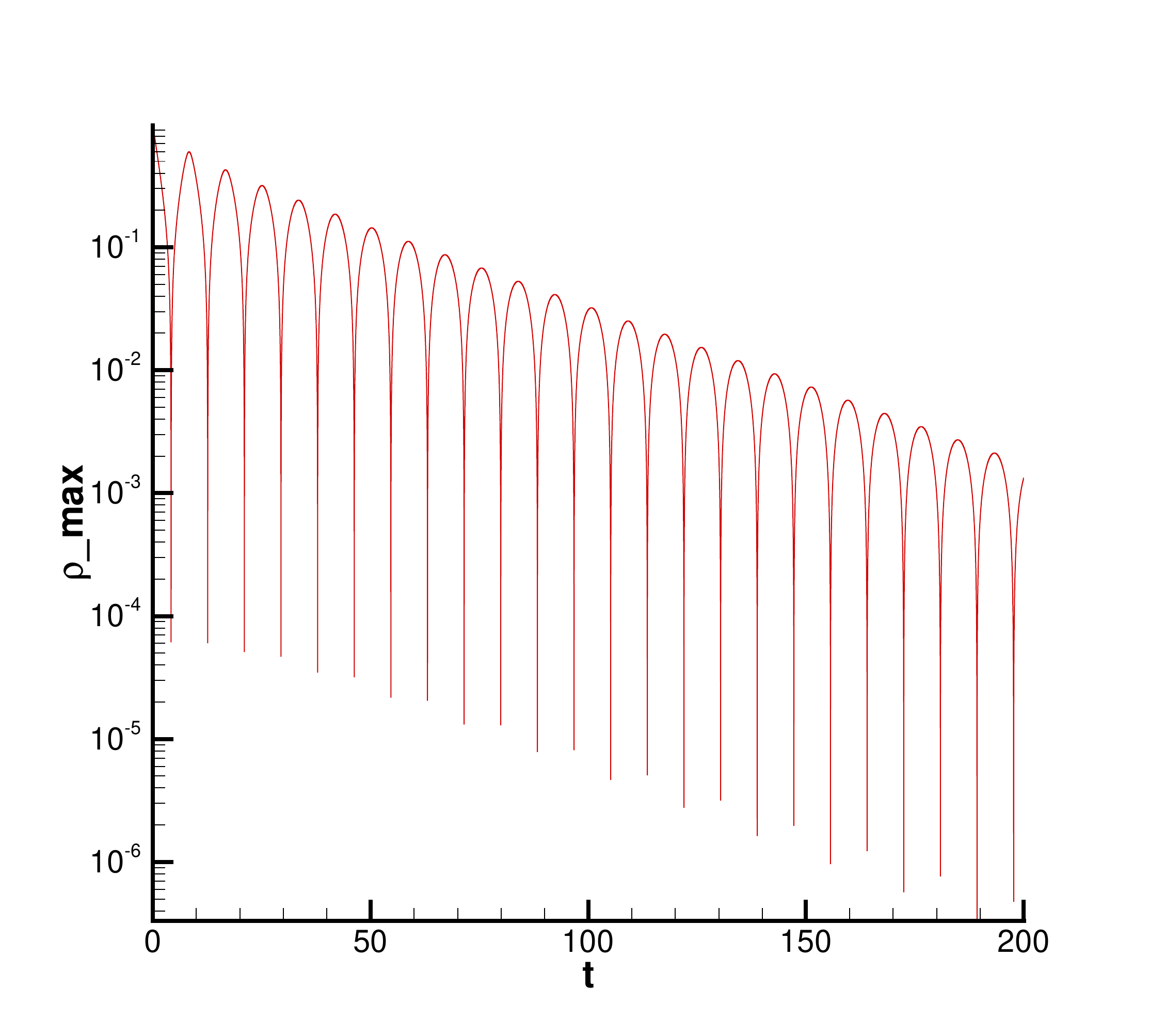}}
  \end{center}
  \caption{Computations of the advection equation for piecewise constant polynomials showing  local  maxima of the  density $\rho_{max}$ as a function of time. The mesh is $40 \times 40$ with $\triangle x=\pi/10$. For the  Maxwellian equilibrium $\triangle v=1/4$, while  for the Lorentzian equilibrium $\triangle v=3/2$.}
\label{recurp0}
\end{figure}

Next we compare the above theoretical prediction with numerical results. In all of the calculations below, we take $A=1$, $k=0.5$ and the mesh size to be $40 \times 40$.
In Figure \ref{recurp0}, we display results for  numerical runs using   piecewise constant polynomials and time discretization using  TVD-RK3. (We use the third order method to minimize  the  time discretization error.)  We plot
$\rho_{max}(t)=\max_x \rho(x, t)$  in Figure \ref{recurp0}.   First, we notice the pattern of $\rho_{max}$ has the expected  periodic structure with damping for  both Maxwellian and  Lorentzian equilibria.
For the Maxwellian distribution,  a simple calculation yields $T_R=50.47$. Similarly, with the formulas above, the smallest damping rate is $-0.49 \times 10^{-2}$, while the biggest is $-9.3\times 10^{-2}$.  For the Lorentzian distribution,  $T_R=8.41$, and  the smallest damping rate is $-2.94 \times 10^{-2}$, while the biggest is $-5.58 \times 10^{-1}$. From Figure \ref{recurp0}, by using the second to the fourth peak, the actual computed value of $T_R$ for Maxwellian  is 50.32 and the  damping rate is  $-1.02 \times 10^{-2}$; while for Lorentzian,  from the second to the tenth peak, $T_R$ is 8.40 and the  computed damping rate is$-3.19 \times 10^{-2}$. These numbers agree well with the theoretical prediction.

In conclusion, our  analysis explains the behavior of the first order DG solution. At $t=T_R$, the numerical density  obtains  a local maximum;  hence,  clearly at this time the numerical solution can no longer be trusted.  The numerical decay deviates from the theoretical decay well before  $t=T_R$.   To achieve a larger $T_R$, according to the formula, we can refine $\triangle v$.  On the other hand, refining $\triangle x$ will not change $T_R$ by much.

\medskip

\textbf{Remark:} Using the same methodology, it is easy to perform a similar analysis for any type of finite difference method. The real part of $s_j$ will be negative if there is numerical dissipation, and the imaginary part will always approximate  $v_j k$ due to the differential operator $v \frac{\partial}{\partial x}$. This means that for such schemes, the recurrence time $T_R$ will always be close to $\frac{2 \pi}{k \triangle v}$.

\subsection{Higher order polynomials}

In this subsection, we consider higher order polynomials. For the $V_h^1$ space, it takes  four point values in each cell to represent a $\mathbb{Q}^1$ polynomial. This technique was developed in \cite{Zhang_2005_CF_DG} for analyzing piecewise linear DG solutions in one dimension. As in Section \ref{elliszero}, we use  a uniform mesh, i.e. $\triangle x_i \equiv \triangle x$ and $\triangle v_j \equiv \triangle v$. Without loss of generality, we consider (\ref{adv.scheme}) for the case of $v \geq 0$ only, then $\widehat{ v f_h}=v f_h^-$, which means we only consider cells $I_{i,j}$ with $j \geq \frac{N_v}{2}+1$.

In each computational cell $I_{i,j}$, we can always use the following form to represent $f_h$:
$$
f_h=f_{i-\frac{1}{4}, j+\frac{1}{4}} \, \chi_1(x,v)+f_{i-\frac{1}{4}, j-\frac{1}{4}}\, \chi_2(x,v)+f_{i+\frac{1}{4}, j+\frac{1}{4}} \, \chi_3(x,v)+f_{i+\frac{1}{4}, j-\frac{1}{4}}\, \chi_4(x,v),
$$
where
\begin{eqnarray*}
\chi_1(x,v)&=&- 4 \left( \frac{x-x_i}{\triangle x_i}-\frac{1}{4}\right)\left( \frac{v-v_j}{\triangle v_j}+\frac{1}{4}\right) \\
\chi_2(x,v)&=& 4 \left( \frac{x-x_i}{\triangle x_i}-\frac{1}{4})\right)\left( \frac{v-v_j}{\triangle v_j}-\frac{1}{4}\right) \\
\chi_3(x,v)&=& 4 \left( \frac{x-x_i}{\triangle x_i}+\frac{1}{4}\right)\left( \frac{v-v_j}{\triangle v_j}+\frac{1}{4}\right) \\
\chi_4(x,v)&=& -4 \left( \frac{x-x_i}{\triangle x_i}+\frac{1}{4}\right)\left( \frac{v-v_j}{\triangle v_j}-\frac{1}{4}\right)
\end{eqnarray*}
are the basis functions in $\mathbb{Q}^1$ and $f_{i\pm1/4, j\pm 1/4}=f_h(x_{i\pm1/4}, v_{j\pm 1/4})$ are the point values. By choosing the test function  in (\ref{adv.scheme}) to be $\varphi_h=\chi_{\ell}$, $\ell=1, 2, 3, 4$, we obtain four  relations. Letting  $f_{ij}= (f_{i-1/4, j+1/4}, f_{i-1/4, j-1/4}, f_{i+1/4, j+1/4}, f_{i+1/4, j-1/4})^T$, then corresponding to each of the four terms in (\ref{adv.scheme}), we have
$$
M \frac{d f_{ij}}{dt}  - B  f_{ij}+ C f_{i,j}- D f_{i-1,j}=0,
$$
where
 \begin{eqnarray}
M &=&\frac{\triangle x \triangle v}{144} \begin{pmatrix} 49&-7 &-7&1\\ -7 &49&1&-7 \\ -7 &1&49&-7\\ 1&-7&-7 &49 \end{pmatrix}\,,
\nonumber \\
B &=&\frac{\triangle v}{12} \begin{pmatrix} -(2\triangle v +7 v_j) & v_j &-(2\triangle v +7 v_j)&v_j\\
v_j &2\triangle v -7 v_j&v_j&2\triangle v -7 v_j \\
 2\triangle v +7 v_j&-v_j&2\triangle v +7 v_j&-v_j\\
 -v_j&-2\triangle v +7 v_j&-v_j &-2\triangle v +7 v_j
 \end{pmatrix}\,,
 \nonumber \\
C &=&\frac{\triangle v}{48} \begin{pmatrix} 2\triangle v +7 v_j & -v_j &-(6\triangle v +21 v_j)&3 v_j\\
-v_j &-2\triangle v +7 v_j&3 v_j&6\triangle v -21 v_j \\
 -6\triangle v -21 v_j& 3 v_j&18 \triangle v +63 v_j&-9 v_j\\
 3 v_j&6 \triangle v -21  v_j&-9 v_j &-18\triangle v +63 v_j
 \end{pmatrix}\,,
 \nonumber
\end{eqnarray}
and
\begin{equation}
D= \frac{\triangle v}{48} \begin{pmatrix} -6\triangle v -21 v_j & 3 v_j &18 \triangle v +63 v_j &-9 v_j\\
3 v_j & 6\triangle v -21 v_j&-9 v_j& -18\triangle v +63 v_j \\
 2\triangle v +7 v_j& - v_j&-6 \triangle v -21 v_j&3 v_j\\
 - v_j&-2 \triangle v +7  v_j&3 v_j & 6\triangle v -21 v_j
 \end{pmatrix}\,.
 \nonumber
\end{equation}

After simple algebraic manipulation, we obtain
$$
\frac{d f_{ij}}{dt} =\frac{\triangle v}{\triangle x} \left(S_m f_{ij}+ T_m f_{i-1,j}\right),
$$
where
 \begin{eqnarray*}
 S_m &=& \begin{pmatrix} -\frac{49}{96}-\frac{7m}{8} & \frac{7}{96} &-\frac{7}{32}-\frac{3m}{8}&\frac{1}{32}\\
-\frac{7}{96} & \frac{49}{96}-\frac{7m}{8}&-\frac{1}{32}& \frac{7}{32}-\frac{3m}{8}\\
\frac{77}{96}+\frac{11 m}{8}& - \frac{11}{96}&-\frac{21}{32}-\frac{9m}{8}&\frac{3}{32}\\
 \frac{11}{96}&-\frac{77}{96}+\frac{m}{8}&-\frac{3}{32} & \frac{21}{32}-\frac{9m}{8} \end{pmatrix}\,, \\
 T_m &=& \begin{pmatrix}  -\frac{35}{96}-\frac{5m}{8} & \frac{5}{96} &\frac{35}{32}+\frac{15m}{8}&-\frac{5}{32}\\
-\frac{5}{96} & \frac{35}{96}-\frac{5m}{8}&\frac{5}{32}& -\frac{35}{32}+\frac{15m}{8}\\
\frac{7}{96}+\frac{ m}{8}& - \frac{1}{96}&-\frac{7}{32}-\frac{3m}{8} &\frac{1}{32}\\
 \frac{1}{96}&-\frac{7}{96}+\frac{m}{8}&-\frac{1}{32} & \frac{7}{32}-\frac{3m}{8} \end{pmatrix}\,,
\end{eqnarray*}
 and $m=2 j-N_v-1=1, 3, 5 \ldots$ are positive and odd integers.
 Therefore, the amplification matrix is given by
 $$
 G_j=\frac{\triangle v}{\triangle x}\left(S_m+T_m e^{-i k \triangle x}\right).
 $$

With the initial condition
 $f_{ij}(0)=Re( A e^{i k x_i} \Upsilon)$, where
 $$
 \Upsilon=\left(e^{-i k \triangle x/4} f_{eq}(v_{j+\frac{1}{4}}), e^{-i k \triangle x/4} f_{eq}(v_{j-\frac{1}{4}}), e^{i k \triangle x/4} f_{eq}(v_{j+\frac{1}{4}}), e^{i k \triangle x/4} f_{eq}(v_{j-\frac{1}{4}})\right)^T\,,
 $$
it is clear  that the general expression for the numerical solution is
 $$
 f_{ij}(t)=Re\left( e^{i k x_i} \sum_{\alpha=1}^4a_{\alpha}  V_{\alpha}\, e^{\eta_{\alpha} t}\right)\,.
 $$
Here  $\eta_{\alpha}$ are the eigenvectors of $G_j$ with    $V_{\alpha}$ the corresponding eigenvectors,    $a_{\alpha}$ are constants such that $f_{ij}(0)=\Upsilon$, and all these quantities  are  dependent on $j$ (or equivalently  $m=2 j-N_v-1$).
The collective behaviors of the eigenvalues $\eta_{\alpha}$ will influence the behavior of the density as a function of time. We focus on the  matrix $\Lambda_m=S_m+T_m e^{-i k \triangle x}$, which with  some algebraic manipulation can be seen to have the form
$$
\Lambda_m=W \otimes V,
$$
where $W$ and  $V$ are the following $2 \times 2$ matrices:
\begin{eqnarray*}
W &=& \begin{pmatrix}-3m-\frac{7}{4} & \frac{1}{4}\\
-\frac{1}{4} & -3m+\frac{7}{4} \end{pmatrix}\,, \\
V &=&\begin{pmatrix}\frac{1}{2}+\frac{5}{24}\, \hat{i} & -\frac{1}{2}-\frac{5}{8}\, \hat{i}\\
-\frac{1}{2}-\frac{1}{24}\,\hat{i} & \frac{1}{2}+\frac{1}{8}\, \hat{i} \end{pmatrix}\,,
\end{eqnarray*}
and $\hat{i}=e^{-ik\triangle x}-1= -i k \triangle x +O(\triangle x^2)$.
This nice structure is due to the tensor product formulations of the mesh and the space $\mathbb{Q}^l$. We   compute the eigenvalues of the matrix $V$, obtaining  $\lambda_1=(3+\hat{i} - \sqrt{9+12\hat{i}+\hat{i}^2})/6=\frac{1}{6} i k \triangle x - \frac{1}{12}  k^2 \triangle x^2 +O(\triangle x^3)$ and $\lambda_2=(3+\hat{i}+ \sqrt{9+12\hat{i}+\hat{i}^2})/6=1-\frac{1}{2} i k \triangle x +O(\triangle x^2)$, and the eigenvalues of $W$, obtaining  $-3m \pm \sqrt{3}$. Hence, by simple linear algebra, the four eigenvalues of the matrix $\Lambda_m$ are obtained
 $$
\begin{pmatrix} \xi_1= (-3m - \sqrt{3}) \lambda_2\\
\xi_2= (-3m + \sqrt{3})  \lambda_2\\
\xi_3= (-3m - \sqrt{3}) \lambda_1\\
\xi_4= (-3m + \sqrt{3}) \lambda_1
\end{pmatrix}\,.
$$
It is  easy to show that the eigenvectors corresponding to these eigenvalues are independent of $m$, since the eigenvectors of $V$ and $W$ are independent of $m$. We conclude that  the eigenvalues of $G_j$ are
$$
 \begin{pmatrix} \eta_1= (-3m - \sqrt{3})  \lambda_2 \triangle v/\triangle x\\
\eta_2= (-3m + \sqrt{3}) \lambda_2 \triangle v/\triangle x  \\
\eta_3= (-3m - \sqrt{3})  \lambda_1 \triangle v/\triangle x\\
\eta_4= (-3m + \sqrt{3}) \lambda_1 \triangle v/\triangle x\end{pmatrix}\,,
$$
and therefore,
\beqan
 \sum_{\alpha=1}^4a_{\alpha}  V_{\alpha}\, e^{\eta_{\alpha} t}
&=& e^{-3m (\lambda_2 \triangle v/\triangle x) t} \left(a_1 V_1\,  e^{- \sqrt{3} (\lambda_2 \triangle v/\triangle x)t}
+a_2 V_2 \, e^{ \sqrt{3} (\lambda_2 \triangle v/\triangle x)t}\right)\\
&& \hspace{.60 in} +e^{-3m (\lambda_1 \triangle v/\triangle x) t} \left(a_3 V_3\, e^{- \sqrt{3} (\lambda_1 \triangle v/\triangle x) t}
+a_4 V_4\, e^{ \sqrt{3} (\lambda_1 \triangle v/\triangle x) t}\right)\,.
\eeqan

Since $\eta_1$ and $\eta_2$ have nontrivial negative real parts, the damping for those two modes will be strong. Consequently, the  main behavior of the density will be dominated by the eigenmodes of $\eta_3$ and $\eta_4$. Recall
$$
\rho(x_{i \pm \frac{1}{4}}, t)=\sum_j \int_{I_j} f_h(x_{i \pm \frac{1}{4}}, v, t) dv= \sum_j  (f_{i \pm \frac{1}{4}, j+\frac{1}{4}}+f_{i \pm \frac{1}{4}, j-\frac{1}{4}}) \triangle v\,,
$$
and,  therefore, the behavior of  $\rho(x_{i \pm {1}/{4}}, t) $ is dominated by
$$
\sum_m e^{-3m (\lambda_1 \triangle v/\triangle x) t}
\left(c_3 \, e^{- \sqrt{3} (\lambda_1 \triangle v/\triangle x) t}
+c_4\,  e^{ \sqrt{3} (\lambda_1 \triangle v/\triangle x) t}\right)\,,
$$
where $c_3$ and  $c_4$ are constants that do not depend on $m$. Since  $\lambda_1=\frac{1}{6} i k \triangle x - \frac{1}{12}  k^2 \triangle x^2 +O(\triangle x^3)$, we have $-3m (\lambda_1 \triangle v/\triangle x) =-\frac{k \triangle v}{2}m i-\frac{m\triangle v \triangle x k^2}{4} + O(\triangle v \triangle x^2)$.
Hence,  with  an argument similar to that of Section \ref{elliszero} for the piecewise constant case, when $t \approx T_R=\frac{2 \pi}{k \triangle v}$ the imaginary parts of all modes will return to $m \pi i$,  and this will correspond to a local maximum of $\rho_{max}$ as a function of time. The remaining term, $c\,  e^{- \sqrt{3} i k \triangle v/6 t}
+d \, e^{ \sqrt{3}i k \triangle v/6 t}$,  corresponds to the envelope of the wave, and the negative real part of the eigenvalues indicates numerical dissipation.

\begin{figure}[htb]
  \begin{center}
    \subfigure[ Maxwellian, $\mathbb{Q}^1$]{\includegraphics[width=3in,angle=0]{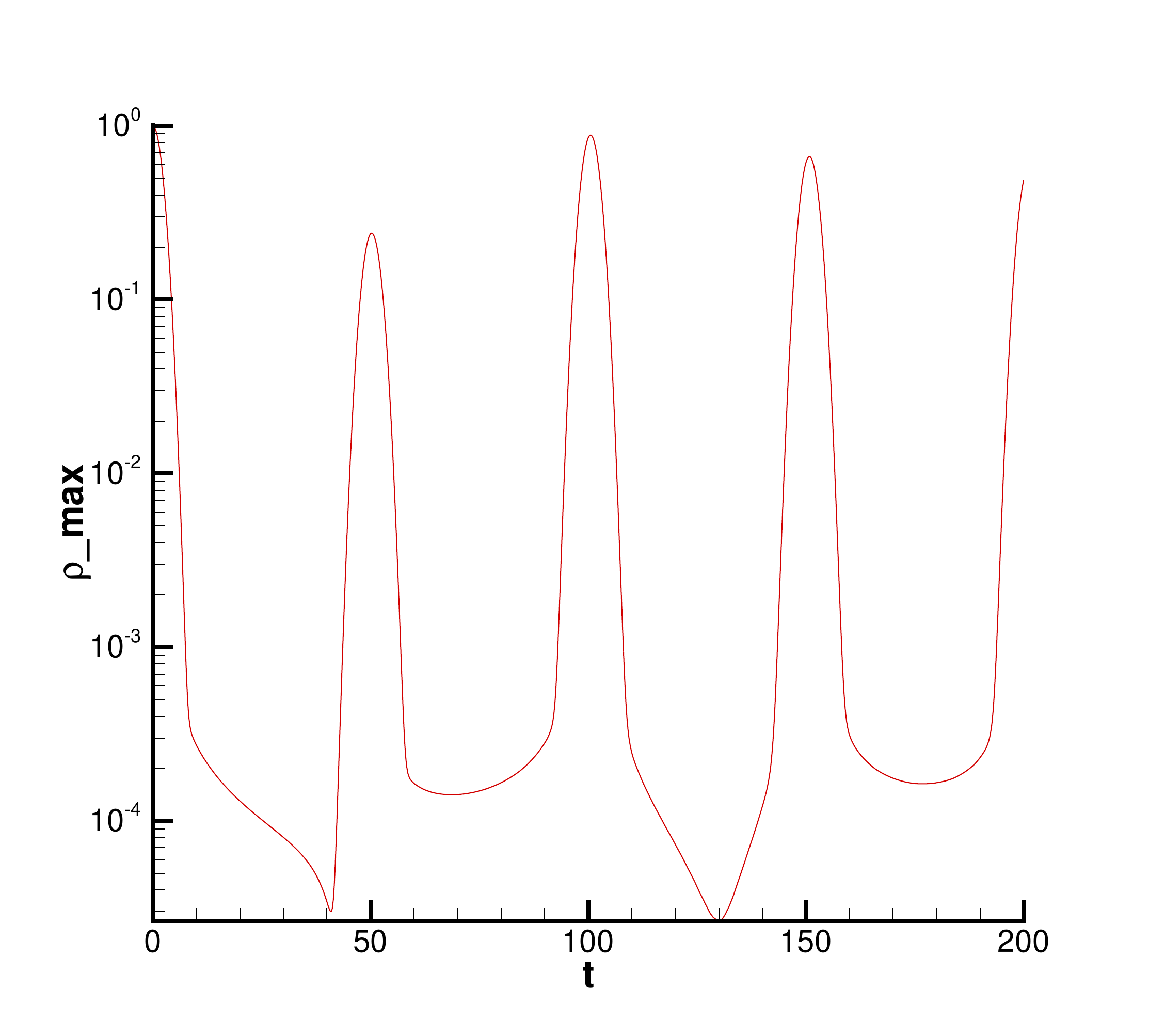}}
    \subfigure[Lorentzian, $\mathbb{Q}^1$]{\includegraphics[width=3in,angle=0]{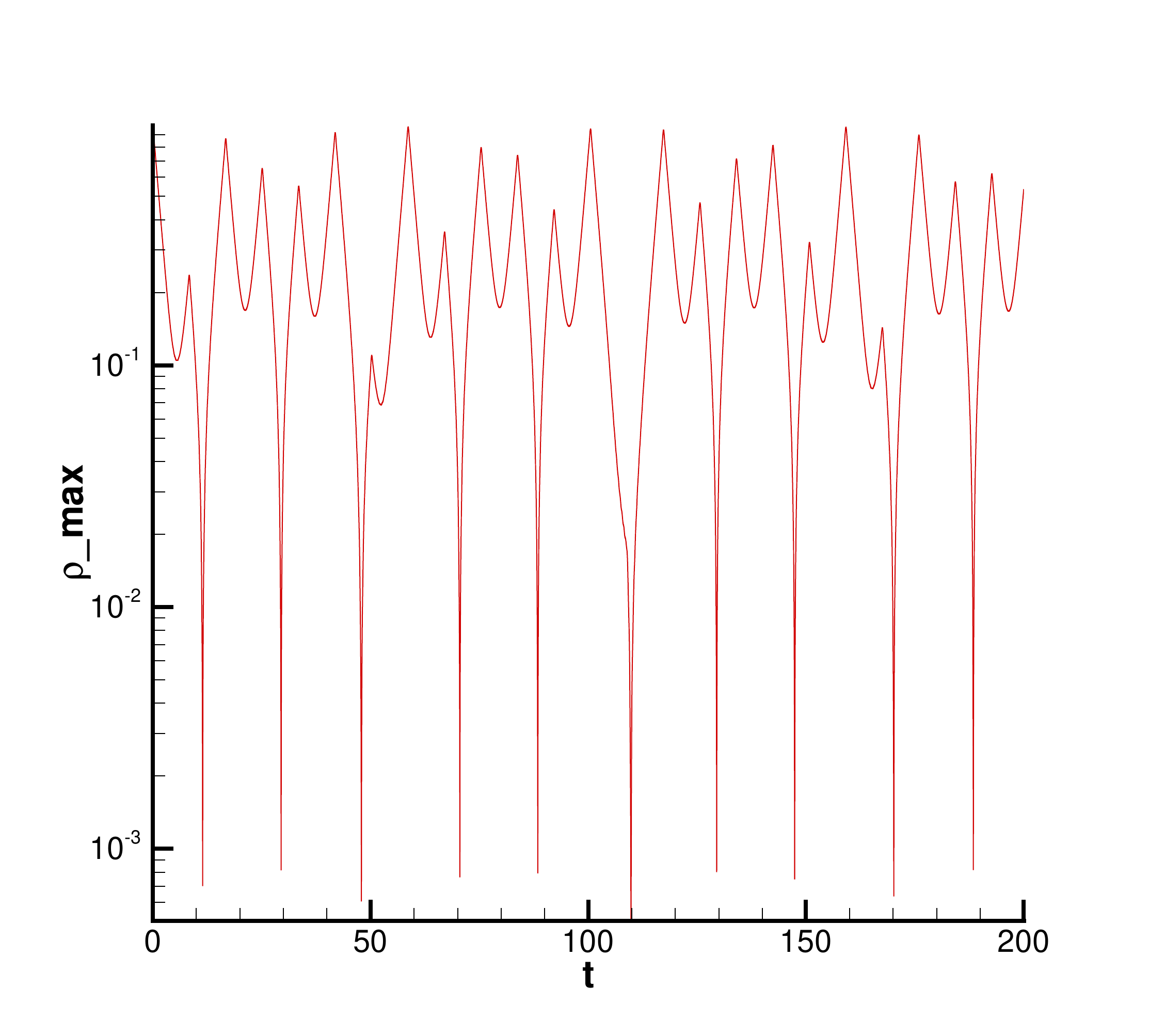}}
    \subfigure[Maxwellian, $\mathbb{Q}^2$]{\includegraphics[width=3in,angle=0]{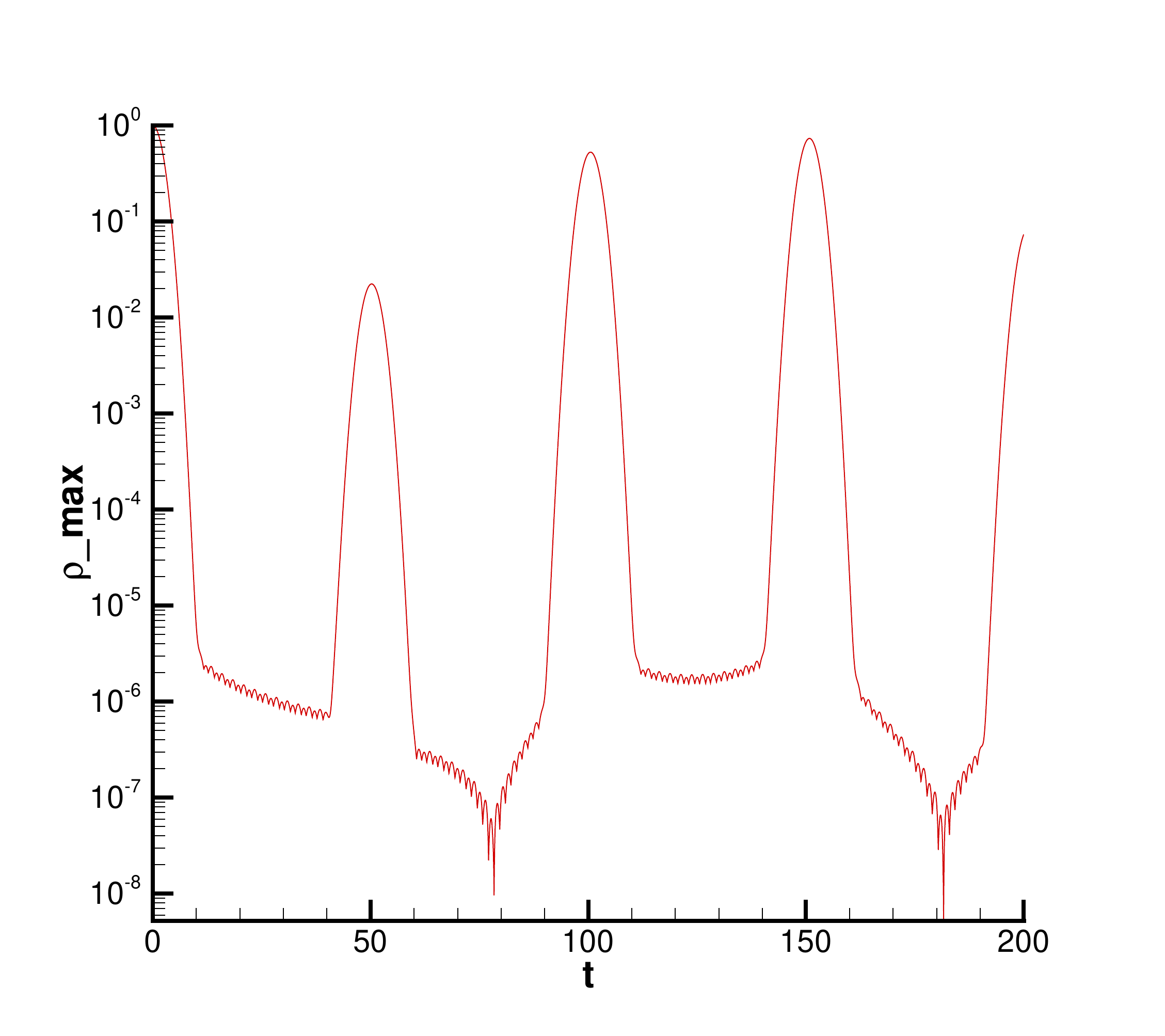}}
    \subfigure[Lorentzian, $\mathbb{Q}^2$]{\includegraphics[width=3in,angle=0]{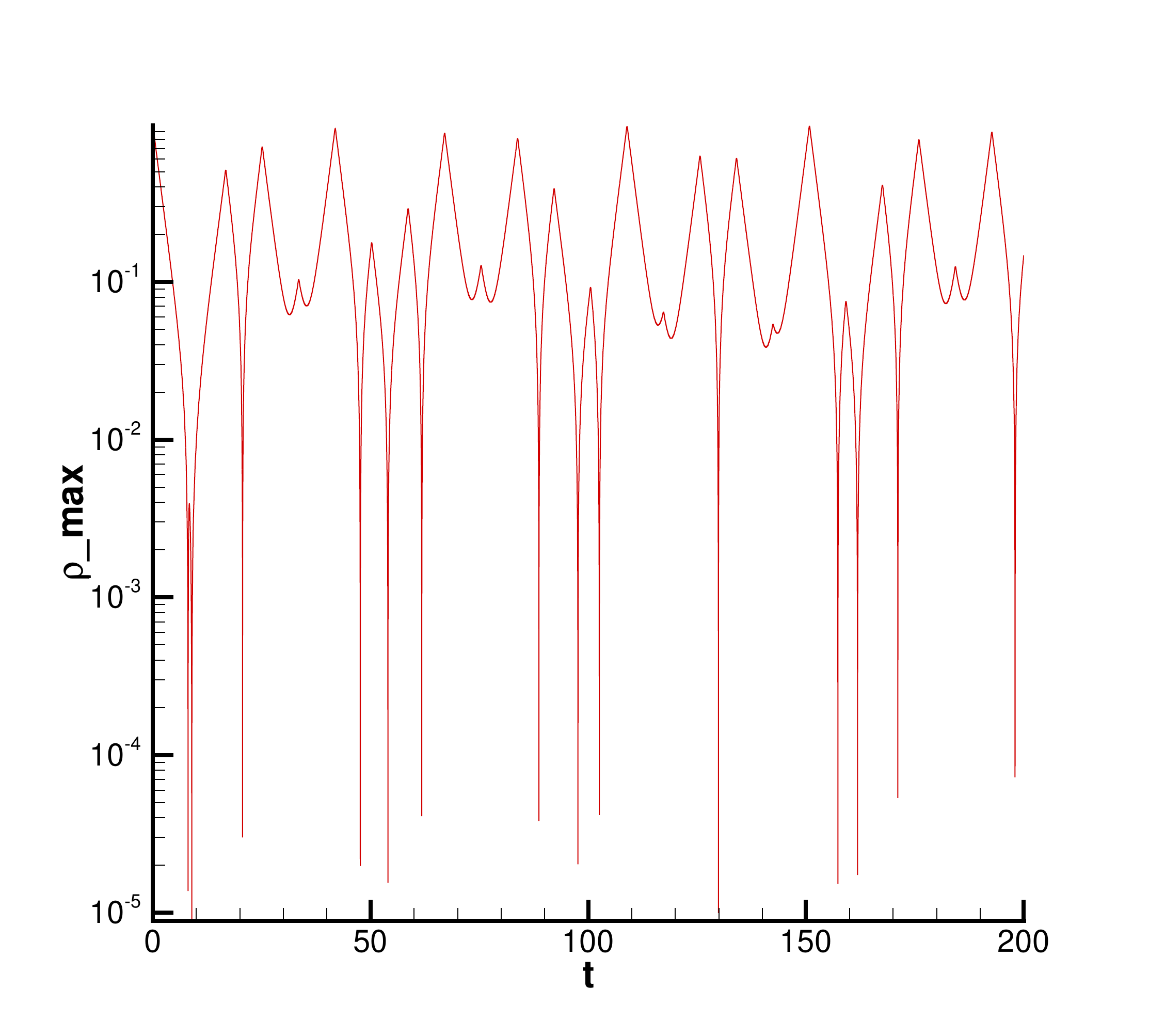}}
  \end{center}
  \caption{Computations of the advection equation for the  polynomial spaces  $\mathbb{Q}^1$ and $\mathbb{Q}^2$ showing  local maxima of the  density $\rho_{max}$ as a function of time.  The mesh is $40 \times 40$ with    $\triangle x=\pi/10$.  For the Maxwellian equilibrium $\triangle v=1/4$, while for the Lorentzian equilibrium $\triangle v=3/2$.}
\label{recurq}
\end{figure}

In Figure \ref{recurq}, we plot the evolution of $\rho_{max}$ as a function of time for  the $\mathbb{Q}^1$ and $\mathbb{Q}^2$ spaces. From Figures \ref{recurq}(a) and  \ref{recurq}(b), we observe the  behavior  predicted by our  analysis for the $\mathbb{Q}^1$ elements. From Figures \ref{recurq}(c) and  \ref{recurq}(d), we find that the solutions using the  $\mathbb{Q}^2$ polynomials share similar structures, except that small oscillations can be observed for the Maxwellian case. Also we note  that the $\mathbb{Q}^2$ discretizations can follow the exact solutions longer in time, in the sense that the minimum value achieved before $\rho_{max}$  starts to  deviates from the exact solution is  on the order of $10^{-6}$ compared to $10^{-4}$ in the $\mathbb{Q}^1$ case. This is expected due to the higher order accuracy of the scheme.  For  the $\mathbb{Q}^2$ polynomials, we deduce that the amplification matrix $G$ is a $9\times9$ matrix.  Thus, for this case there,   there will be nine  eigenvalues and  more modes than for the $\mathbb{Q}^1$ space. In Table \ref{recurtime}, we verify  the recurrence time $T_R$ numerically;  good agreement between the predicted values and the observed values are seen.

\begin{table} [htb]
\begin{center}
\caption {The location of local maxima of the density $\rho_{max}$ compared with the predicted  recurrence time $T_R$. The $T_R$ values for the Maxwellian equilibrium are computed using the average of the first three peaks, while for  the Lorentzian they are   computed using the average of the first seven peaks.}
\bigskip
\begin{tabular}{|c|c|c|c|}
\hline  &Predicted $T_R=\frac{2 \pi}{k \triangle v}$ & Numerical value of $\mathbb{Q}^1$ & Numerical value of  $\mathbb{Q}^2$ \\
\hline Maxwellian  &50.26548245743669  &50.265482457450& 50.265482457450 \\
\hline Lorentzian  & 8.37758040957278  & 8.37787960887962& 8.37787960887760\\
\hline
\end{tabular}
\label{recurtime}
\end{center}
\end{table}

Note, the trace $Tr(S_m)=-4 m, m=1, 3, 5 \ldots$, while a similar calculation for cells when $v <0$ yields $Tr(S_m)=4 m,  m=1, 3, 5 \ldots$. Therefore, we conclude that our semi-discrete algorithm has an incompressible vector field and thus possesses a version of Liouville's theorem on conservation of phase space volume.  Liouville's theorem  for  finite difference and Fourier discretization of fluid and plasma equations is well known  and has been used in statistical theories of turbulence (see e.g.\ \cite{burgers,lee, kraichnan80,jung}).   We also note that we have performed the analysis for the semi-discrete DG schemes. For fully discrete RKDG schemes, one could use a similar method,  as proposed in \cite{Zhong_11} for the wave equation, to write the fully discrete amplification matrices, but we do not pursue this in this paper.

\begin{figure}[htb]
  \begin{center}
        \subfigure[ Maxwellian, $\mathbb{P}^1$]{\includegraphics[width=3in,angle=0]{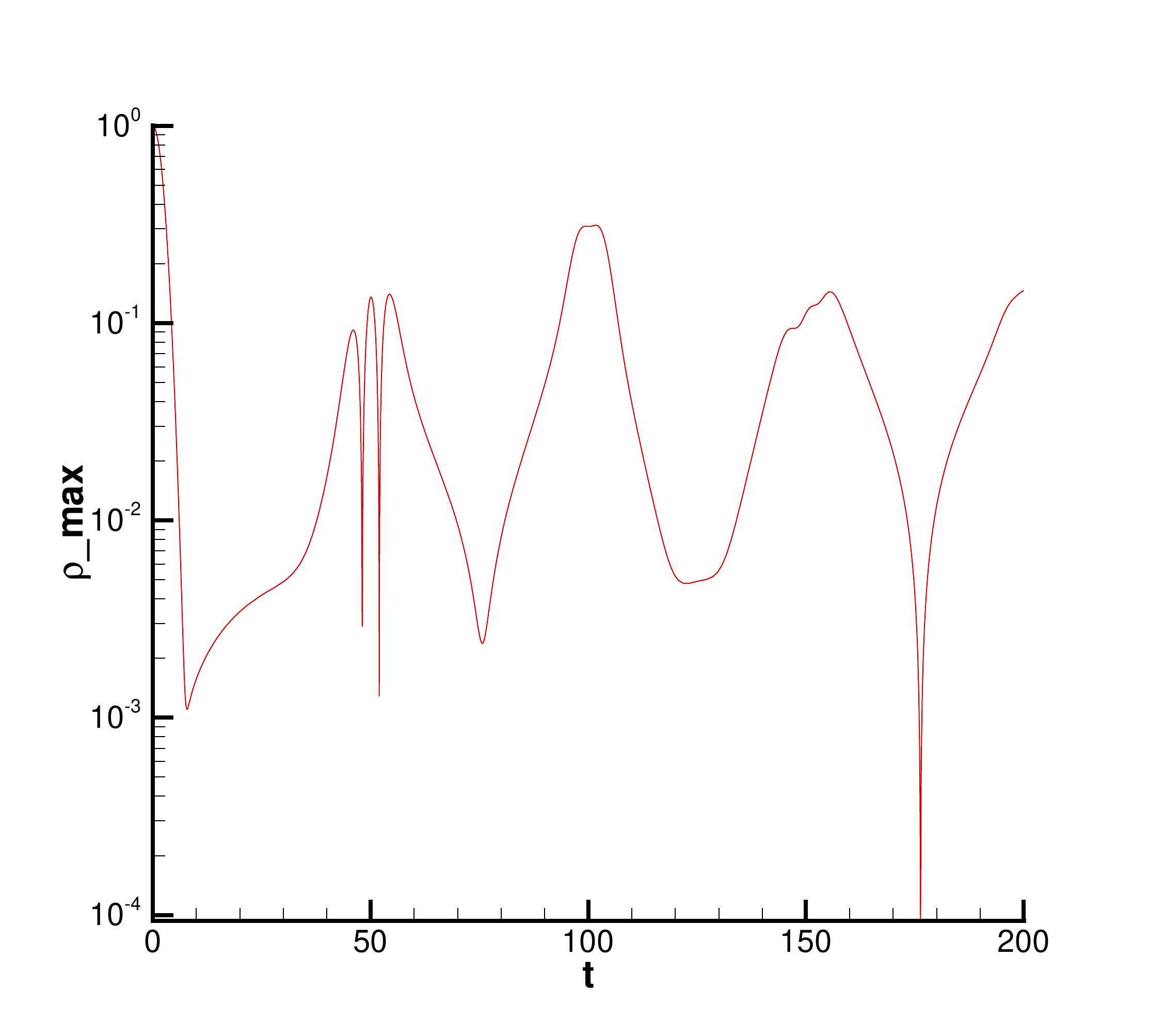}}
    \subfigure[Lorentzian, $\mathbb{P}^1$]{\includegraphics[width=3in,angle=0]{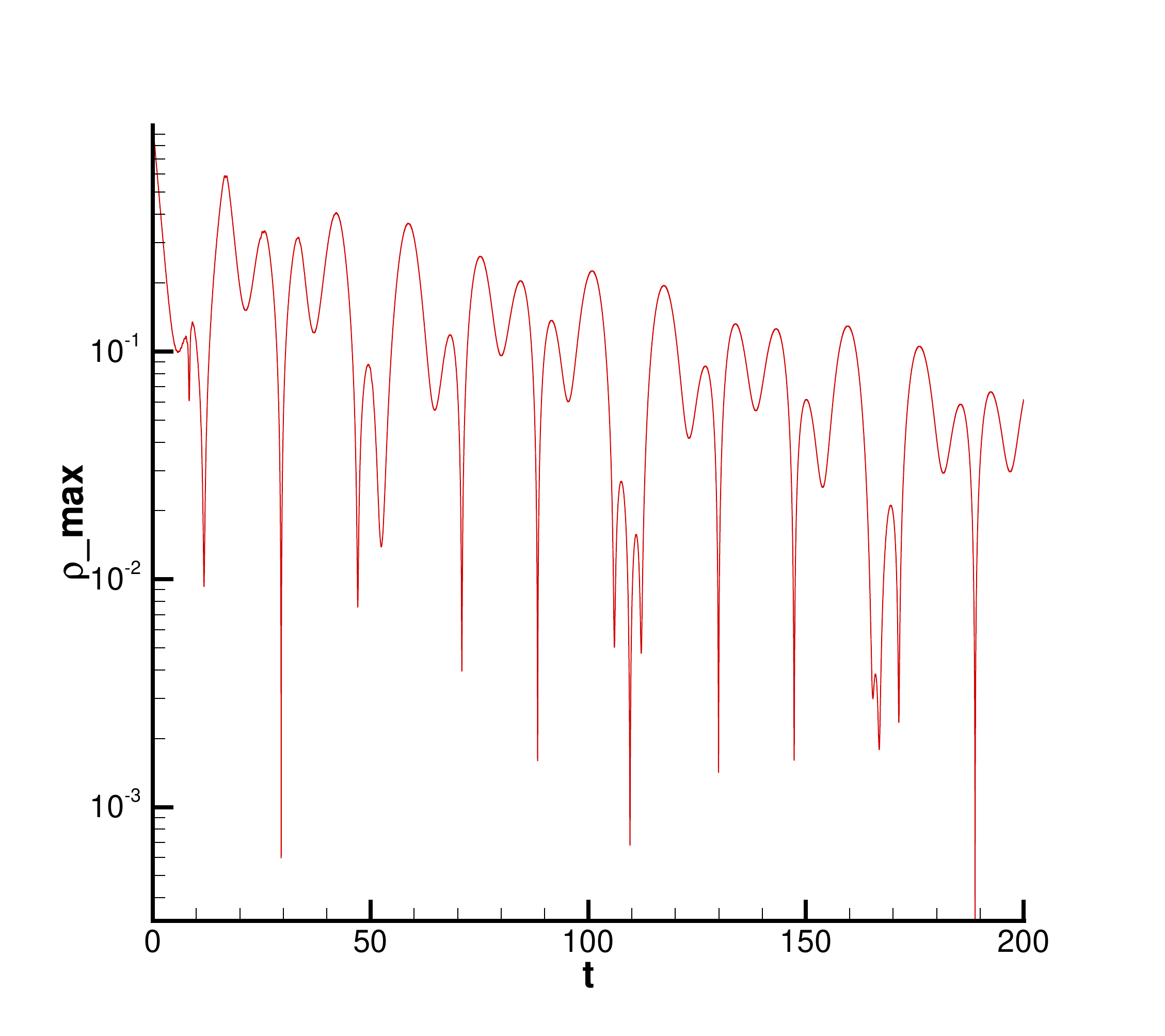}}
    \subfigure[Maxwellian, $\mathbb{P}^2$]{\includegraphics[width=3in,angle=0]{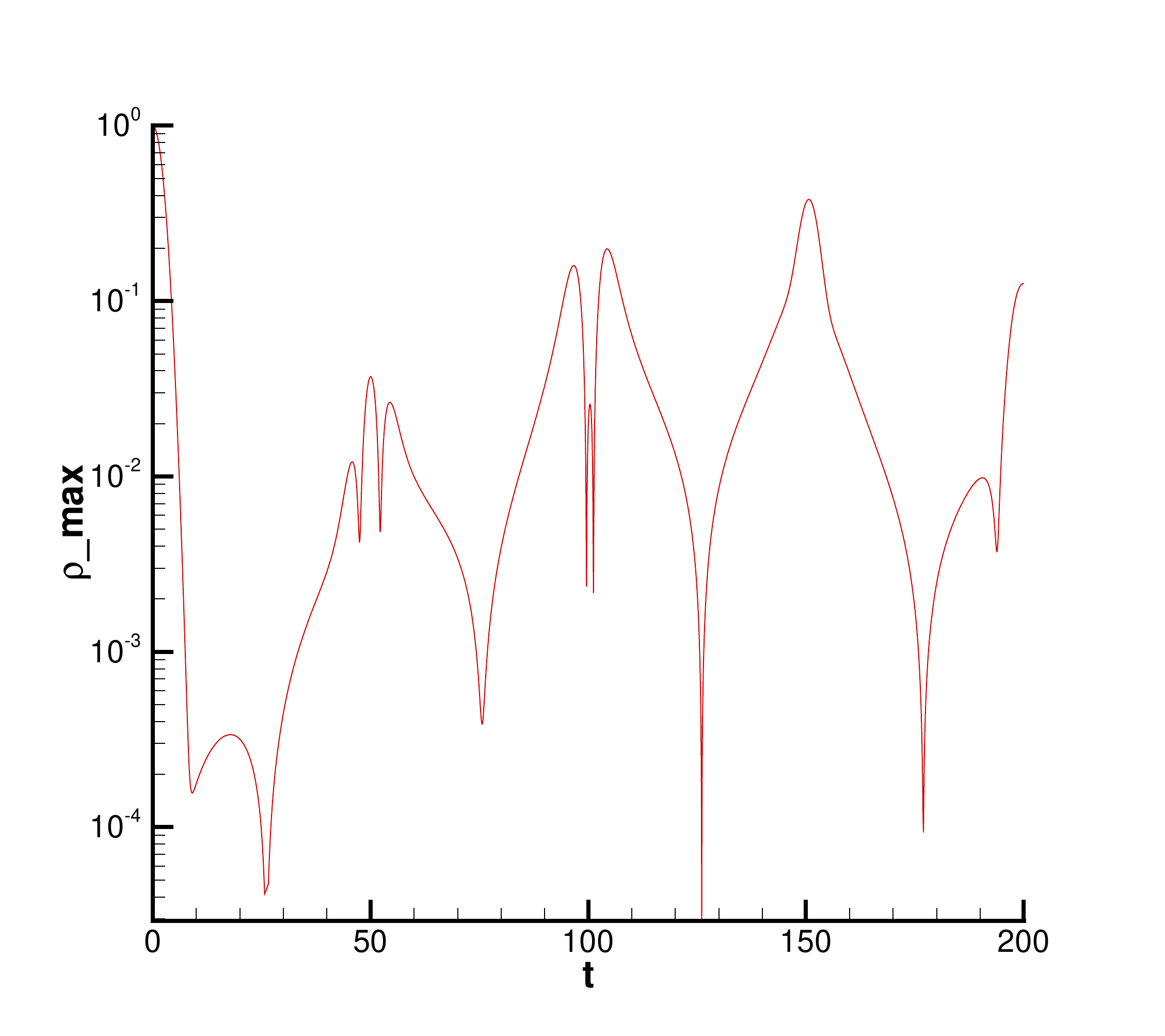}}
    \subfigure[Lorentzian, $\mathbb{P}^2$]{\includegraphics[width=3in,angle=0]{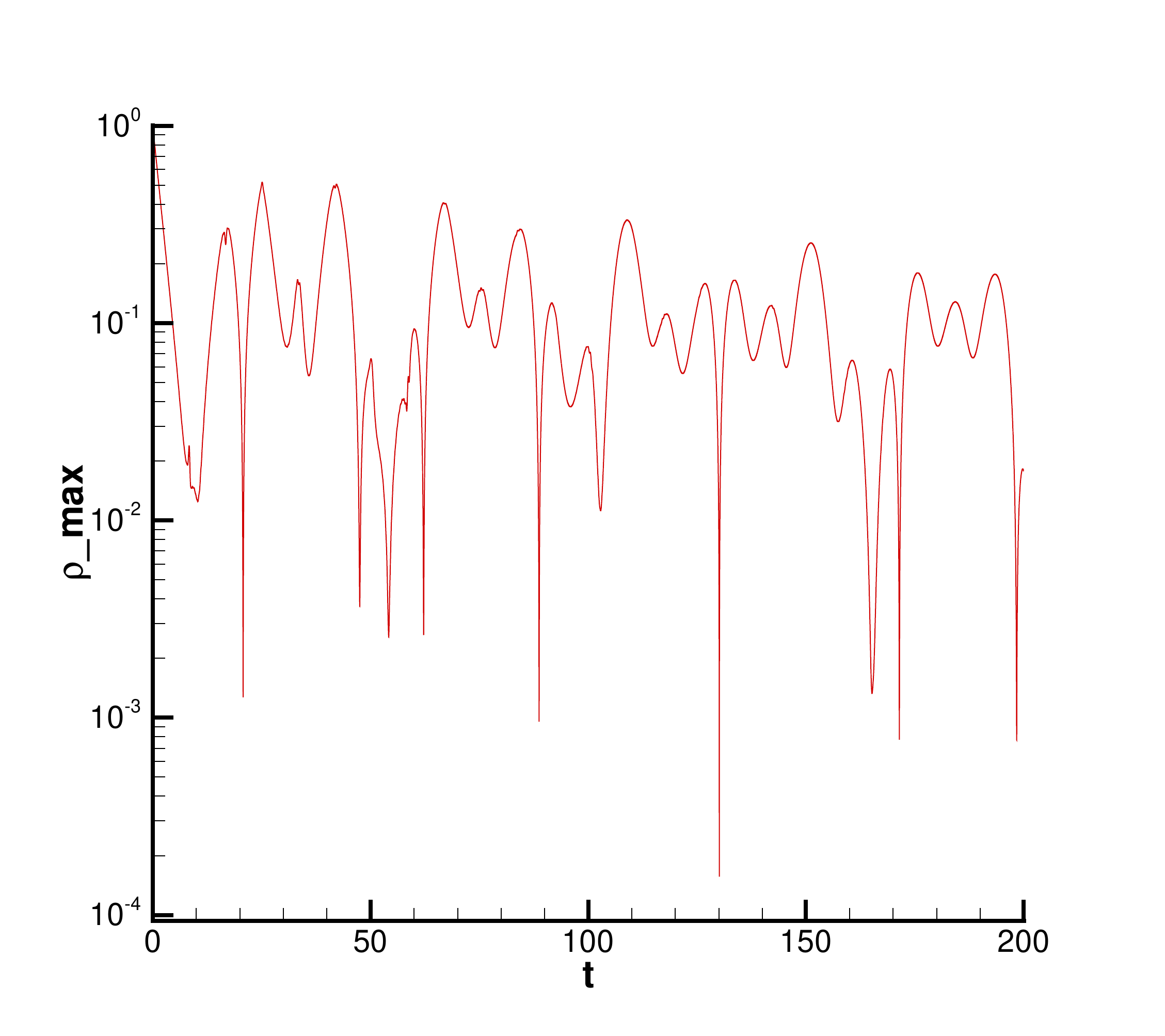}}
  \end{center}
  \caption{Computations of the advection equation for the  polynomial spaces  $\mathbb{P}^1$ and $\mathbb{P}^2$ showing  local maxima of the  density $\rho_{max}$ as a function of time.  The mesh is $40 \times 40$ with    $\triangle x=\pi/10$.  For the Maxwellian equilibrium $\triangle v=1/4$, while for the Lorentzian equilibrium $\triangle v=3/2$.}
  \label{recurp}
\end{figure}

We close  this section with a few comments.  An analysis  similar to that of this section for   $\mathbb{P}^1$  elements yields a $3\times3$ matrix; however, this basis does not yield the  nice form possessed by $\mathbb{Q}^1$ because of the loss of the tensor structure.  Figure \ref{recurp}   shows the temporal behavior of $\rho_{max}$ using the $\mathbb{P}^l$ elements. Observe that,  although the local maxima still are located near  $T_R=\frac{2 \pi}{k \triangle v}$, there appear   to be  several small local maxima instead of one main  maximum, and overall  the long time dissipation seems to be stronger than that for $\mathbb{Q}^l$. We also noted that the $\mathbb{P}^2$ basis follow the exact solution longer than $\mathbb{P}^l$, but shorter than $\mathbb{Q}^2$ cases, because for $\mathbb{P}^2$ elements, the minimum value that the solution achieves before it  deviates from the exact solution is on the order of $10^{-4}$.   In summary, we conclude that increasing the polynomial order does not change $T_R$ by much. However, higher order accuracy seems to improve the time that the numerical solution can follow the exact solution.  For $\mathbb{Q}^l$ elements,  the amplification matrix can be written as a tensor product of two small matrices, and this made  possible our direct analysis for the recurrence time. For $\mathbb{P}^l$ elements, we lose this  tensor structure, and the solution is  more dissipative.  

Finally we remark that since the linearized equation involves an operator $E f'_{eq}(v)$, where the electric field depends on the distribution function $f$ on all cells, it is not trivial to generalize the analysis to the LVP system.  However, it was proven  in \cite{Morrison_92, Morrison_00} that there exists a generalization
  of the Hilbert transform that maps the solution of the advection equation to the solution of this LVP system, so we expect similar type of recurrence behavior for the LVP system, and this is verified by  numerical calculation in Section \ref{linearVP}.

\section{Vlasov numerical results}
\label{numerical}

Now we turn to some numerical tests of  our method  for both the VP and LVP systems.  For the LVP system we consider the standard tests of linear and nonlinear Landau damping,  which have been studied in many references in the contexts of various  numerical techniques  since \cite{chengknorr_76} (see \cite{Heath} for an extended list),  but we also consider a test that heretofore does not appear to have been done, viz., we monitor the linearized energy that is conserved  by the LVP system \cite{KO58,Morrison_92,Morrison_00}.   Similarly, for  the nonlinear VP system we consider the standard tests of nonlinear Landau damping and a symmetric version of the two-stream instability (also see \cite{Heath} for references).  In addition, for the VP system  we consider an example  that is initialized by a driving electric field, resulting in a dynamically accessible initial condition as described in \cite{Morrison_89, Morrison_90,Morrison_92}, which has been observed to approach nonli
 near BGK \cite{BGK} states that have been termed KEEN waves  in
Ref.~\cite{Afeyan_03, Johnston_09}   (see also \cite{Heath_thesis, Valentini}).

\subsection{Linearized VP system}
\label{linearVP}

Associated with the LVP system of (\ref{eq: linearvp}) is the  well-known plasma dispersion function \cite{fried},
\begin{equation}
\varepsilon(k,\omega)=1-\frac{1}{k^2} \int_{-\infty}^{+\infty} \frac{f_{eq}'(v)}{v-\omega/k} \, dv,
\label{dispersion}
\end{equation}
which (with the appropriate choice of contour) will be used to benchmark the accuracy of the Landau damping rate and oscillation frequency obtained from our   DG solver  with choices for the various polynomial spaces.
The LVP system conserves not only the total charge  and momentum, but also the linear energy \cite{KO58,Morrison_92,Morrison_00}, which is defined as
\begin{equation}
H_L= -\frac{1}{2} \int_{\Omega} \frac{v f^2}{f'_{eq}}\,  dxdv + \frac{1}{2} \int_{\Omega_x} E^2\,  dx\, .
\label{hl}
\end{equation}
As noted above,  we monitor this quantity and check for its conservation.  In addition,  we monitor the shift of energy to the first term of (\ref{hl}) as the second decays in time in accordance with  Landau damping, consistent with the discussion of \cite{Morrison_92}.

\subsubsection*{Linear Landau damping}

For this  classical test problem,   we choose the usual  initial condition  $f_0(x,v)=A\cos(kx)f_M(v)$, with $A=0.01$ and  $k=0.5$.  For the Maxwellian distribution function the dispersion relation becomes
\[
\varepsilon(k,\omega)=1+\frac{1}{k^2} \left\{ 1+ \frac{\omega}{\sqrt{2}k} Z\left(\frac{\omega}{\sqrt{2}k}\right) \right\},
\]
where the plasma $Z$-function is defined as
$$
Z(z)=\frac{1}{\sqrt{\pi}} \int_{-\infty}^{\infty} e^{-t^2} \frac{dt}{t-z} = 2 i e^{-z^2} \int_{-\infty}^{iz} e^{-t^2} dt.
$$
{}From this relation,   the predicted damping rate is computed to be 0.153359 and the predicted oscillation frequency to be 1.41566.

\begin{figure}[htb]
  \begin{center}
   \subfigure[ $\mathbb{P}^0, 40 \times 40$ mesh]{\includegraphics[width=2in,angle=0]{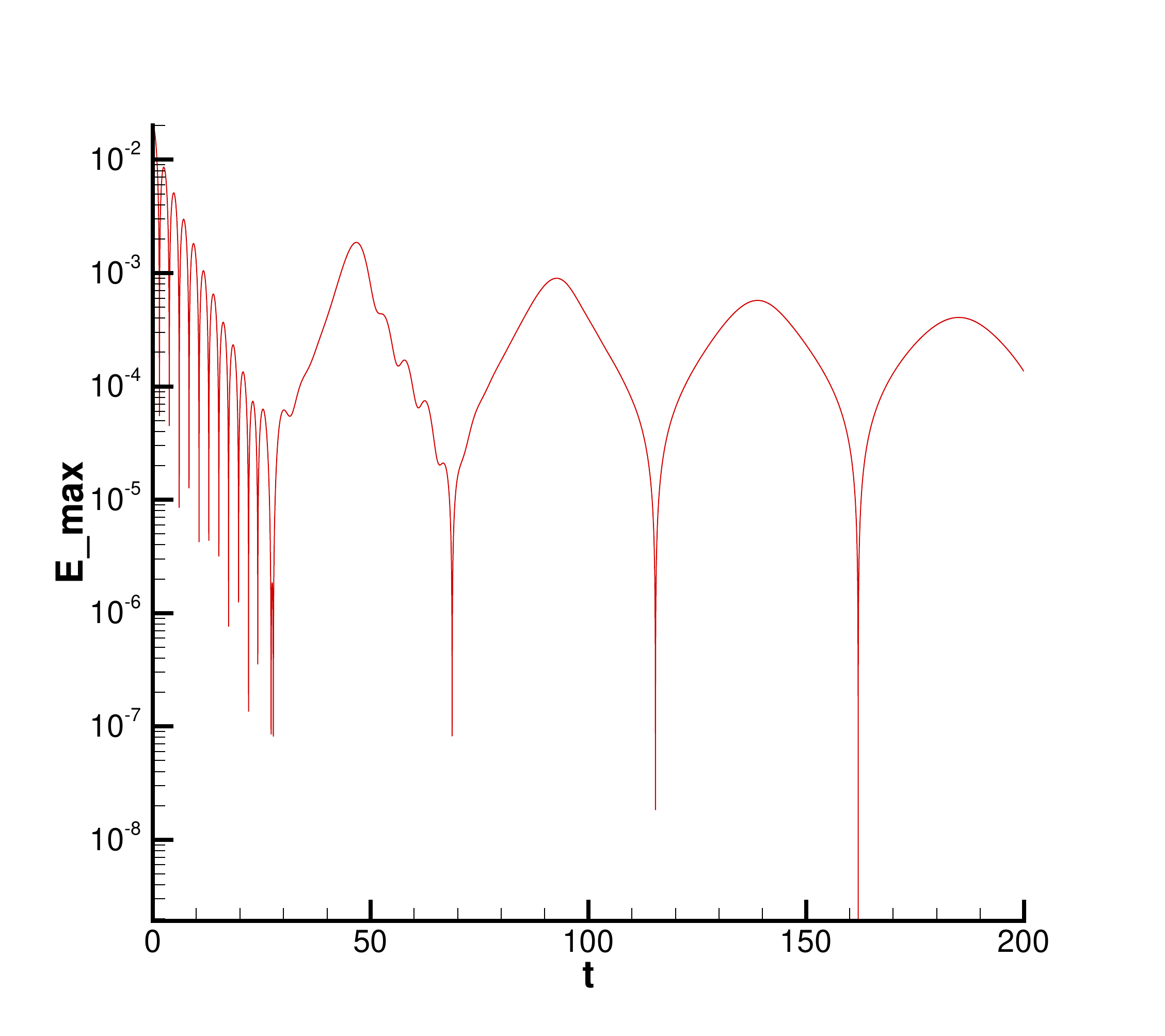}}
      \subfigure[ $\mathbb{P}^1, 40 \times 40$ mesh]{\includegraphics[width=2in,angle=0]{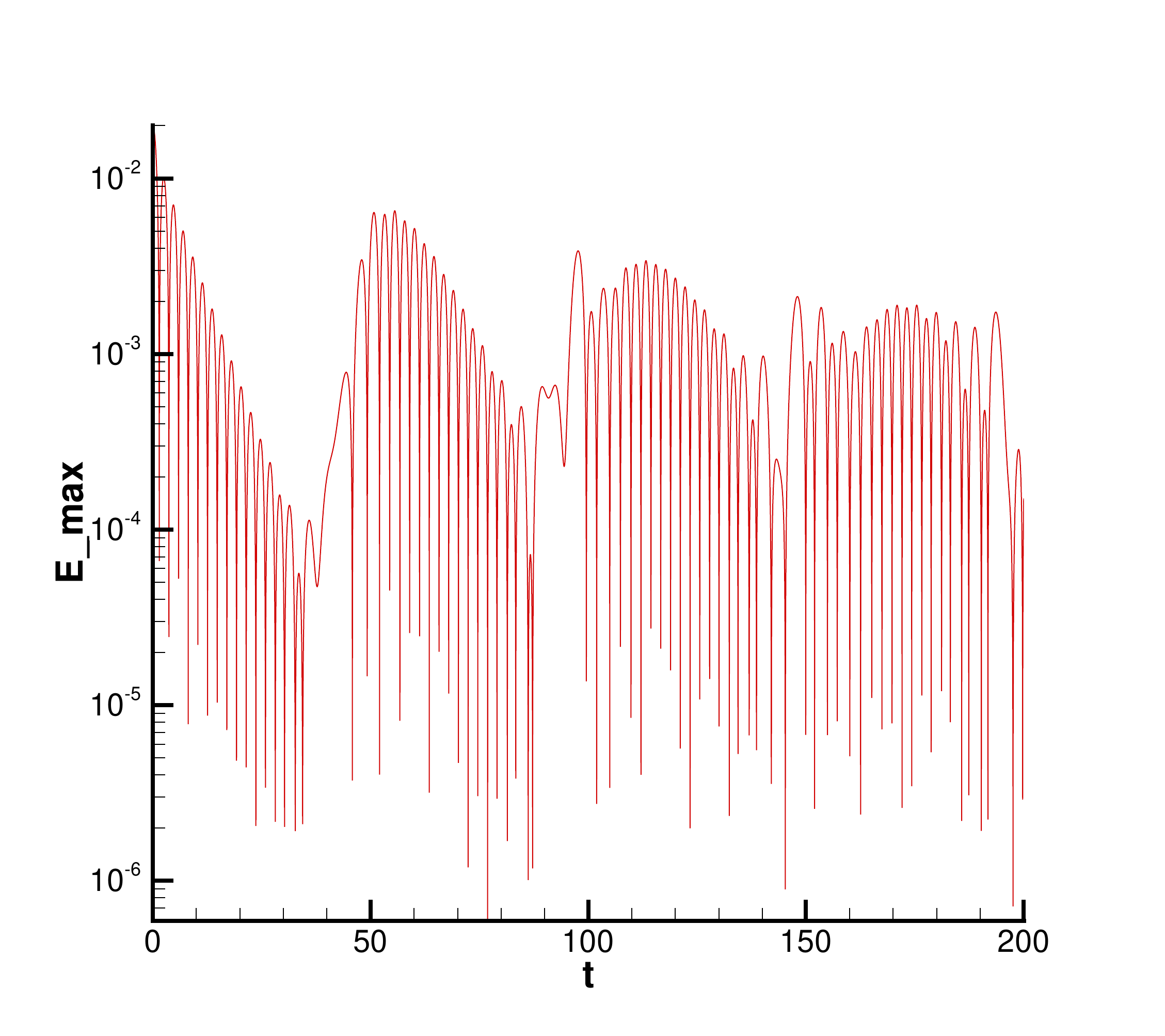}}
      \subfigure[ $\mathbb{P}^2, 40 \times 40$ mesh]{\includegraphics[width=2in,angle=0]{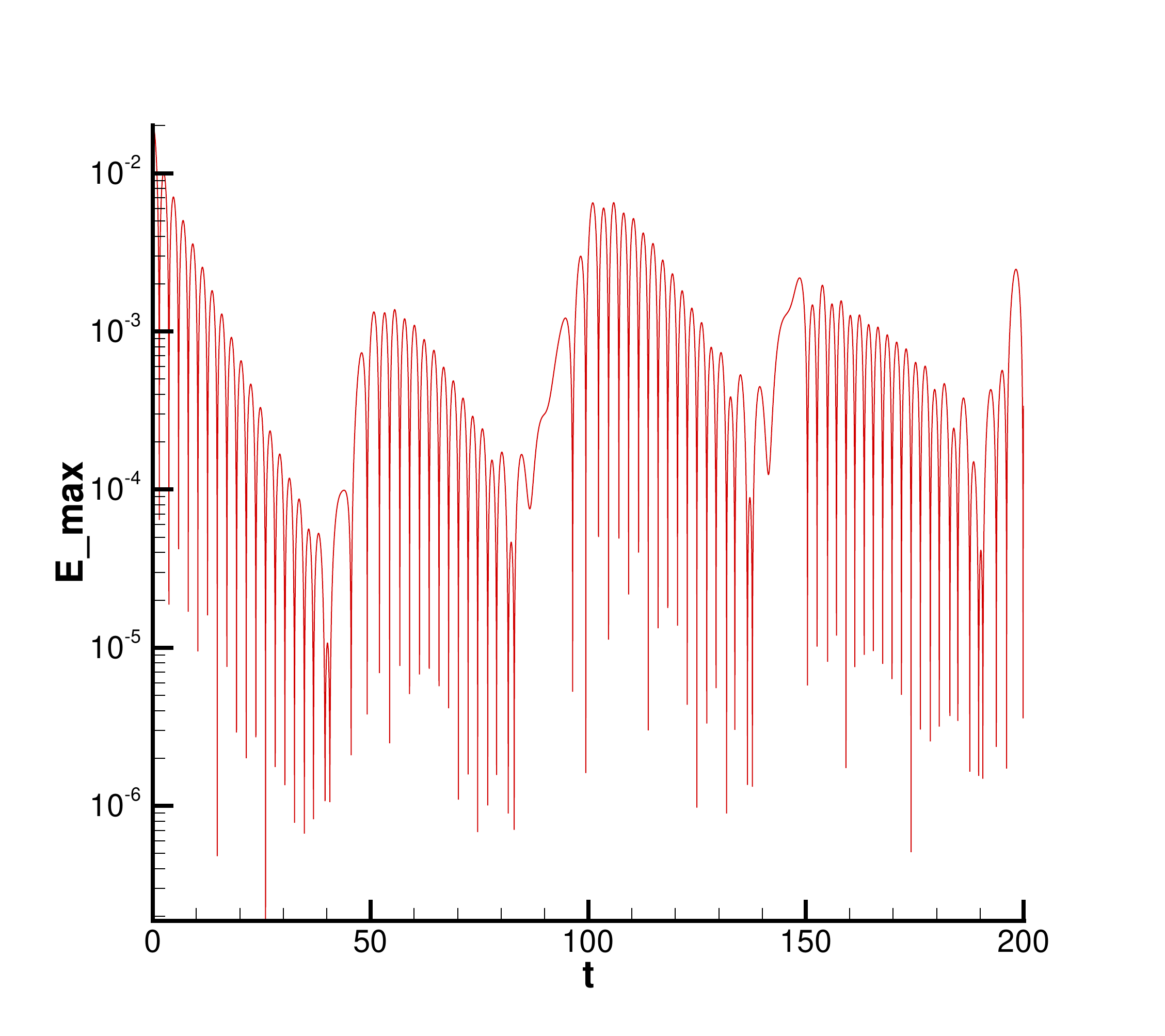}}
      \subfigure[ $\mathbb{P}^0, 80 \times 80$ mesh]{\includegraphics[width=2in,angle=0]{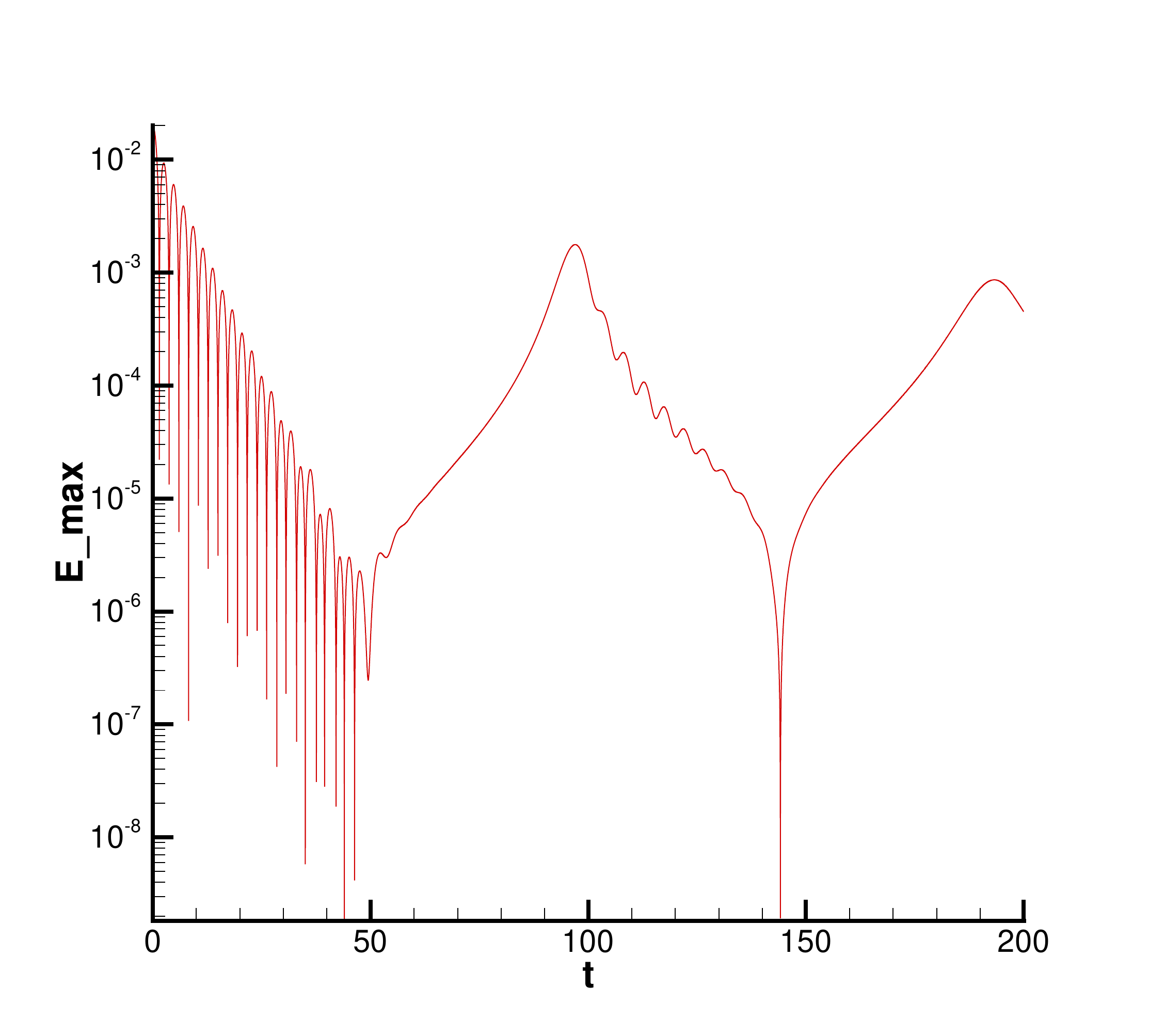}}
    \subfigure[ $\mathbb{P}^1, 80 \times 80$ mesh]{\includegraphics[width=2in,angle=0]{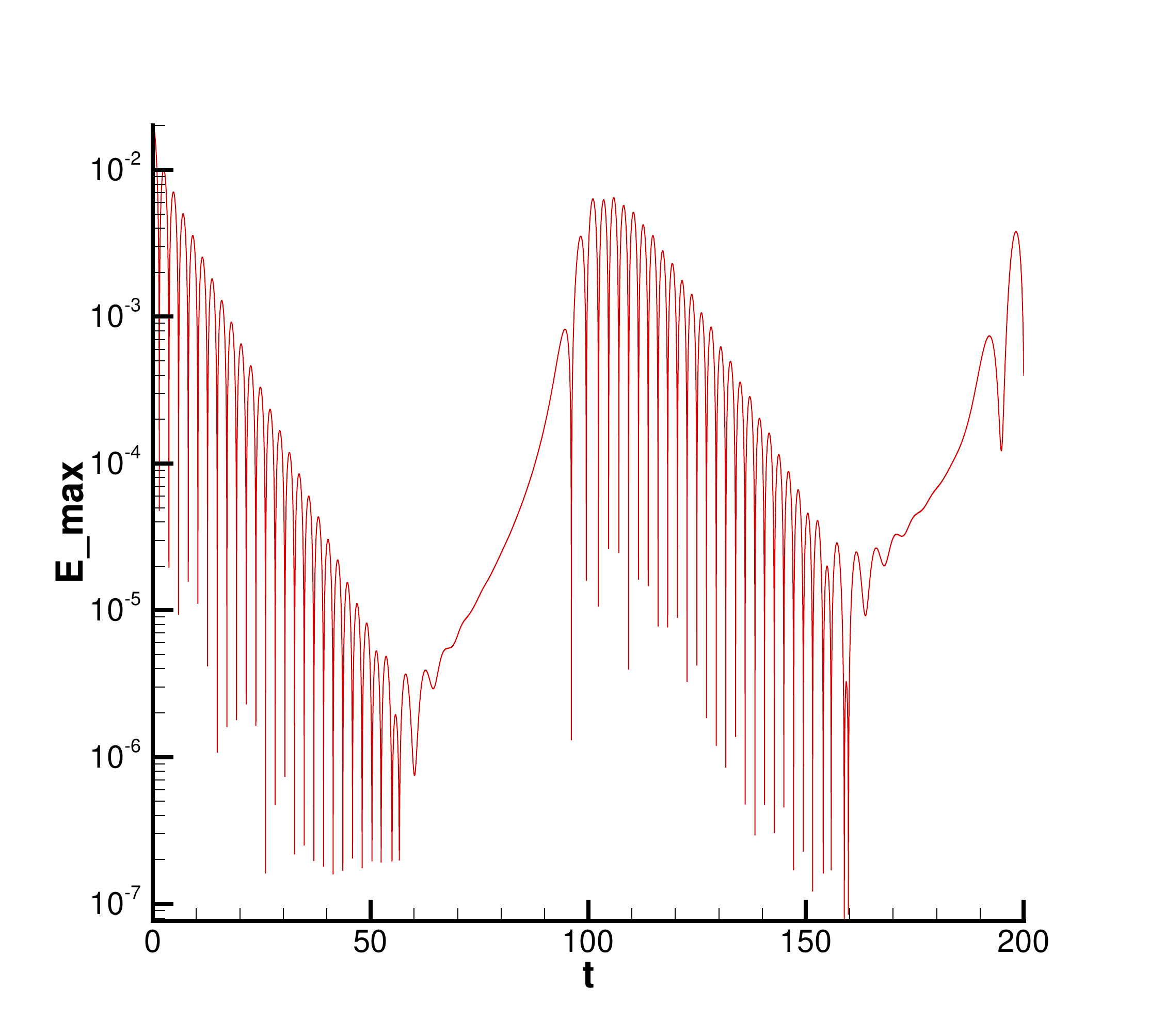}}
   \subfigure[ $\mathbb{P}^2, 80 \times 80$ mesh]{\includegraphics[width=2in,angle=0]{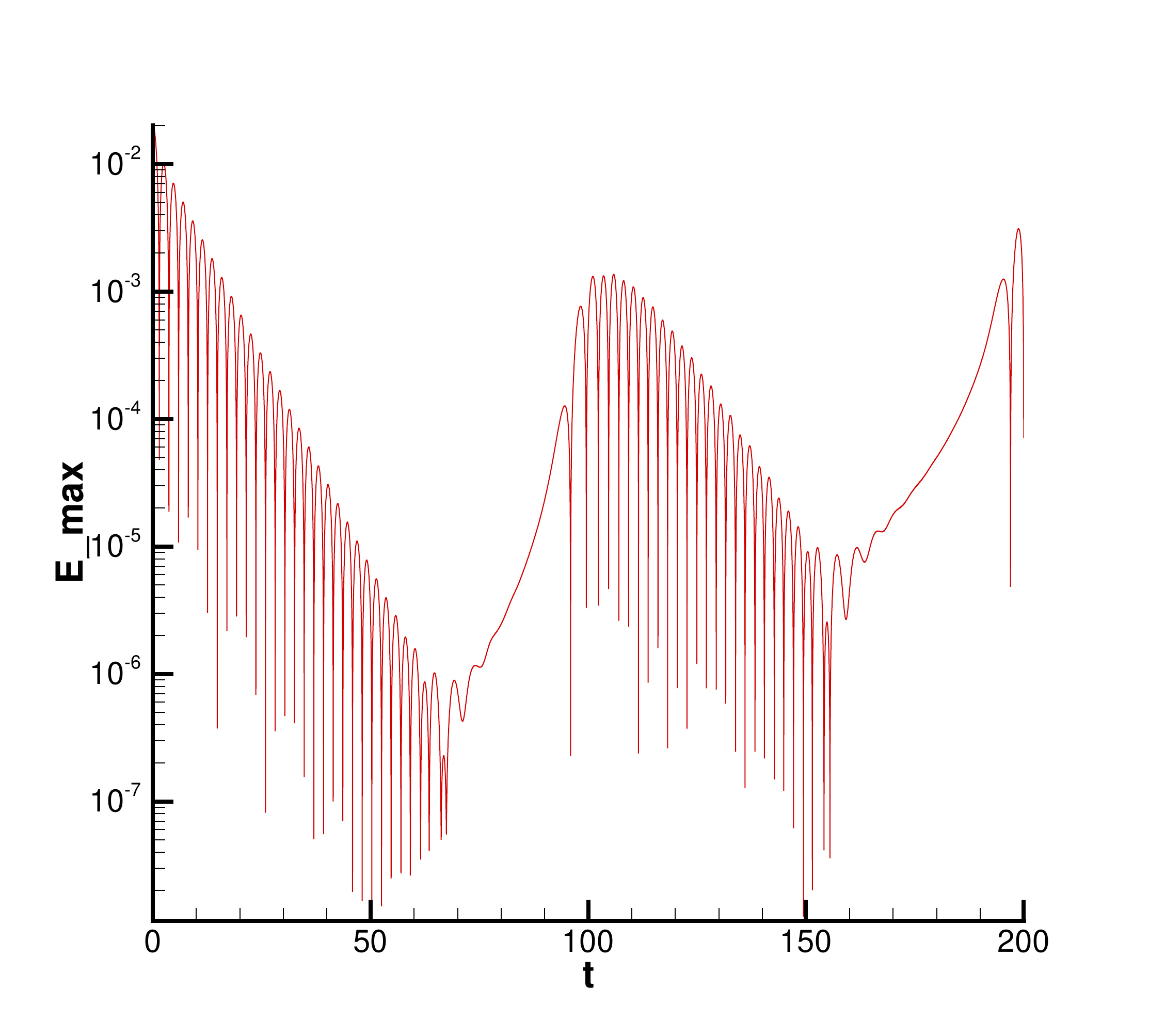}}
         \subfigure[ $\mathbb{Q}^1, 40 \times 40$ mesh]{\includegraphics[width=2.in,angle=0]{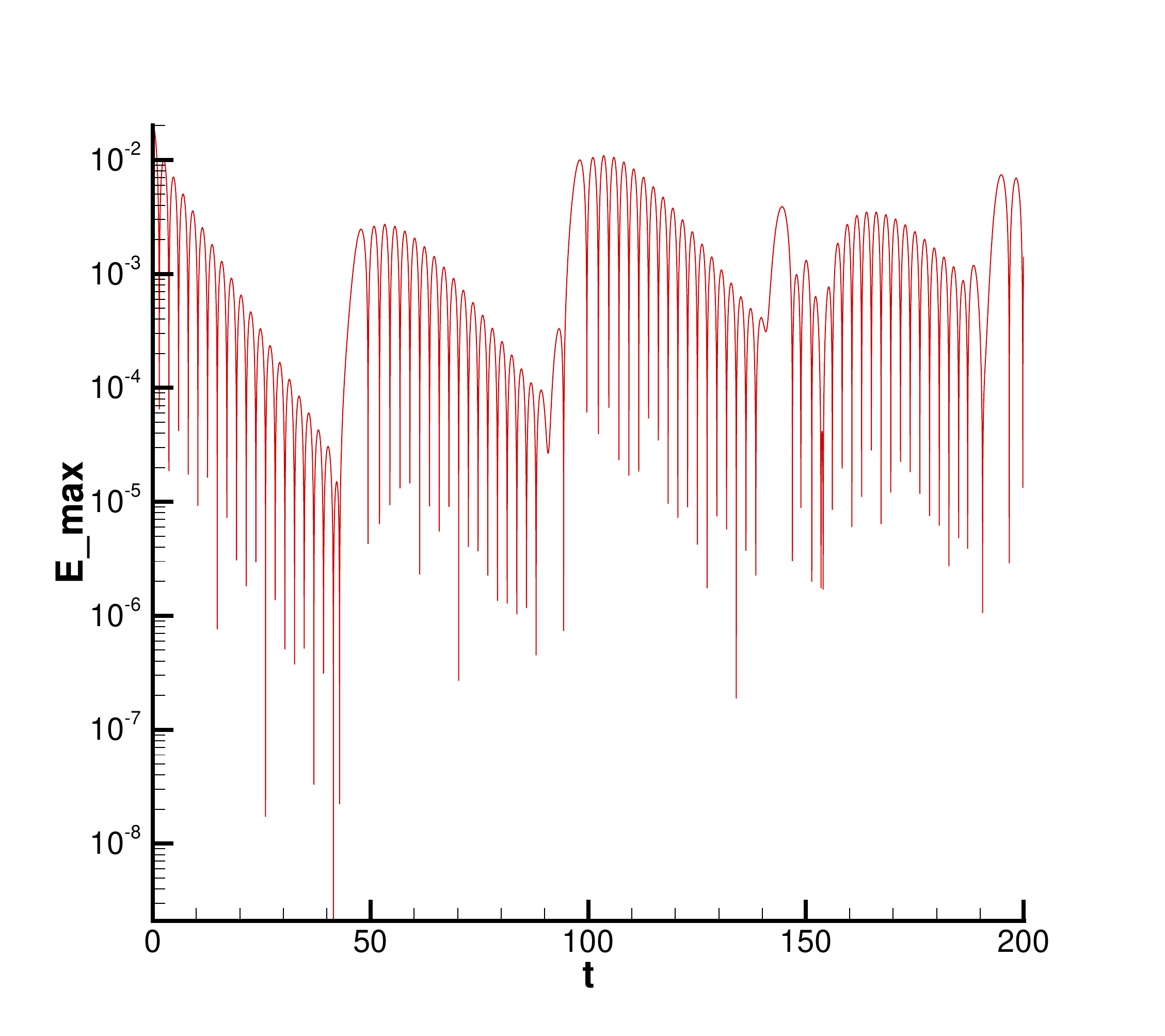}}
      \subfigure[ $\mathbb{Q}^2, 40 \times 40$ mesh]{\includegraphics[width=2.in,angle=0]{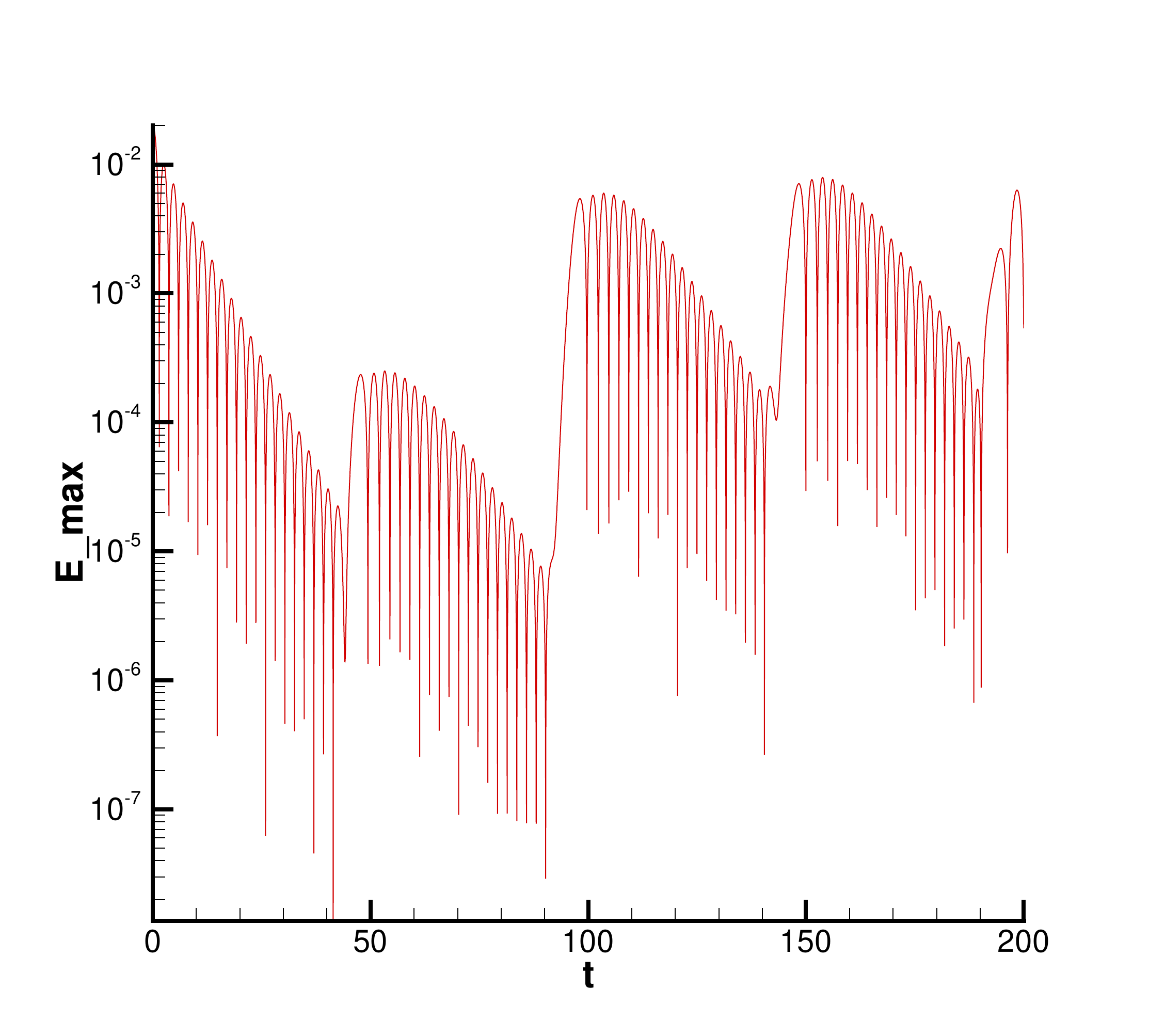}}\\
         \subfigure[ $\mathbb{Q}^1, 80 \times 80$ mesh]{\includegraphics[width=2.in,angle=0]{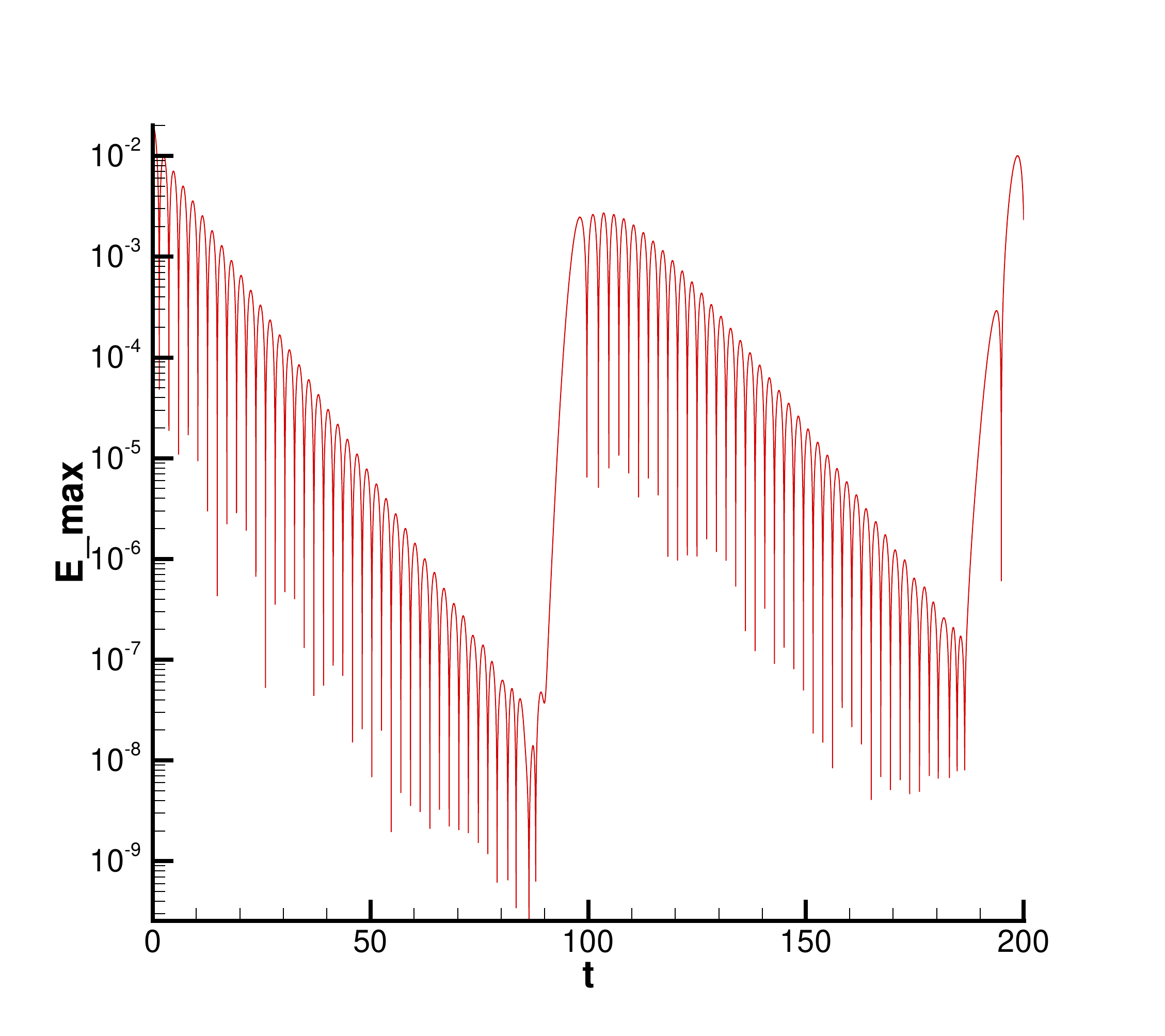}}
   \subfigure[ $\mathbb{Q}^2, 80 \times 80$ mesh]{\includegraphics[width=2.in,angle=0]{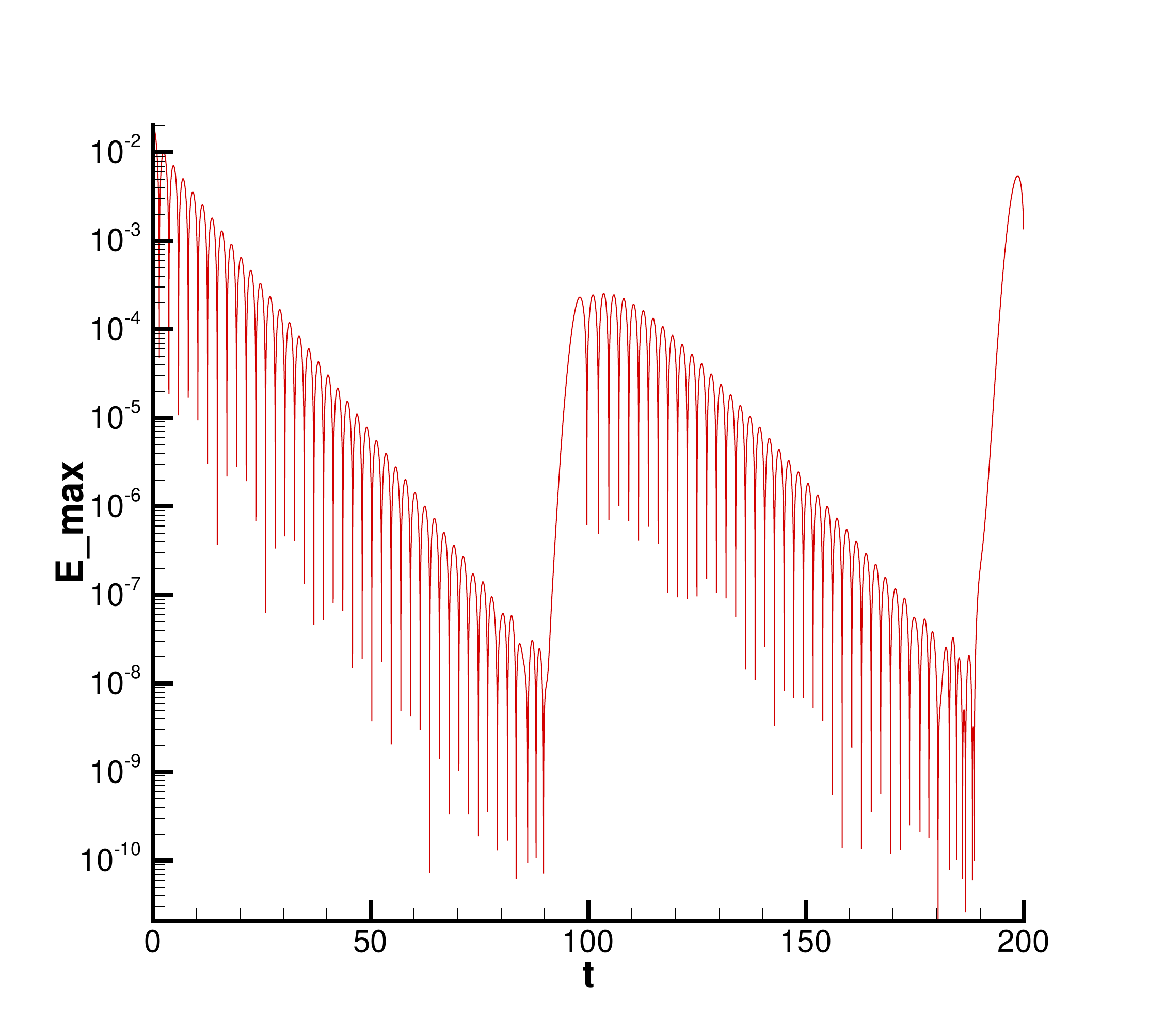}}
  \end{center}
  \caption{Depiction of linear Landau damping showing recurrence in the maxima  of the electric field, $E_{max}$, as a function of time for various polynomial spaces.}
\label{linearlog}
\end{figure}

In Figures \ref{linearlog}, we plot the evolution of the maximum of the electric field $E_{max}$ using various polynomial spaces. In Table \ref{tdamp}, we compare the theoretical and numerical values of damping rate and frequency as a measurement of accuracy.  We see that refining the mesh always gives  better approximations. The piecewise constant polynomials $\mathbb{P}^0$ give much larger error compared to higher order polynomials. While the difference between the $\mathbb{P}^l$ and $\mathbb{Q}^l$ spaces is not significant.  Observe from Figure \ref{linearlog}  how  similar the recurrence behavior is for this LVP problem to that of  the advection equation.

\begin{table} [htb]
\begin{center}
\caption {The damping rate and frequency for linear Landau damping. The numerical values are computed using the fourth to the tenth peak and the predicted value is obtained from the plasma dispersion function (\ref{dispersion})  with a Maxwellian equilibrium.}
\bigskip
\begin{tabular}{|c|c|c|c|c|c|c|c|}
\hline
& Predicted value & Mesh &$\mathbb{P}^0$ & $\mathbb{P}^1$ &$\mathbb{P}^2$ & $\mathbb{Q}^1$ & $\mathbb{Q}^2$\\
\hline
\multirow{2}{*}{Damping rate} & \multirow{2}{*}{0.153359} & $40\times 40$ & 0.227489&0.153536&0.153375&0.153425&0.153379 \\
\cline{3-8} & & $80\times 80$ & 0.191702&0.153366& 0.153363&0.15369 & 0.153363\\
\hline
\multirow{2}{*}{Frequency} & \multirow{2}{*}{1.41566} & $40\times 40$ & 1.38249 &1.41643 &1.41643 &1.41643 &1.41643\\
\cline{3-8} & & $80\times 80$ & 1.40056 & 1.41576& 1.41576 & 1.41576 & 1.41576\\
\hline
\end{tabular}
\label{tdamp}
\end{center}
\end{table}

\begin{figure}[htb]
  \begin{center}
   \subfigure[ $40 \times 40$ mesh]{\includegraphics[width=2.5in,angle=0]{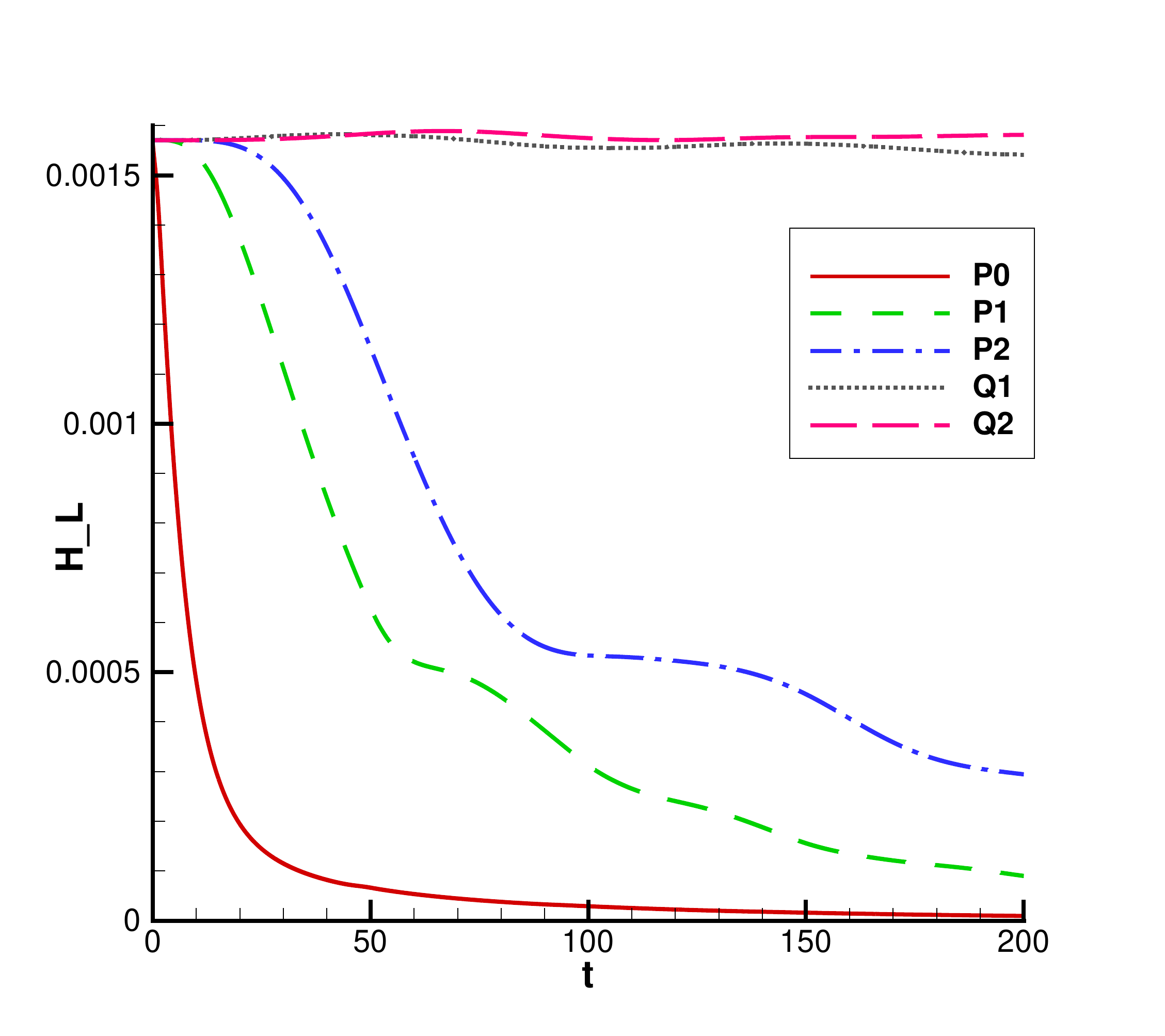}}
    \subfigure[ $80 \times 80$ mesh]{\includegraphics[width=2.5in,angle=0]{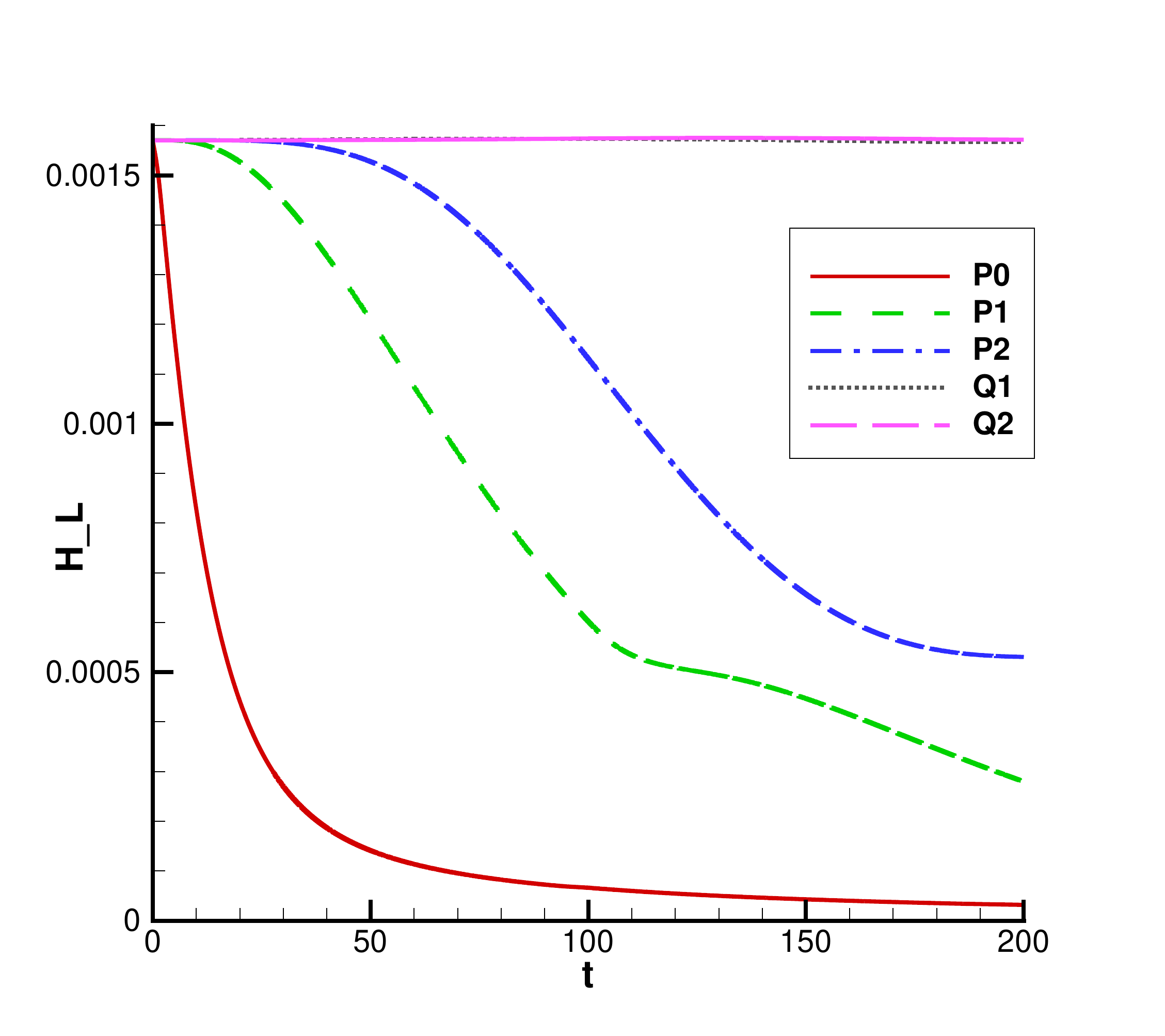}}
  \end{center}
  \caption{Evolution of  the linear energy $H_L$ of (\ref{hl}) as a function of time,  while the Vlasov system undergoes linear Landau damping.  Various polynomial spaces and mesh sizes were used,  as indicated. }
\label{linearene}
\end{figure}

As for conservation properties, the charge and momentum are well conserved as predicted by Propositions 1 and 3.  However, the linear energy $H_L$ demonstrates different behaviors depending on the polynomial spaces. Figure \ref{linearene} shows that  $H_L$ decays significantly for all $\mathbb{P}^l$ spaces even upon mesh refinement. On the other hand, the $\mathbb{Q}^l$ seems to conserve it much better. We note that  $\mathbb{Q}^1$ conserves $H_L$ much better than $\mathbb{P}^2$, although the former is a subspace of the later.

\begin{figure}[htb]
\begin{center}
\includegraphics[width=3.in,angle=0]{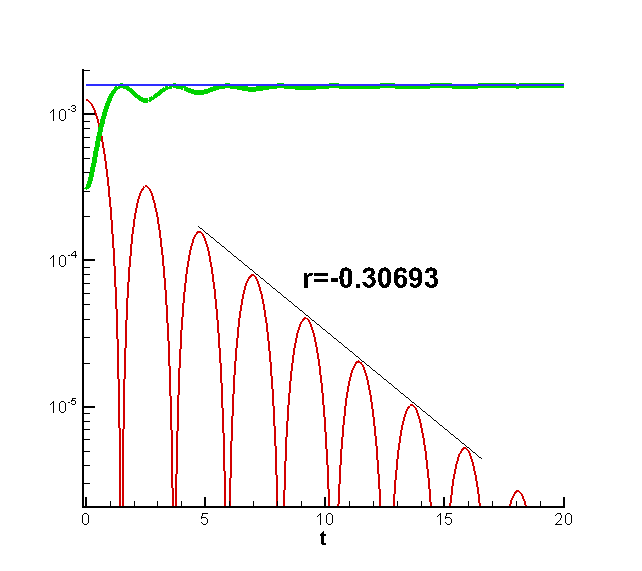}
\end{center}
  \caption{Evolution of  the electrostatic energy (red), linear energy (blue) and the first term in the linear energy (green)  as a function of time,  while the Vlasov system undergoes linear Landau damping.  Here $\mathbb{Q}^2$ was used with a $80 \times 80$ mesh. }
\label{linearele}
\end{figure}

 Also,  note from Figure \ref{linearele}  that the electrostatic energy for both choices of polynomial spaces damps at a rate given by twice the Landau damping rate.  This is to be expected for  the linear theory,  since after integration over space the oscillatory component is removed and $E\sim \exp(-2\gamma t)$.    Therefore, if the  energy is conserved numerically this damped electrostatic energy must be  converted into the relative kinetic energy that is represented by the first term of (\ref{hl}).   Thus, conservation of $H_L$ serves as a  global measure of the ability of an algorithm to resolve fine scales in velocity space.  That this transference must take place for the linear VP system was proven in Section IV of Ref.~\cite{Morrison_92}.

\subsection{Nonlinear VP system}
\label{NLVP}

In this section, we consider the nonlinear VP system.  As noted above, we benchmark the solver against three test cases: the nonlinear Landau damping, two-stream instability,  and  an external drive problem with dynamically accessible initial condition.

The $n$-th Log Fourier mode for the electric field $E(x,t)$ \cite{Heath} is defined as
$$
logF\!M_n(t)=\log_{10} \left(\frac{1}{L} \sqrt{\left|\int_0^L E(x, t) \sin(knx) \, dx \right|^2 +
\left|\int_0^L E(x, t) \cos(knx)\,  dx \right|^2} \right).
$$
We will use this quantity to plot data from our various runs.

\subsubsection*{Nonlinear Landau damping}

\begin{figure}[htb]
  \begin{center}
      \subfigure[ $logF\!M_1$, $\mathbb{P}^2$ no limiter]{\includegraphics[width=2in,angle=0]{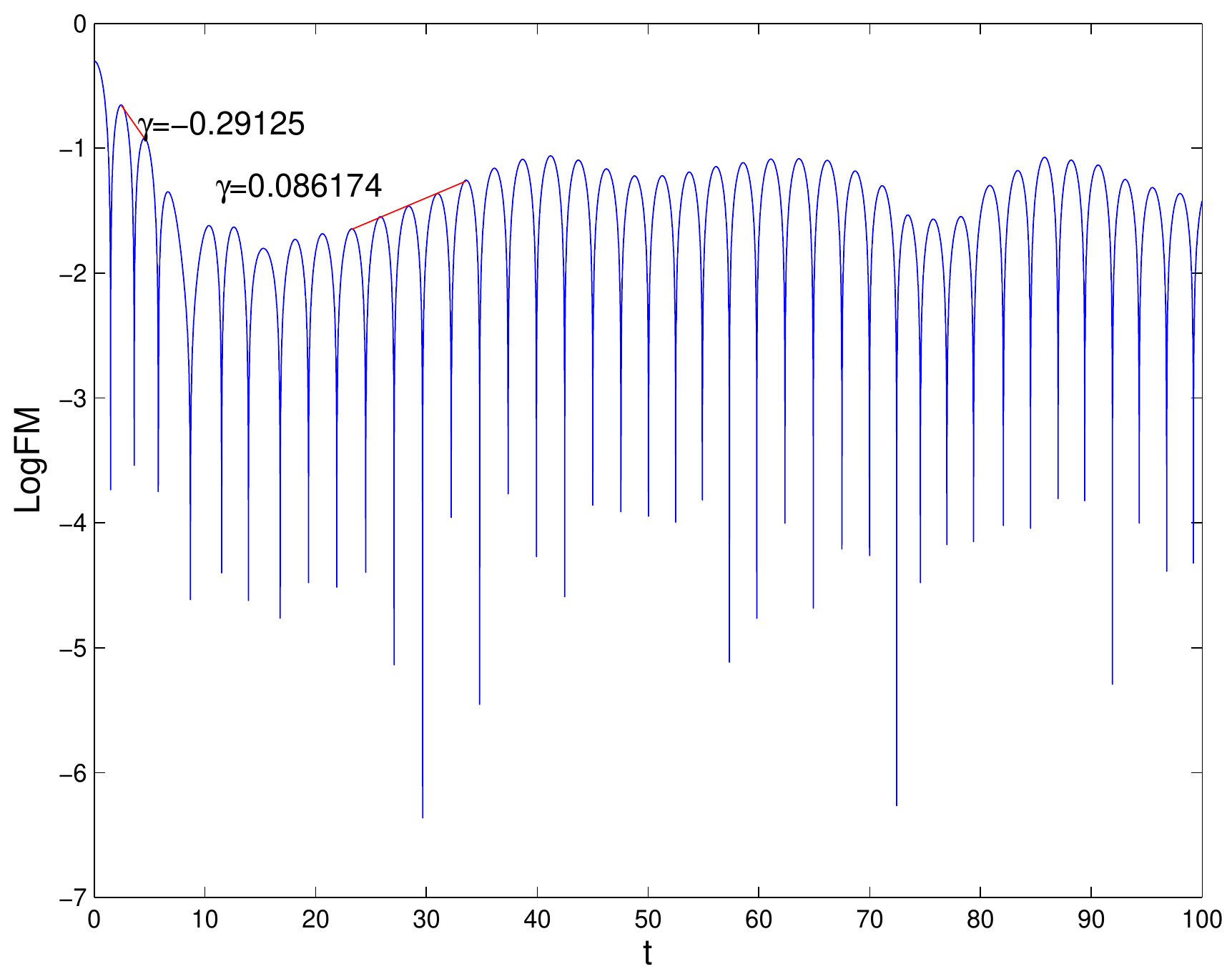}}
               \subfigure[ $logF\!M_1$, $\mathbb{P}^2$ with limiter]{\includegraphics[width=2in,angle=0]{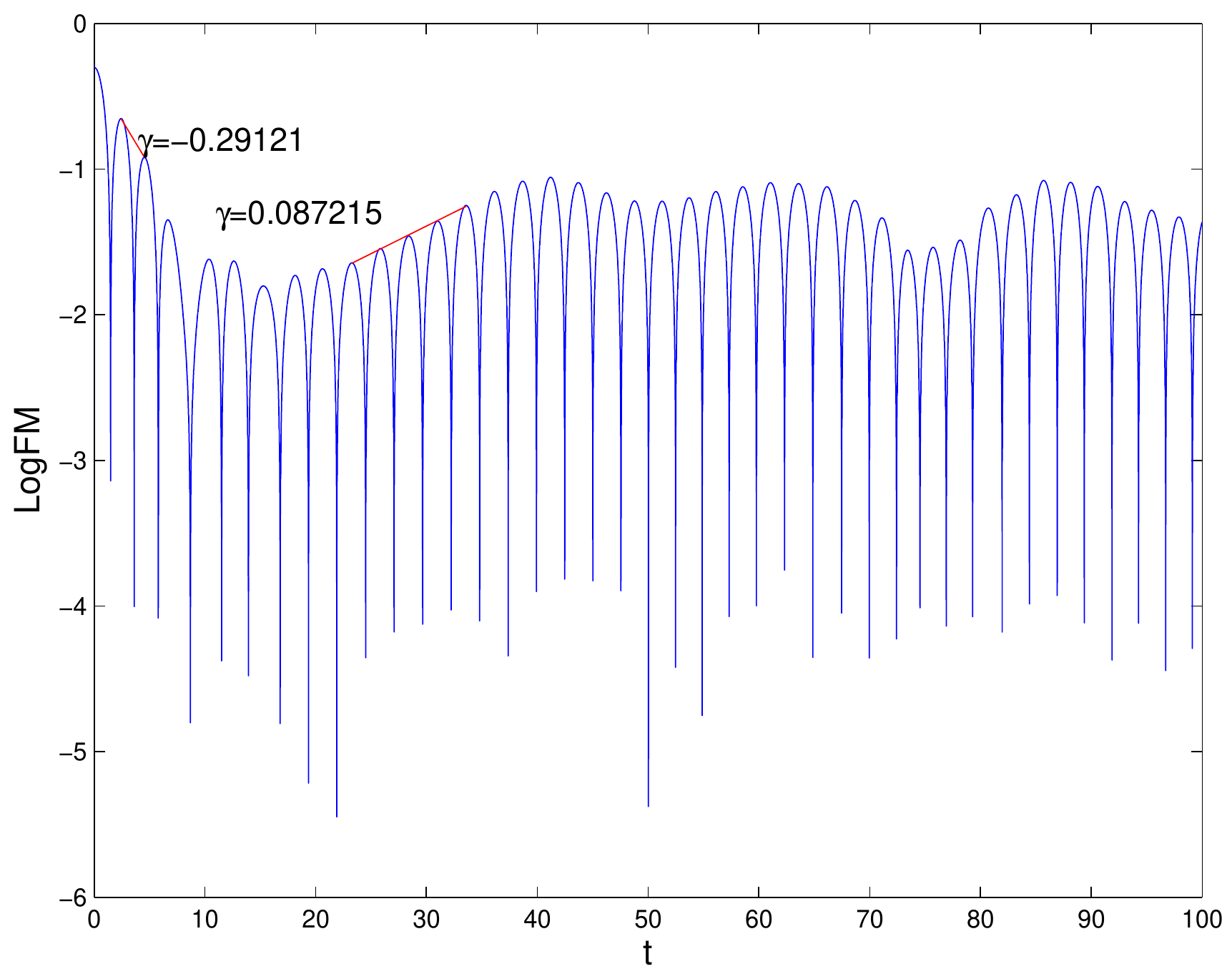}}
                        \subfigure[ $logF\!M_1$, $\mathbb{Q}^2$ no limiter ]{\includegraphics[width=2in,angle=0]{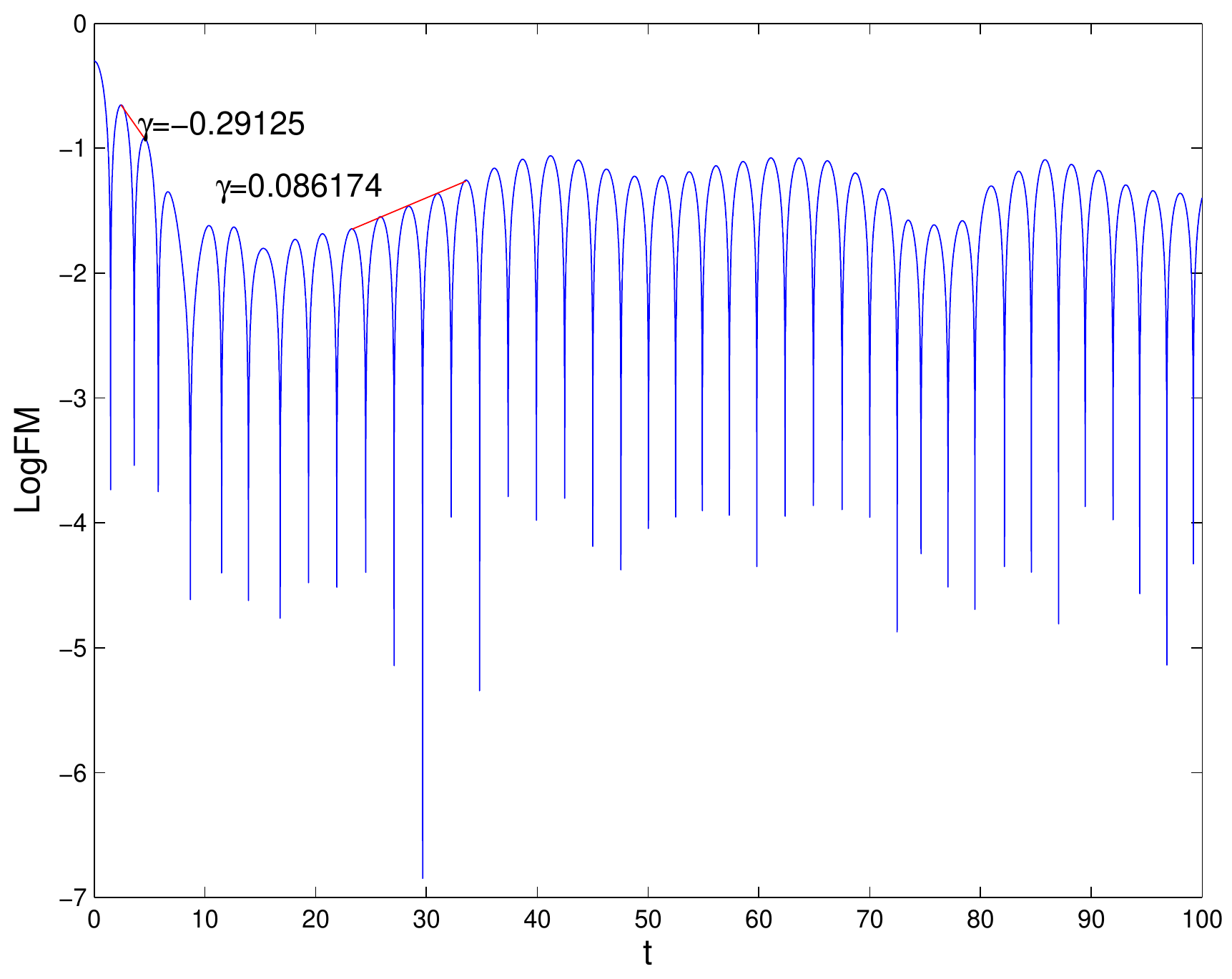}}
   \subfigure[ $logF\!M_2$, $\mathbb{P}^2$ no limiter]{\includegraphics[width=2in,angle=0]{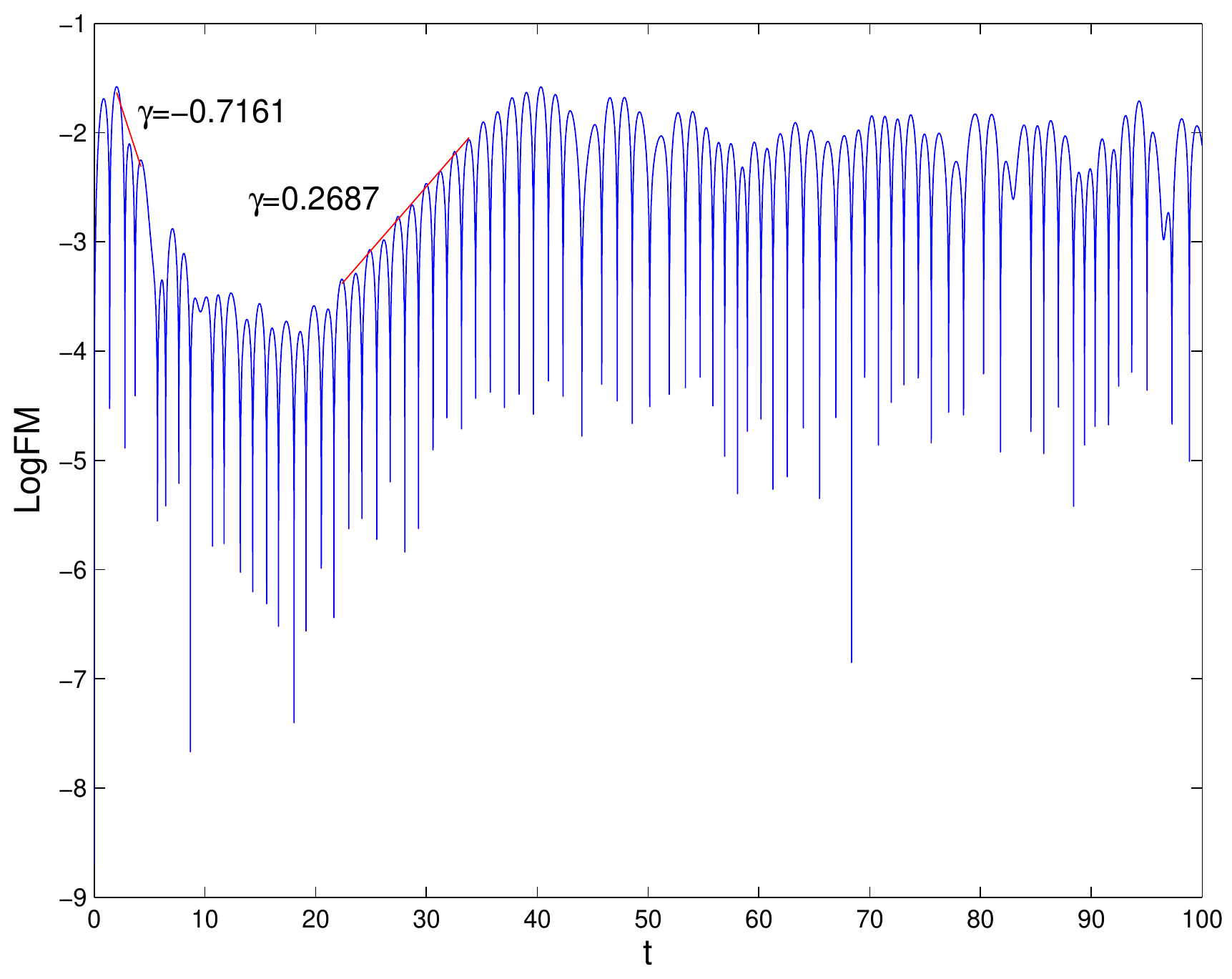}}
      \subfigure[ $logF\!M_2$, $\mathbb{P}^2$ with limiter]{\includegraphics[width=2in,angle=0]{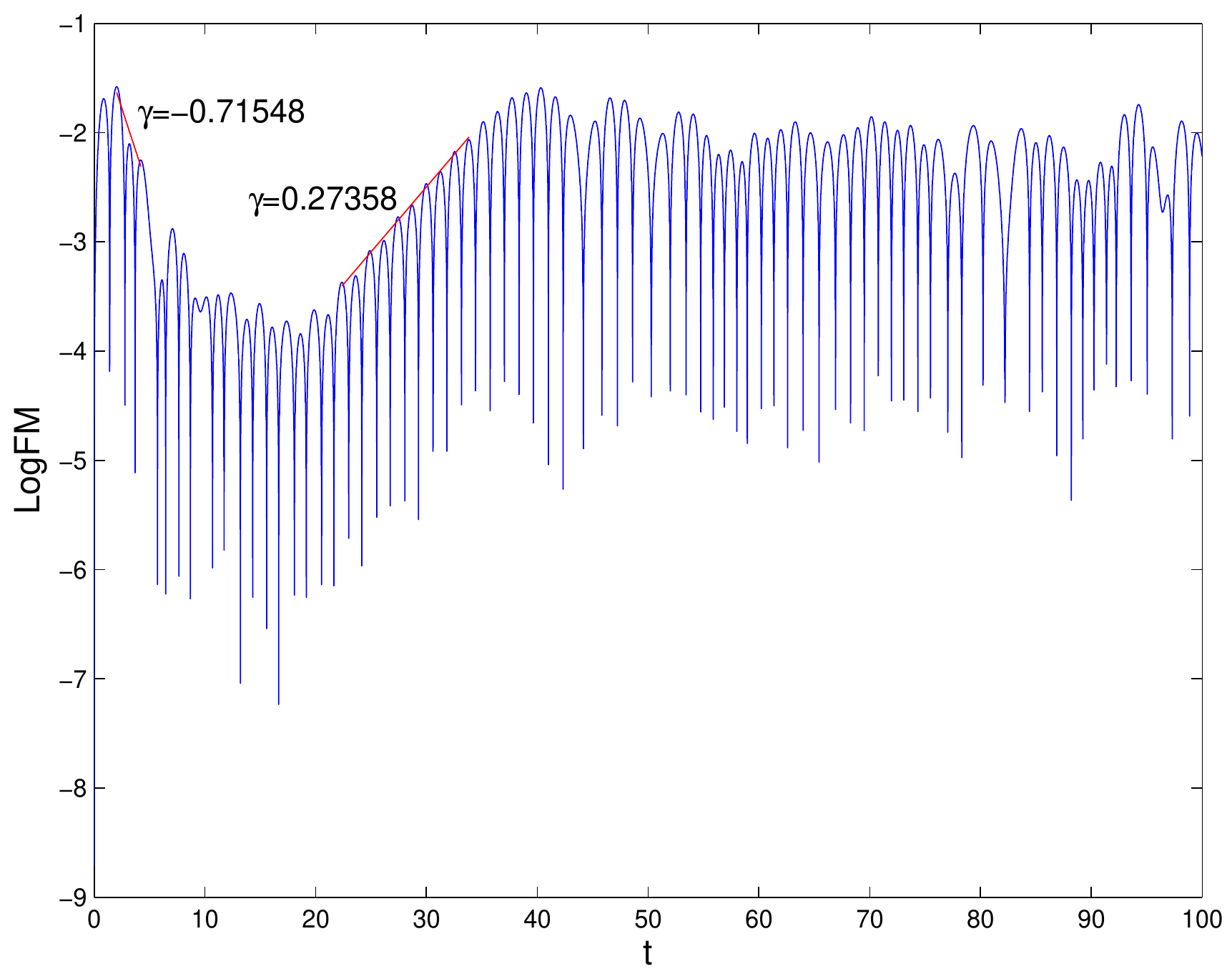}}
         \subfigure[ $logF\!M_2$, $\mathbb{Q}^2$ no limiter]{\includegraphics[width=2in,angle=0]{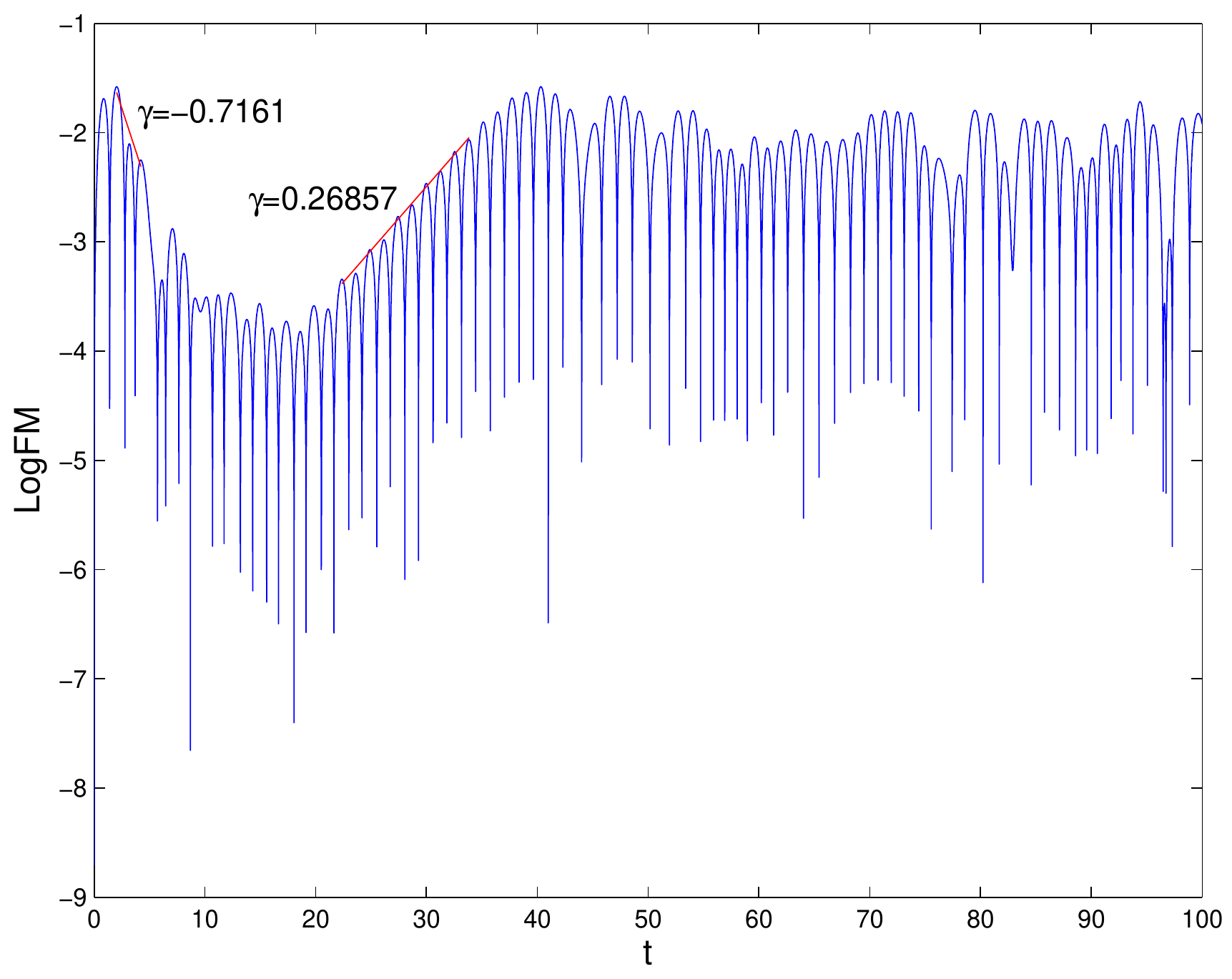}}
      \subfigure[ $logF\!M_3$, $\mathbb{P}^2$ no limiter]{\includegraphics[width=2in,angle=0]{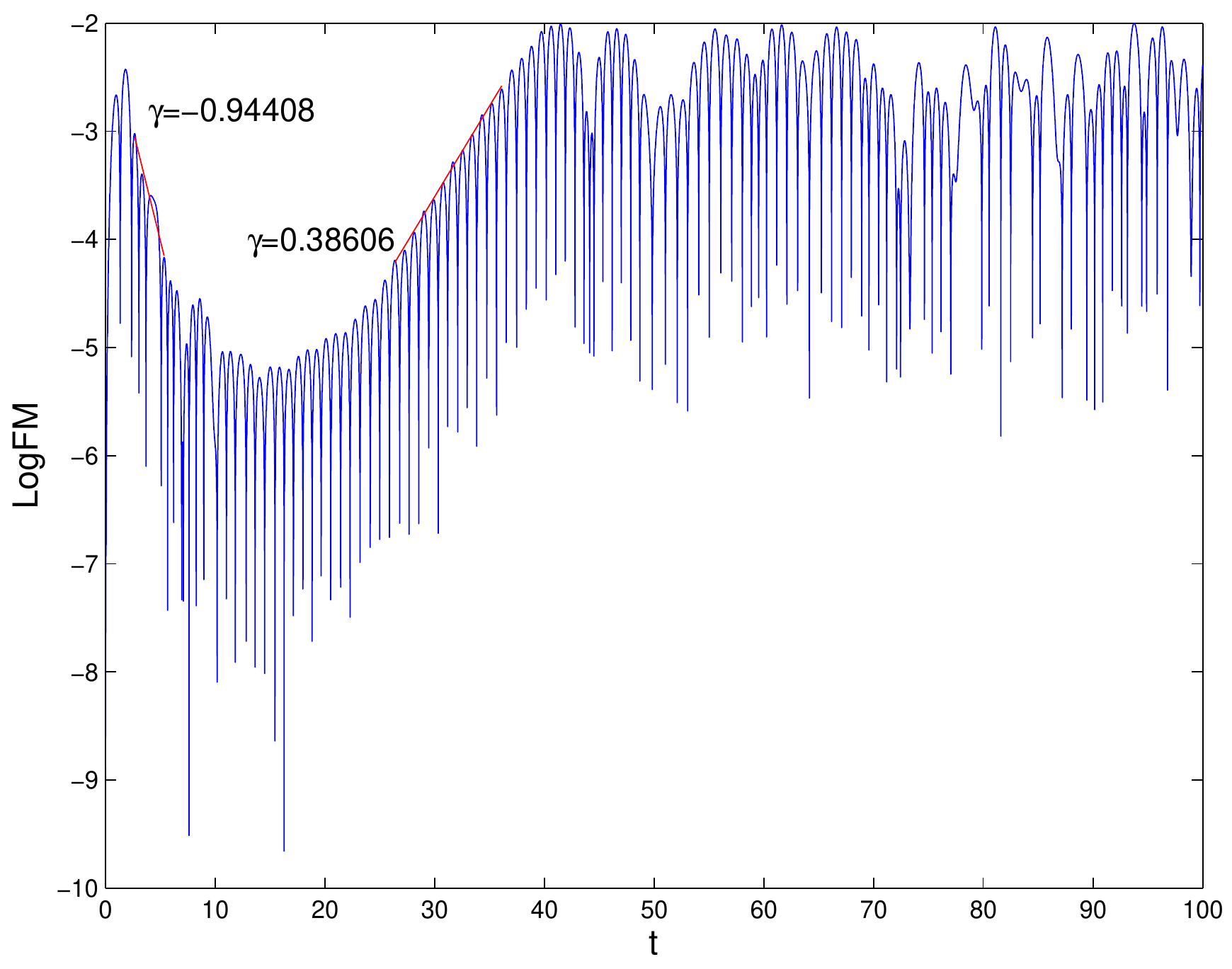}}
            \subfigure[ $logF\!M_3$, $\mathbb{P}^2$ with limiter]{\includegraphics[width=2in,angle=0]{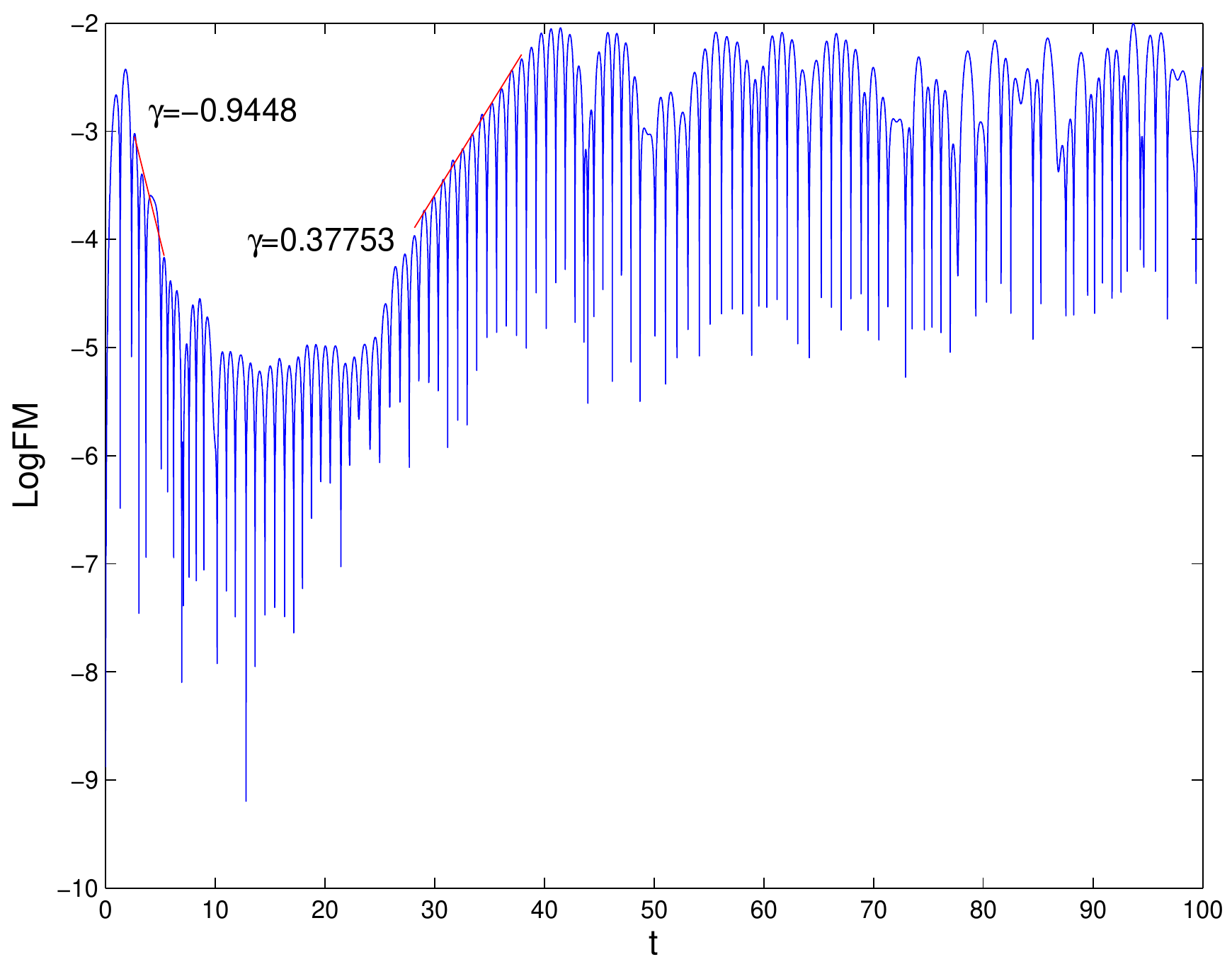}}
                  \subfigure[ $logF\!M_3$, $\mathbb{Q}^2$ no limiter]{\includegraphics[width=2in,angle=0]{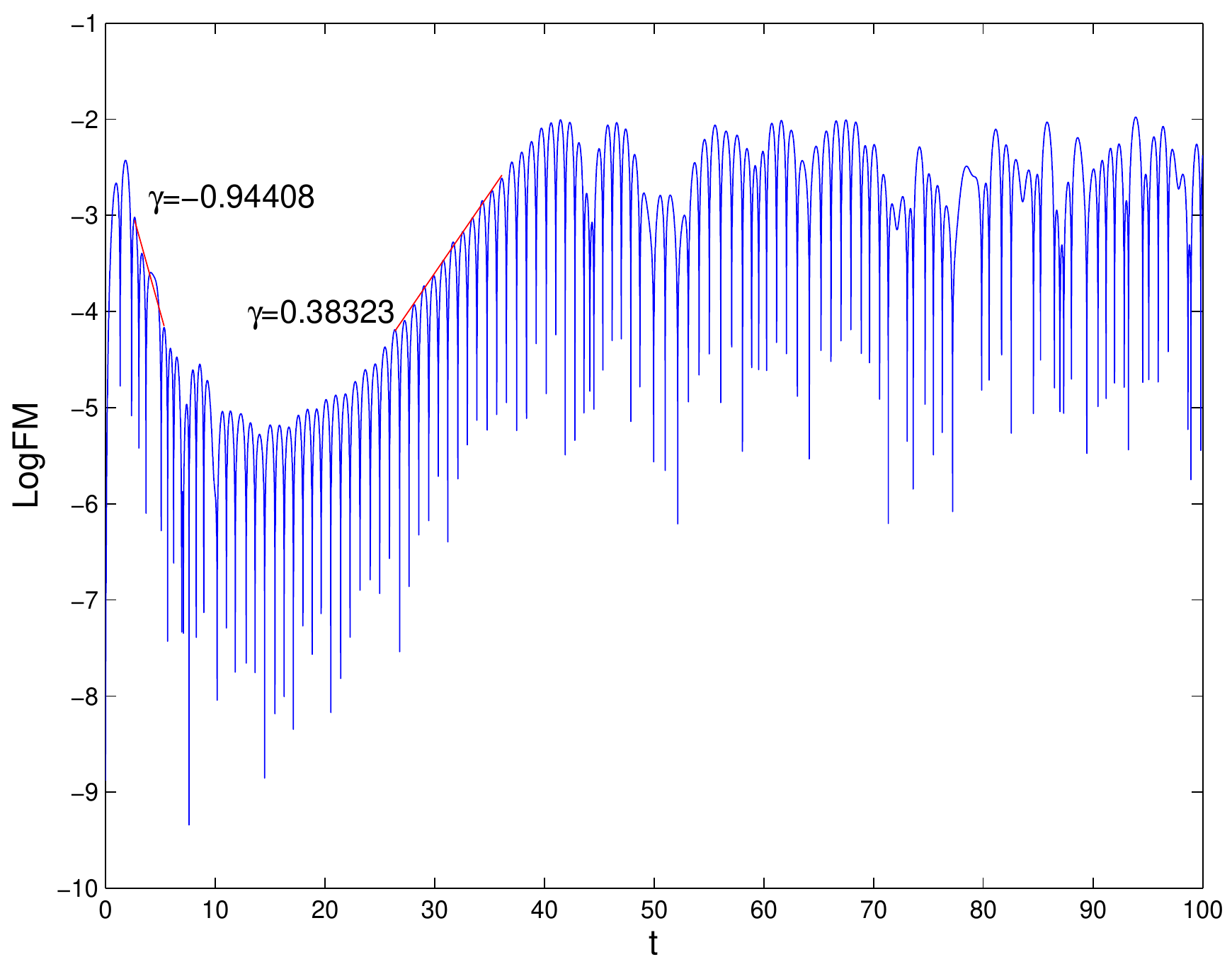}}
   \subfigure[ $logF\!M_4$, $\mathbb{P}^2$ no limiter]{\includegraphics[width=2in,angle=0]{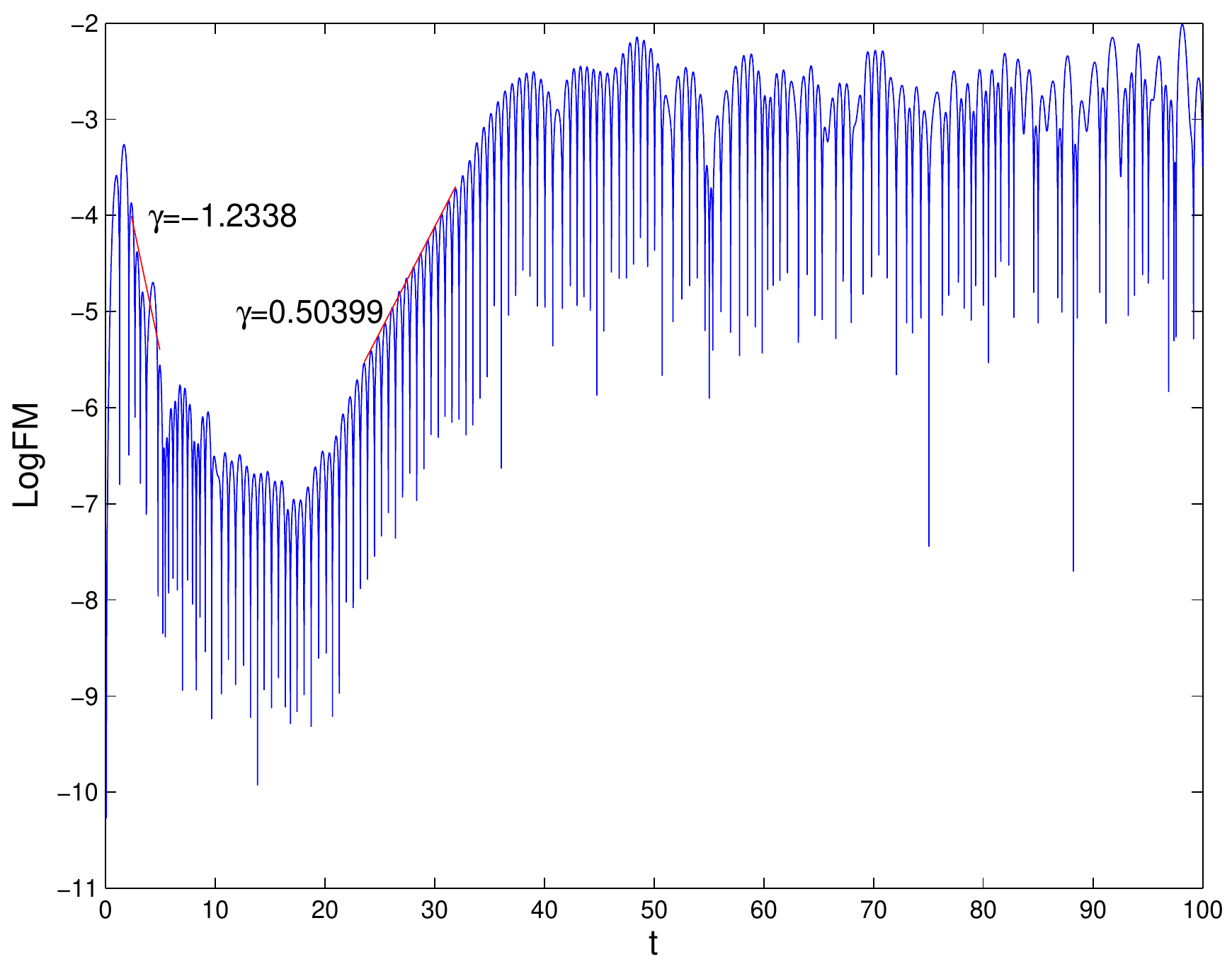}}
   \subfigure[ $logF\!M_4$, $\mathbb{P}^2$ with limiter]{\includegraphics[width=2in,angle=0]{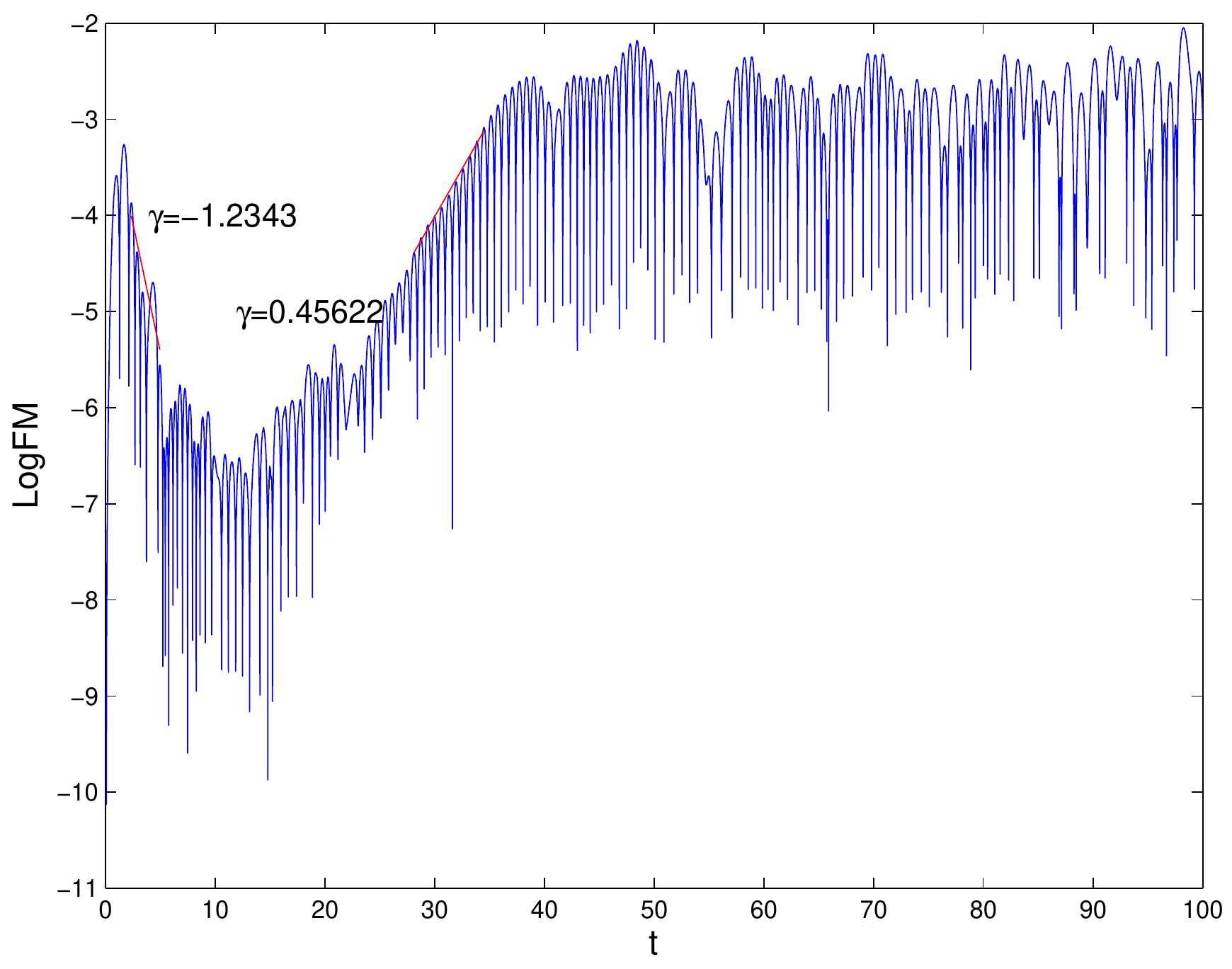}}
   \subfigure[ $logF\!M_4$, $\mathbb{Q}^2$ no limiter]{\includegraphics[width=2in,angle=0]{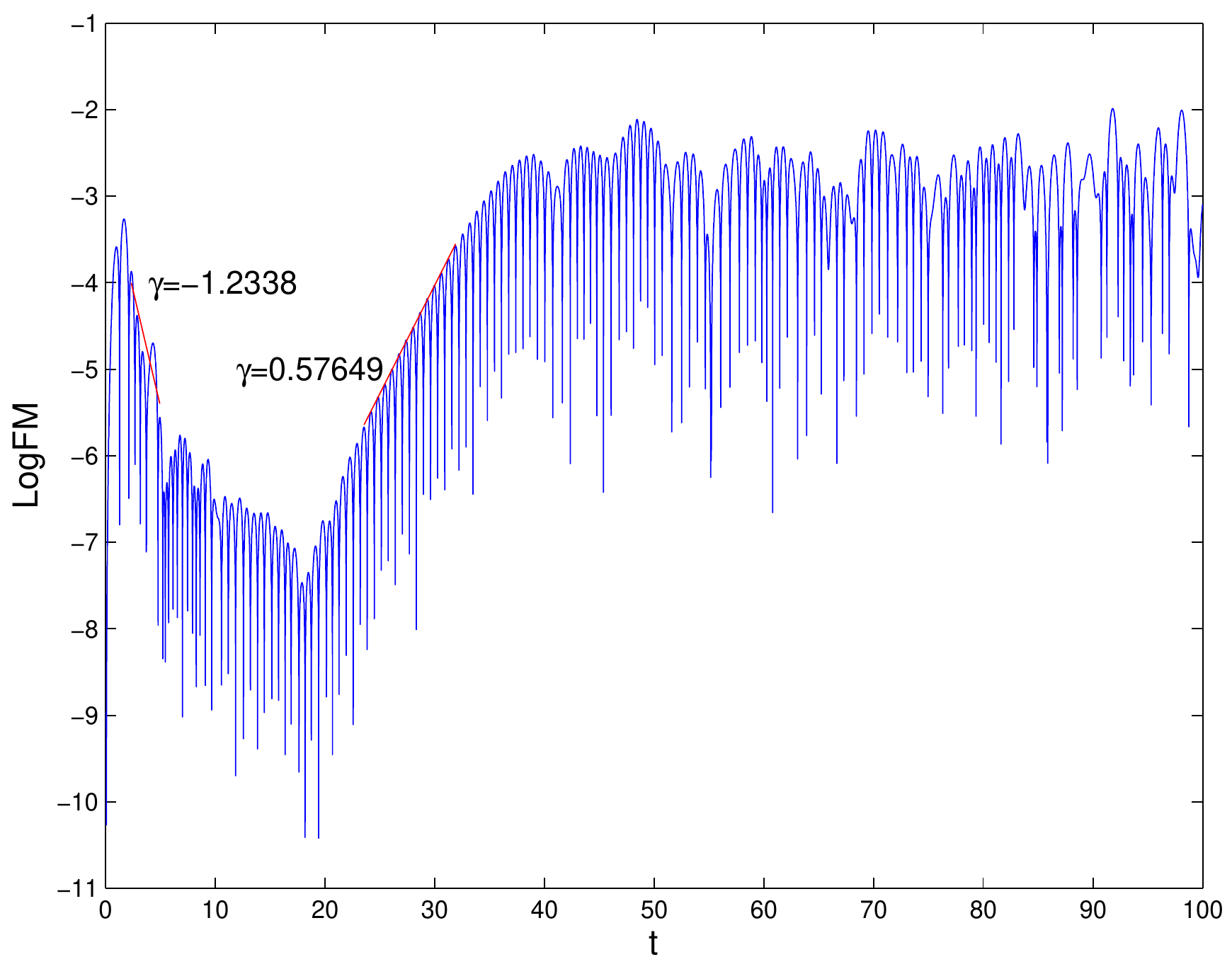}}
  \end{center}
  \caption{Evolution of  the first four Log Fourier mode as a function of time for nonlinear Landau damping.  Various values of the numerical damping/growth rate  are  marked on the graphs. Here the $\mathbb{P}^2$ space with the positivity-preserving limiter was used on a $100 \times 200$ mesh. The predicted recurrence time $T_R$ for $logF\!M_1$ is 209.44, for $logF\!M_2$ is 104.72, for $logF\!M_3$ is 69.81, and for $logF\!M_4$ is 52.36.}
\label{landaulfmp2n}
\end{figure}

For this case we choose $f_0(x,v)=f_M(v) (1+ A \cos (kx))$ with $A=0.5$, $k=0.5$, $L=4 \pi$, and $V_c=6$.  We implement the scheme on a $100 \times 200$ mesh and integrate up to  $T=100$ using  three methods: $\mathbb{P}^2$, $\mathbb{P}^2$ with the positivity-preserving limiter, and $\mathbb{Q}^2$. In Figure \ref{landaulfmp2n}, we plot the evolution of the first four Log Fourier modes as a function of time. All three methods give qualitatively similar  results that compare well with other calculations in the literature. We observe initial damping (until $t \approx 15$),   followed by exponential growth (until $t \approx 40$), and finally saturation of the modes.  Note the predicted recurrence times $T_R$ for each of the modes are as follows:  for $logF\!M_1$,  $T_R= 209.44$;  for $logF\!M_2$, $T_R=104.72$; for $logF\!M_3$,  $T_R=69.81$;  and for $logF\!M_4$,  $T_R=52.36$.    Since the bounce time is about 40, we have some confidence that the solution is resolved at least up u
 ntil nonlinearity becomes important.  Although, the role played by $T_R$ for the nonlinear evolution is not clear since nonlinearity could remove the fine scales generated by linear phase mixing.

\begin{figure}[htb]
  \begin{center}
      \subfigure[ charge]{\includegraphics[width=2.5in,angle=0]{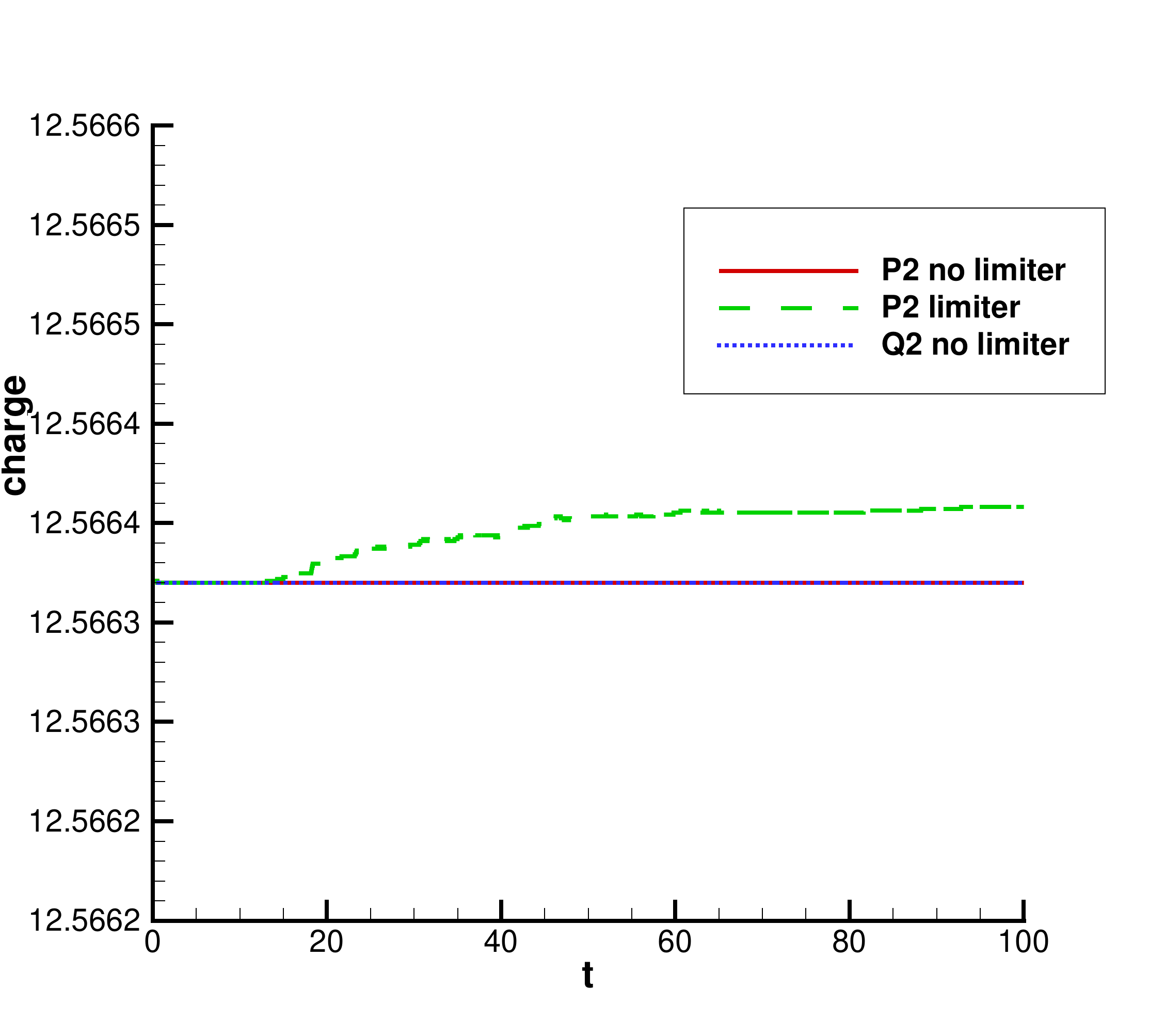}}
   \subfigure[ momentum]{\includegraphics[width=2.5in,angle=0]{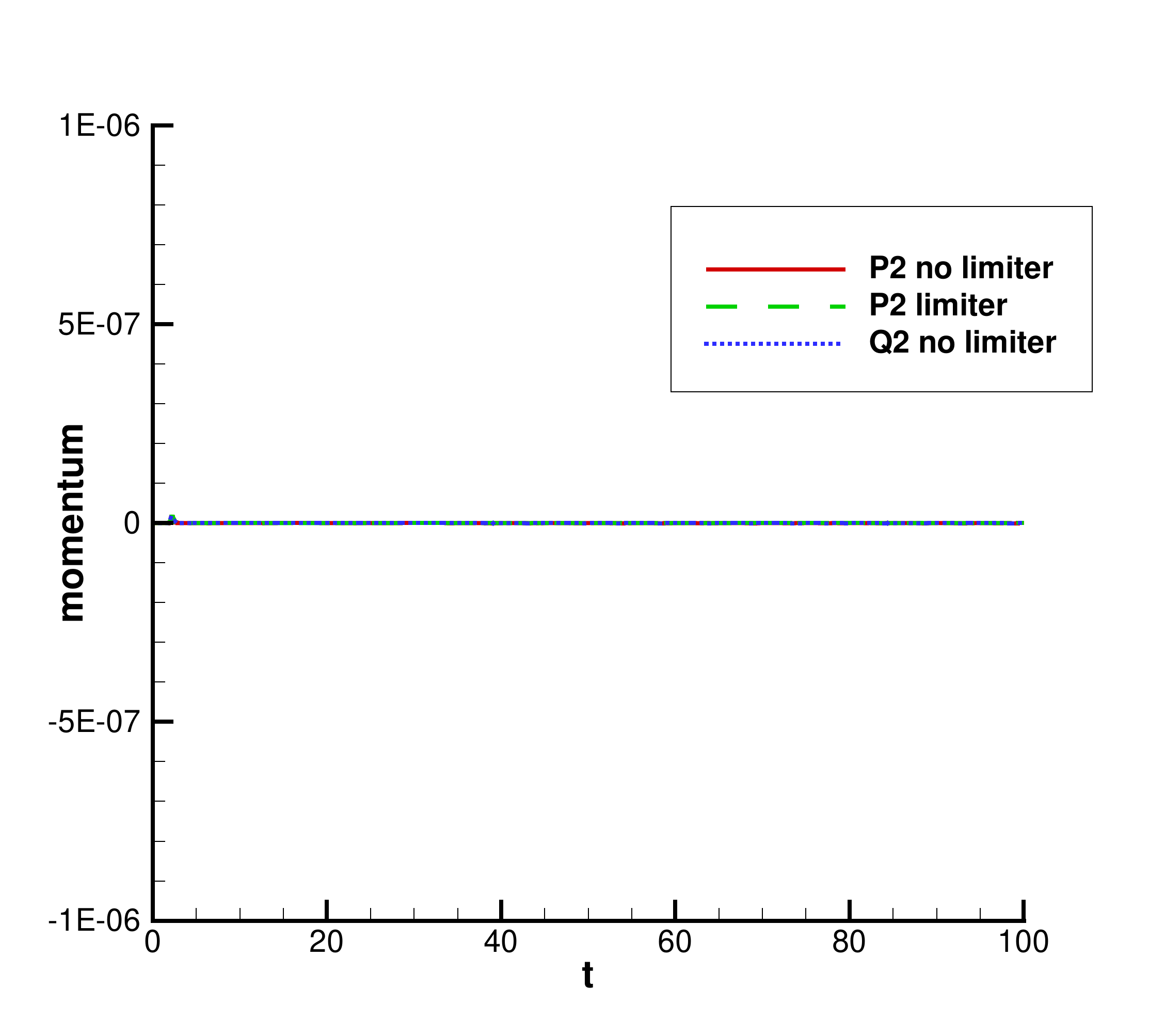}}
            \subfigure[ enstrophy]{\includegraphics[width=2.5in,angle=0]{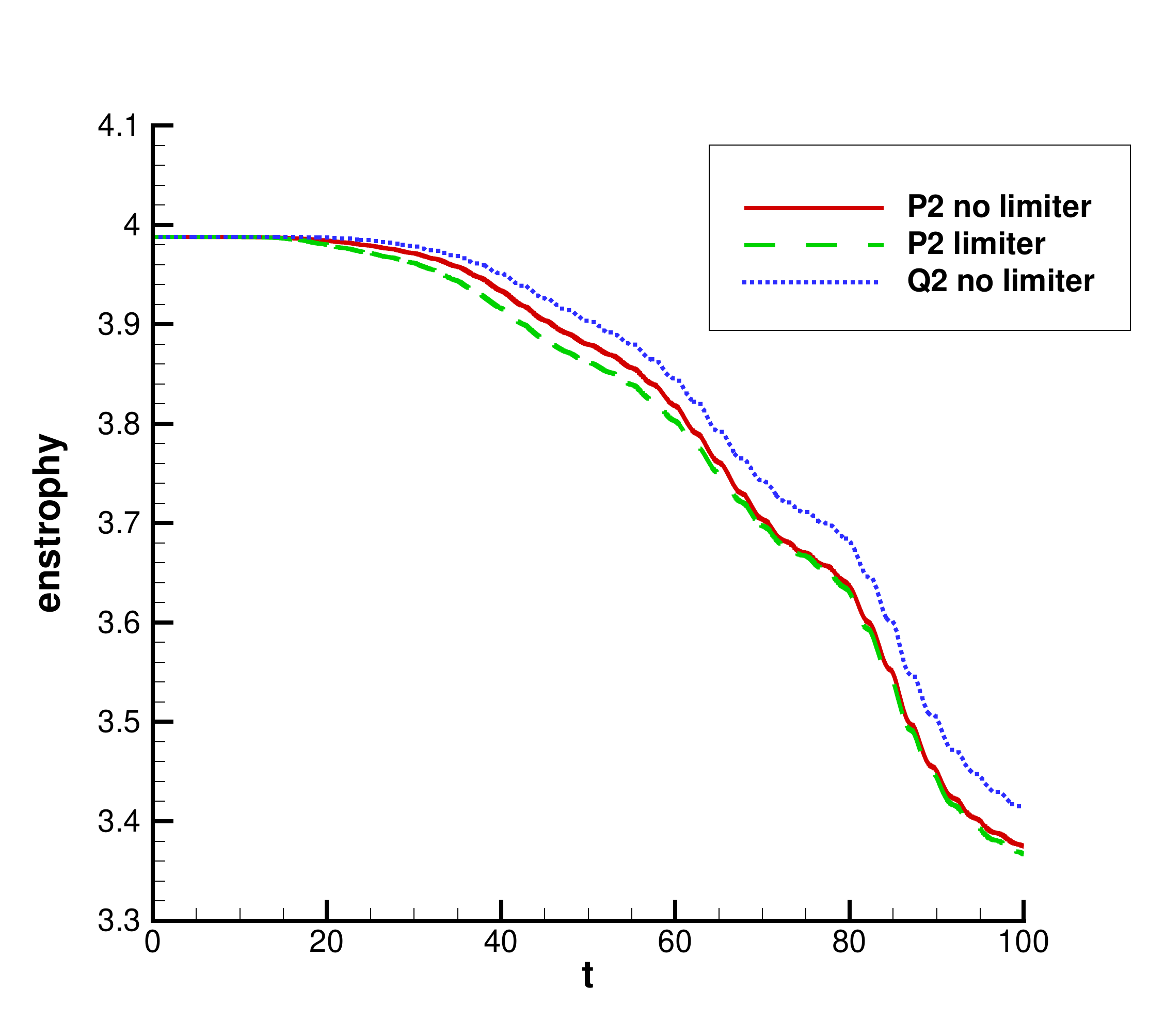}}
         \subfigure[ total energy]{\includegraphics[width=2.5in,angle=0]{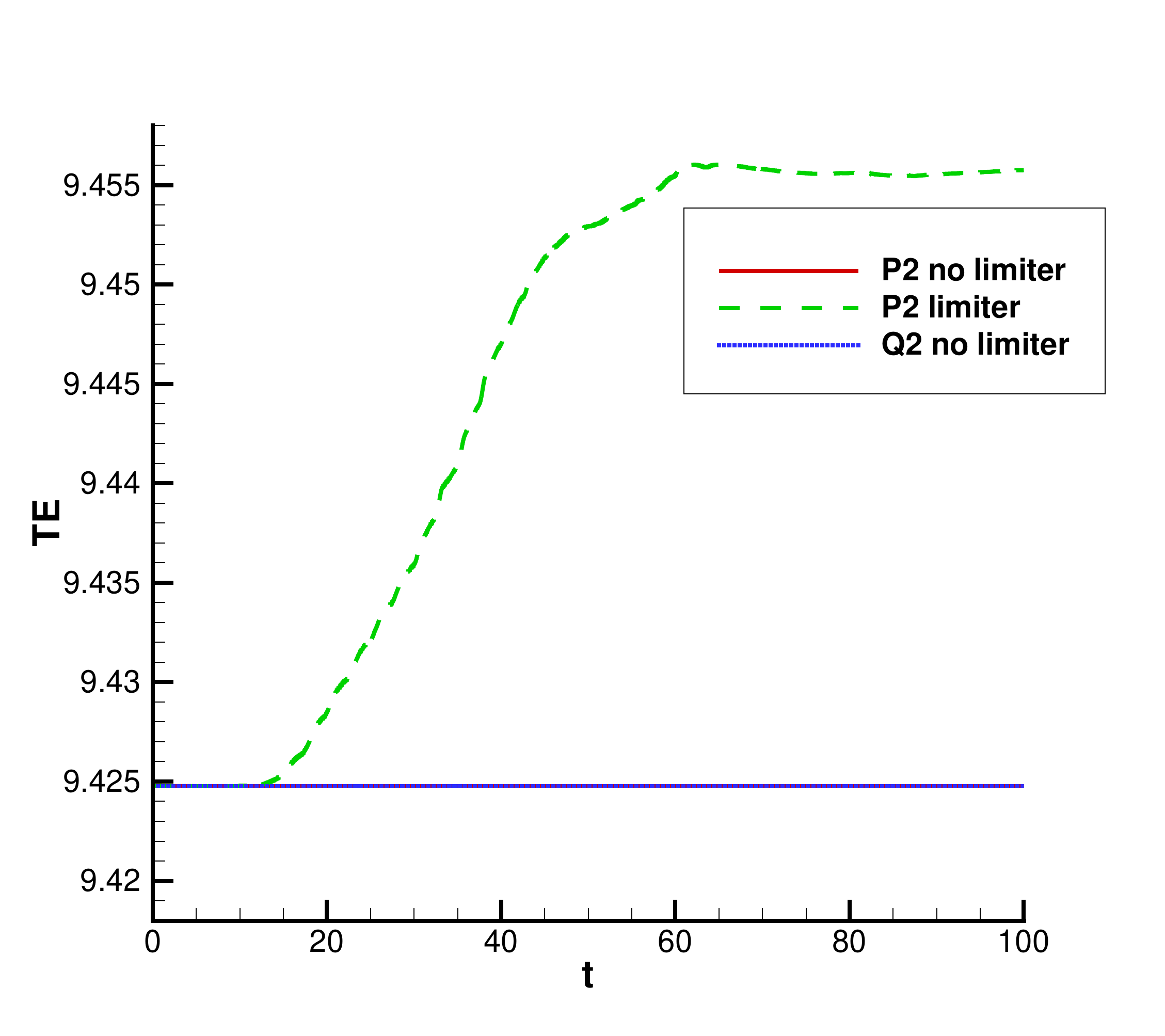}}
  \end{center}
  \caption{ Evolution of  conserved quantities as a function of time during the course of nonlinear Landau damping for various computational methods. A mesh of   $100 \times 200$ was used.}
\label{landaumacp2n}
\end{figure}

In Figure \ref{landaumacp2n}, we plot the conserved quantities of Section \ref{cons}. The charge and momentum are well  conserved for all methods, while the enstrophy has decayed by about 15\% at $T=100$ for all three methods. This result agrees with our analysis in Section 2. We remark that the limiter has an effect  on charge conservation, due to its modification of
 the solution on the boundary. The total energy is conserved much better without the positivity-preserving limiter. When we use the limiter, the total energy grows by about 0.3\% at $T=100$.

\subsubsection*{Two-stream instability}

\begin{figure}[htb]
  \begin{center}
        \subfigure[ $logF\!M_{1,2,3,4}$, $\mathbb{P}^2$ no limiter]{\includegraphics[width=2in,angle=0]{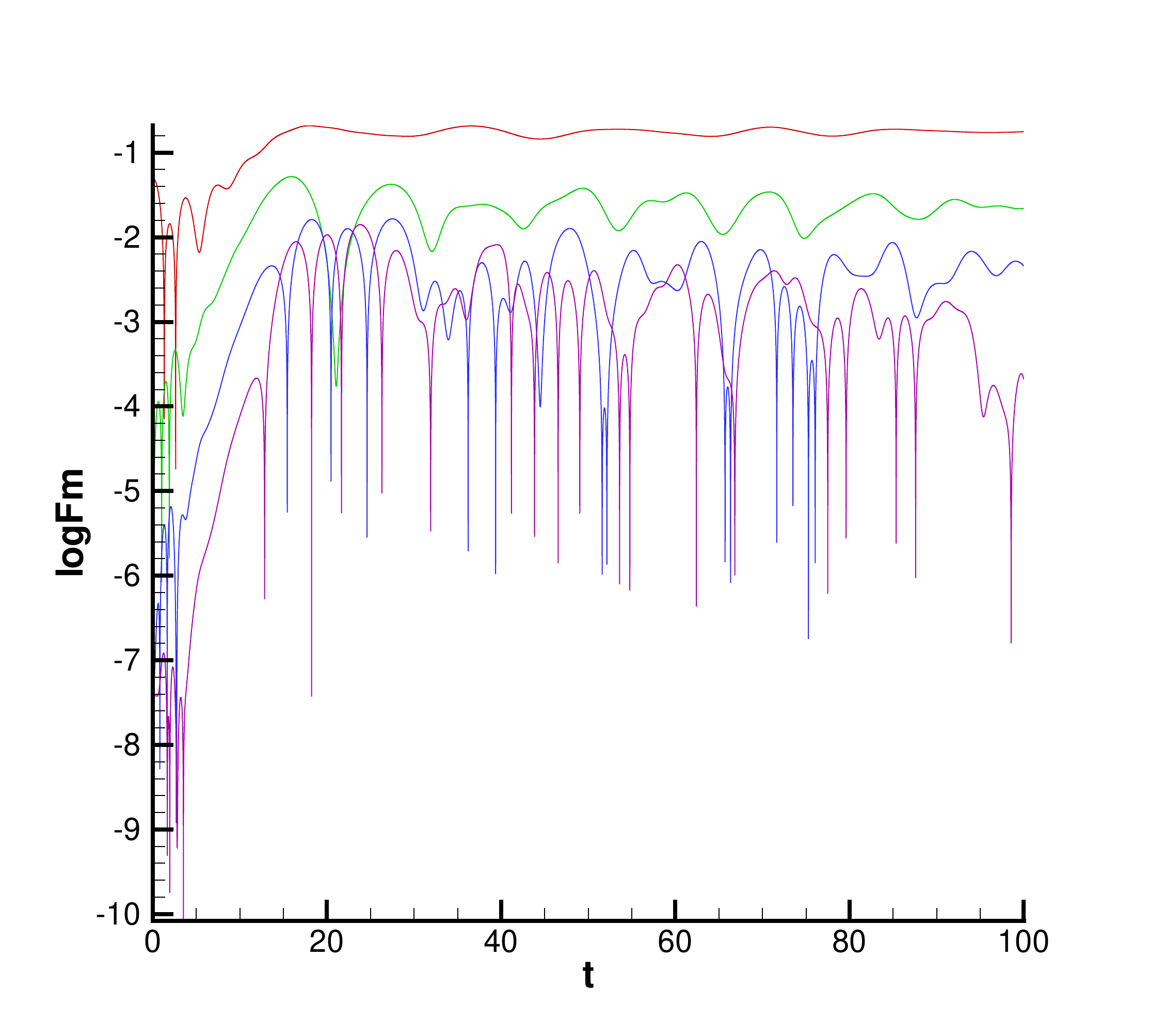}}
                    \subfigure[ $logF\!M_{1,2,3,4}$, $\mathbb{P}^2$ with limiter]{\includegraphics[width=2in,angle=0]{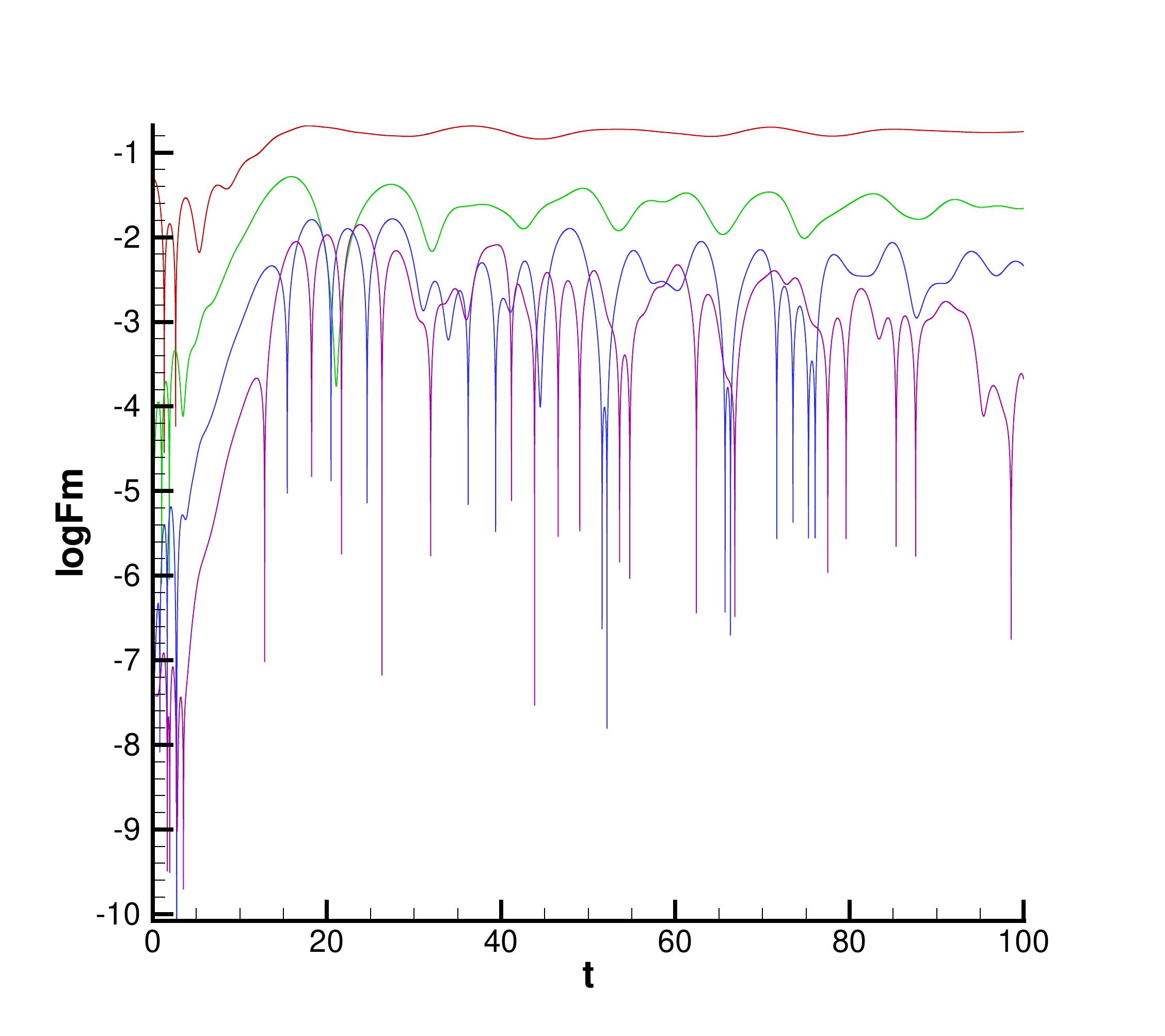}}
                  \subfigure[ $logF\!M_{1,2,3,4}$, $\mathbb{Q}^2$ no limiter]{\includegraphics[width=2in,angle=0]{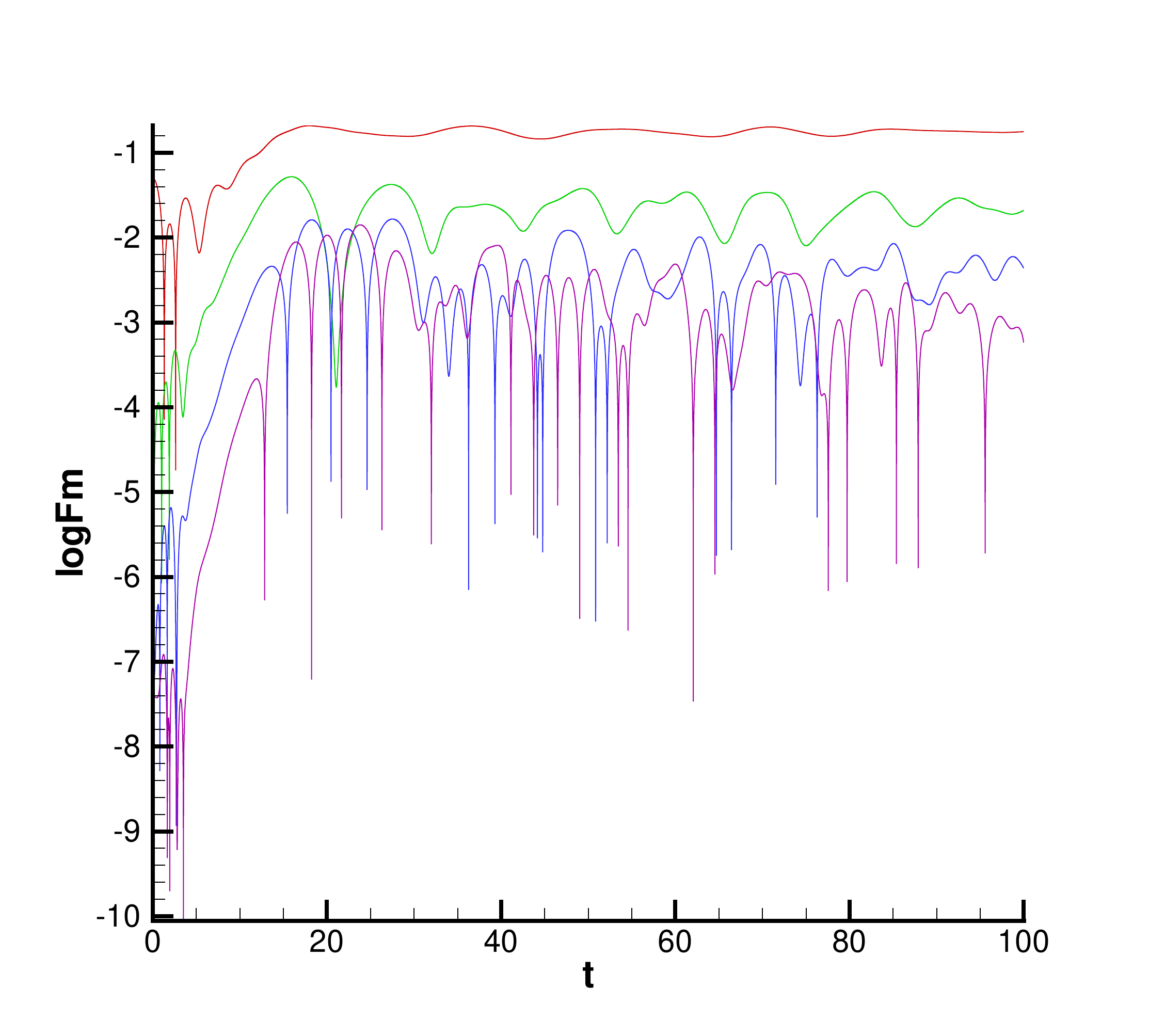}}
      \subfigure[ enstrophy]{\includegraphics[width=2.5in,angle=0]{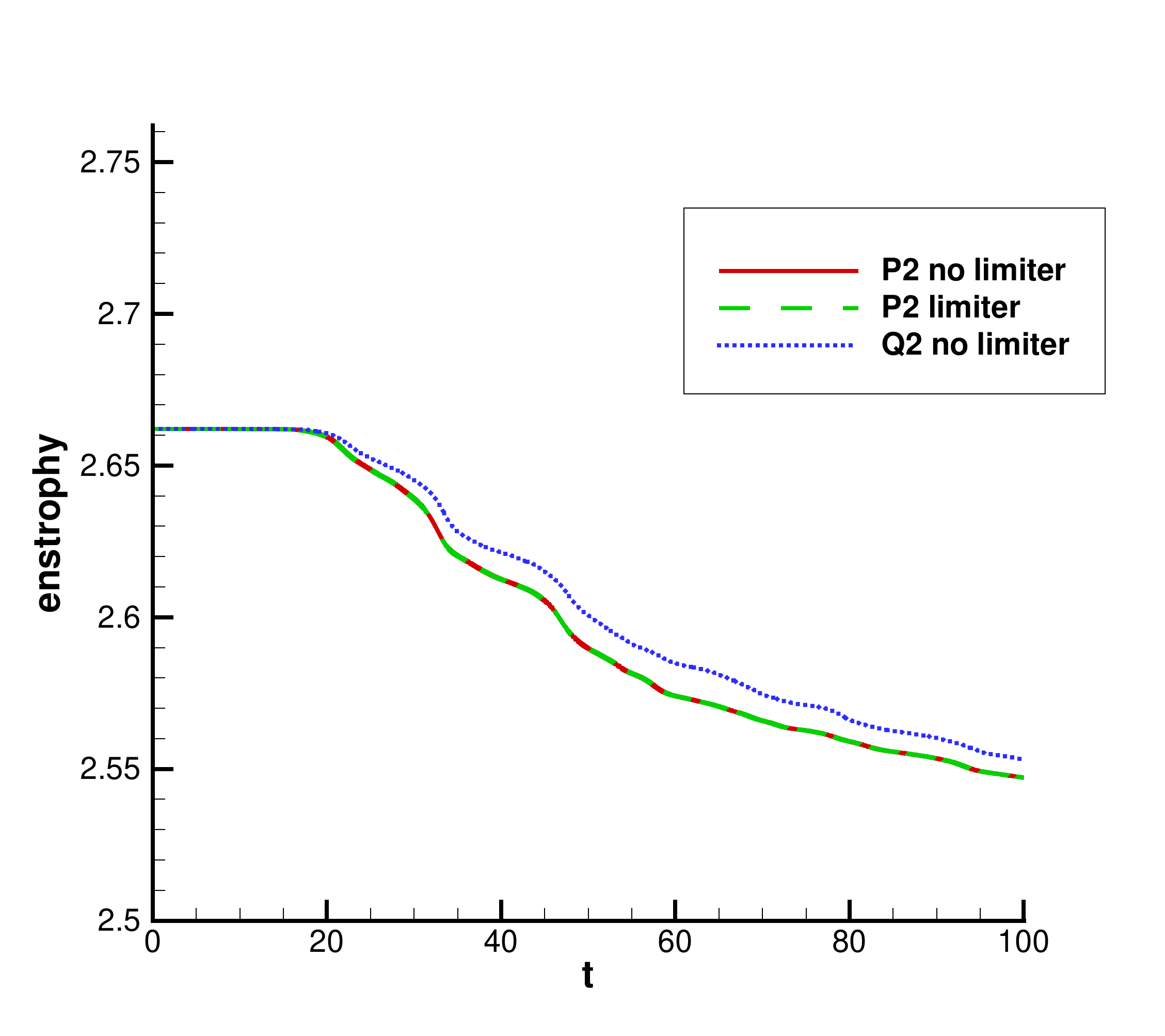}}
         \subfigure[ total energy]{\includegraphics[width=2.5in,angle=0]{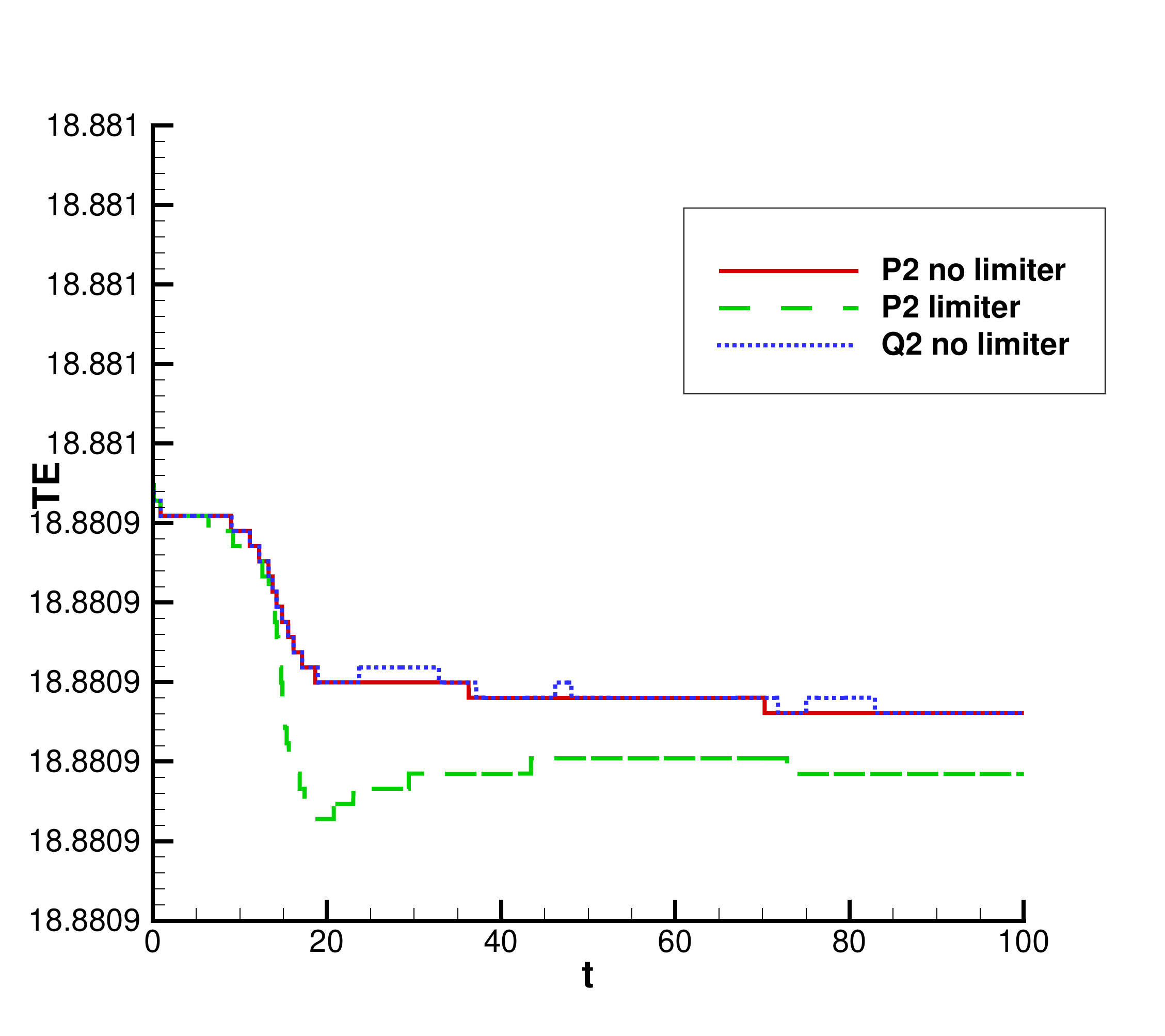}}
  \end{center}
  \caption{Depiction of the first four Log Fourier modes during the nonlinear evolution of the two-stream instability.  Also depicted is the evolution of  energy and enstrophy as a function of time for various methods.  }
\label{tsmac}
\end{figure}

\begin{figure}[htb]
  \begin{center}
        \subfigure[ $t=0$]{\includegraphics[width=3in,angle=0]{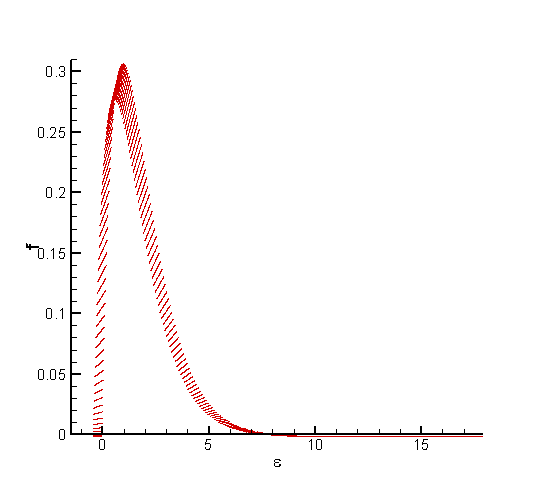}}
    \subfigure[ $t=50$]{\includegraphics[width=3in,angle=0]{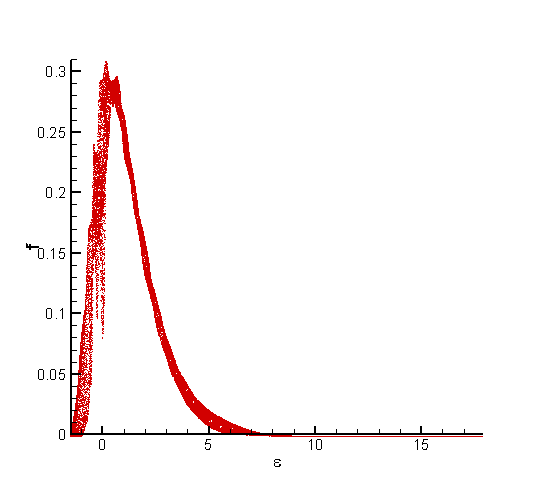}}
      \subfigure[ $t=100$]{\includegraphics[width=3in,angle=0]{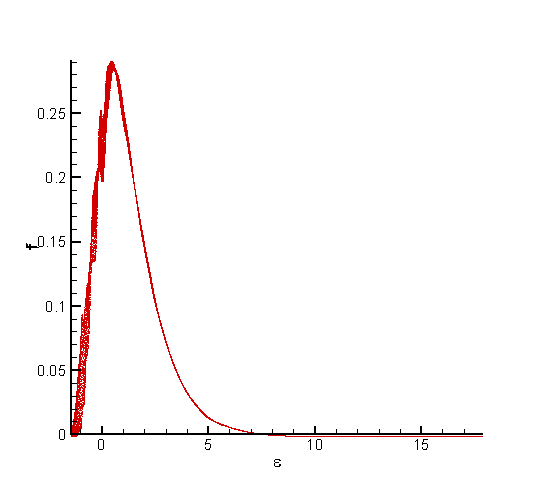}}
      \subfigure[ $t=200$]{\includegraphics[width=3in,angle=0]{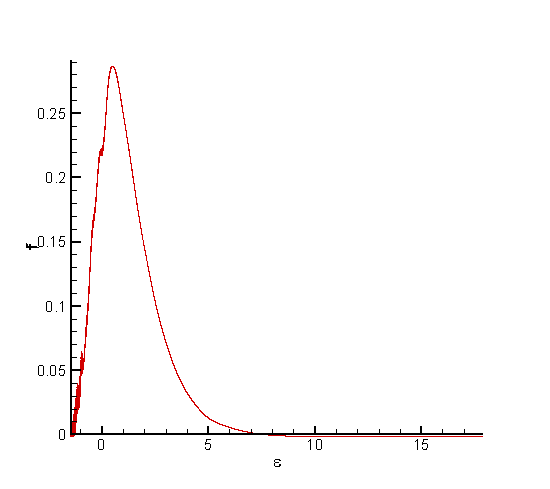}}
  \end{center}
  \caption{Plots of the distribution $f(x,v,t)$ versus $\epsilon=\frac{v^2}{2}-\Phi(x, t)$ for the two-stream instability at the times $t$ indicated, showing saturation to a BGK state.  Here the $\mathbb{P}^2$ space with the positivity-preserving limiter was used on a $100 \times 200$ mesh. }
\label{tsbgk}
\end{figure}

For this case we choose  $f_0(x,v)=f_{TS}(v) (1+ A \cos (kx))$, where $f_{TS}(v)=\frac{1}{\sqrt{2 \pi}} v^2 e^{-v^2/2}$, $A=0.05$, $k=0.5$, $L=4 \pi$, and $V_c=6$. The mesh size we take is $100 \times 200$. In Figure \ref{tsmac}, we plot the evolution of  conserved quantities. For this example, charge  and momentum are well conserved by all methods,  so are not plotted. The enstrophy decays by about 4\% at $T=100$, while the  total energy is well conserved even with the limiter. The plots of the log Fourier modes show an early exponential growth followed by  oscillation. Figure \ref{tsbgk} provides evidence that the system has relaxed into a  BGK mode.  Here,  the relation defined by the ordered pair $(\epsilon=v^2/2+\Phi(x, T), f(x, v, t))$) is plotted at various times $t$. The use of this  kind of  plot  as a diagnostic  was first reported in \cite{Heath} for electrostatic VP equations and later in \cite{Cheng_jeans} for the gravitational VP equations. Here, the evolution
 clearly indicates convergence to a BGK equilibrium.

\subsubsection*{Dynamically accessible excitations--KEEN waves}

Motivated by experiments  performed for understanding aspects of laser-plasma interaction  \cite{Mont},  several  authors  have considered numerical solution of the VP system with a transitory external driving electric field (see, e.g., \cite{Afeyan_03, Johnston_09} ), rather than just specifying an ad hoc initial condition for $f$, as is usually done.  Such drive generated initial conditions are examples of those  proposed and discussed in  \cite{Morrison_89, Morrison_90,Morrison_92},  where they were termed dynamically accessible (DA) initial conditions.   DA initial conditions are important because they have a Hamiltonian origin and preserve phase space constraints.  Moreover, since ultimately any perturbation of charged particles within the confines of VP theory is in fact caused by an electric field,  it is physically very natural to consider DA initial conditions.  We consider two numerical examples and compare our results with those of  \cite{Afeyan_03, Johnston_09}, w
 here the authors observed saturation to nonlinear traveling BGK-like states that they termed KEEN waves, standing for kinetic electrostatic electron nonlinear waves.   We note that the calculations of \cite{Afeyan_03} were duplicated in  \cite{Heath_thesis} and allied work was given in \cite{Valentini,Valentini2}.

Specifically, the system is driven by a single prescribed frequency and wavelength, where  the driven Vlasov equation,
\[
f_t+v f_x - (E+E_{ext})f_v=0\,,
\]
is solved.  Here, $E_{ext}(x, t)=A_d(t) \sin(kx-\omega t)$ is the external field, where $A_d$  is  a temporal envelope that is ramped up to a plateau and then ramped down to zero. For our  two examples we  consider the following two ramping profiles:
\begin{equation}
A^{J}_d(t)= \left \{
\begin{array}{ll}
A_m\,  \sin(t \pi/100)  & \textrm{if} \quad 0< t  < 50\\ \\
A_m
 &\textrm{if} \quad 50\leq  t <150\\ \\
 A_m\,  \cos\big((t-150)\pi/100\big)
 &\textrm{if} \quad 150< t  < 200\\ \\
 0 &\textrm{if} \quad 200< t <T
 \end{array}
\right .\,,
\label{rampJ}
\end{equation}
with $A_m= 0.052$ as used  in  \cite{Johnston_09} and
\begin{equation}
A^{A}_d(t)= \left \{
\begin{array}{ll}
A_m\,  \frac{1}{1+e^{-40 (t-10)}}
 & \textrm{if} \quad 0< t  < 60\\ \\
A_m \left(1-\frac{1}{1+e^{-40 (t-110)}}\right)
 &\textrm{if} \quad 60\leq  t <T
\end{array}
\right .\,,
\label{rampA}
\end{equation}
with $A_m=0.4$ as used in \cite{Afeyan_03}.   
In practice, the  system is initialized on  $f(0,x,v)=f_M(v)$, then ramped according to (\ref{rampJ}) or (\ref{rampA})  to prepare the DA initial condition.  The system is then evolved after $E_{ext}$ is  turned off  and seen to approach asymptotic states.
For both cases the computational domain is of size $[0, 2 \pi/k] \times [-8,8]$, and we  take $k=0.26$ and  $\omega=0.37$.

%

%
%

\begin{figure}[htb]
  \begin{center}
\includegraphics[width=3in,angle=0]{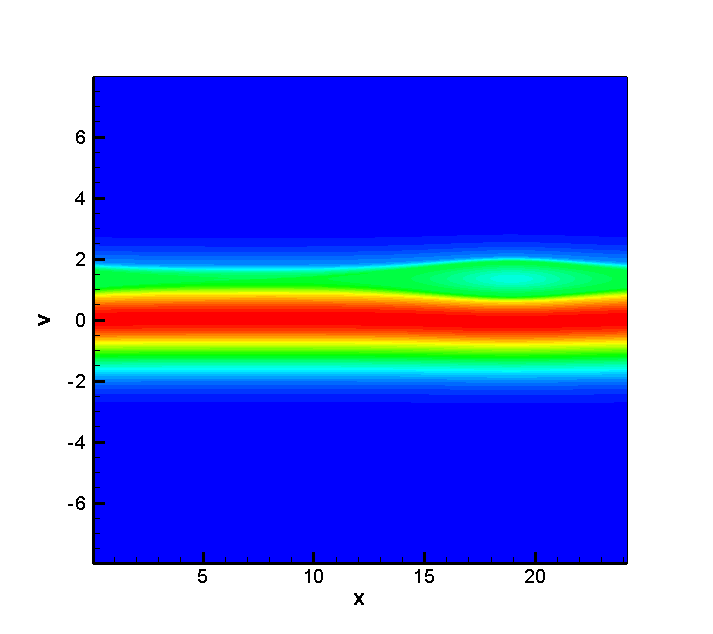}
  \end{center}
  \caption{Phase space contour at $T=1000$ for with DA initial condition with drive $A_d^{J}$.   The plot suggests saturation to a moving BGK-like state. Here the $\mathbb{Q}^1$ element  was used on a $200 \times 400$ mesh.}
\label{Jcont}
\end{figure}
\begin{figure}[htb]
  \begin{center}
\includegraphics[width=3in,angle=0]{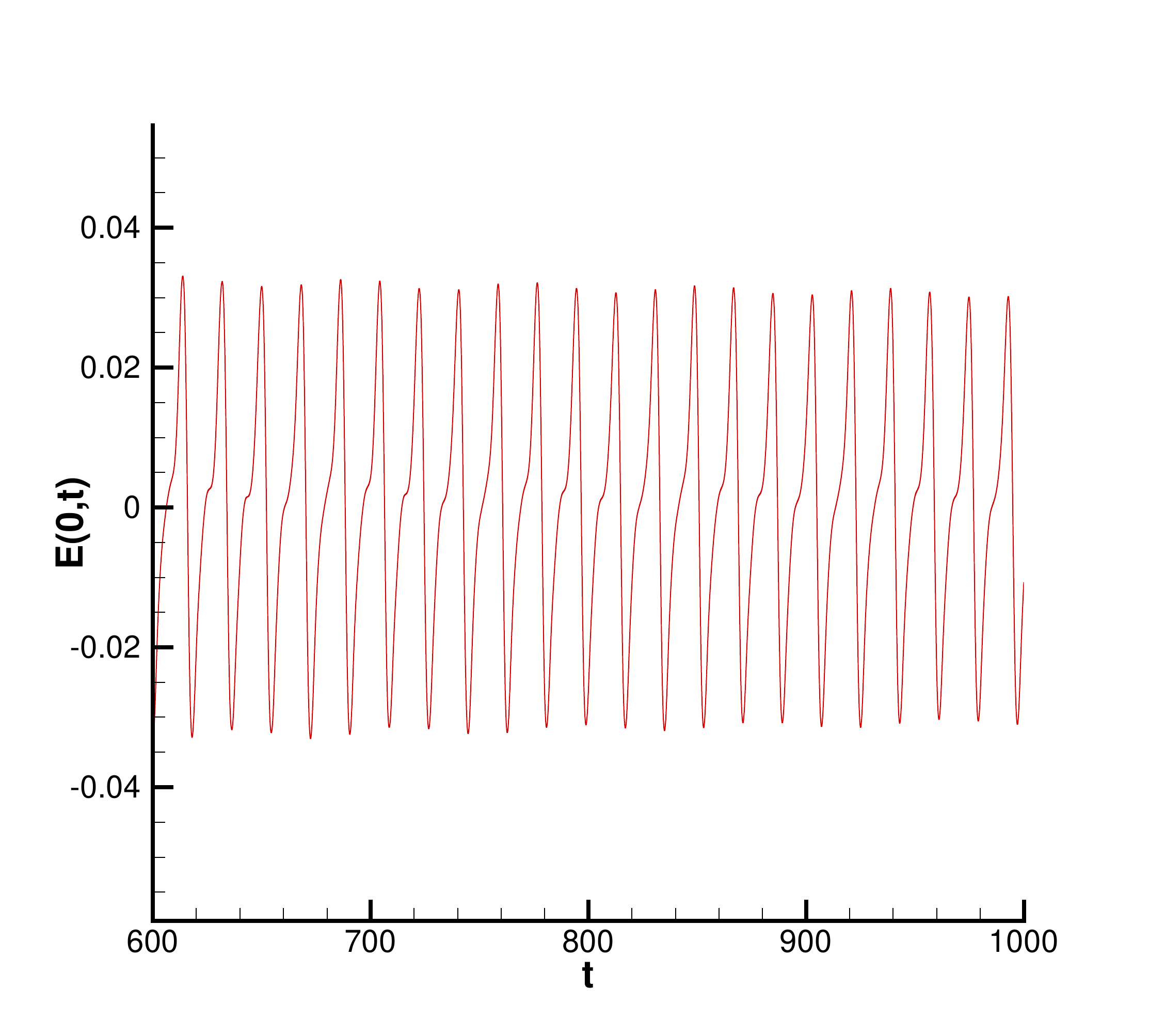}
  \end{center}
  \caption{The electric field E(0,t) at the center of the spacial domain at late times for the drive $A_d^{J}$.  The  periodicity matches the propagation of the BGK-like state through the domain. Here the $\mathbb{Q}^1$ element  was used on a $200 \times 400$ mesh.}
\label{JE}
\end{figure}
\begin{figure}[htb]
  \begin{center}
\includegraphics[width=3in,angle=0]{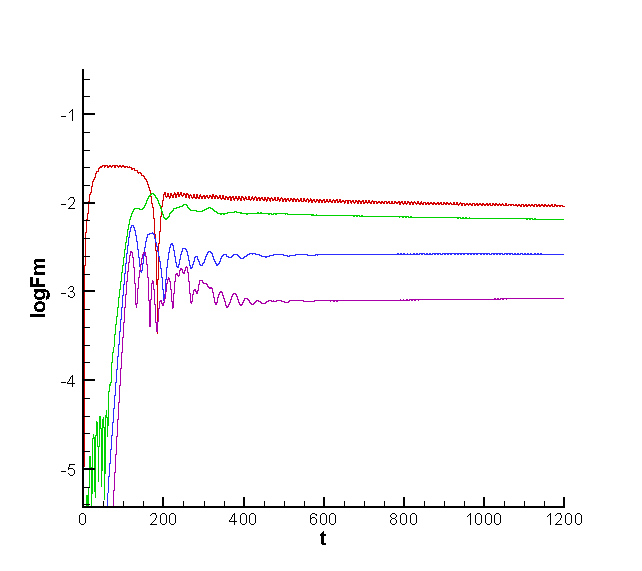}
  \end{center}
  \caption{The first four Log Fourier Modes for the drive $A_d^{J}$, indicating  saturation. Here the $\mathbb{Q}^1$ element  was used on a $200 \times 400$ mesh.}
\label{JlogF}
\end{figure}

Following  \cite{Johnston_09} with the drive $A_d^{J}$ of (\ref{rampJ}) with $A_m=0.4$  we obtain for latter times a translating BKG-like state, a snapshot of which is depicted in the phase space portrait of Fig.~\ref{Jcont}.  This structure moves through the spatial domain giving rise to the central periodic electric field signal, $E(0,t)$,  depicted in Fig.~\ref{JE}.   The period of this signal coincides with the propagation speed of the BKG-like state, which in agreement with  \cite{Johnston_09} is about 1.35.  Figure \ref{JlogF}  shows the first four Fourier modes and indicates saturation.

\begin{figure}[htb]
  \begin{center}
        \subfigure[ $t=0$]{\includegraphics[width=2in,angle=0]{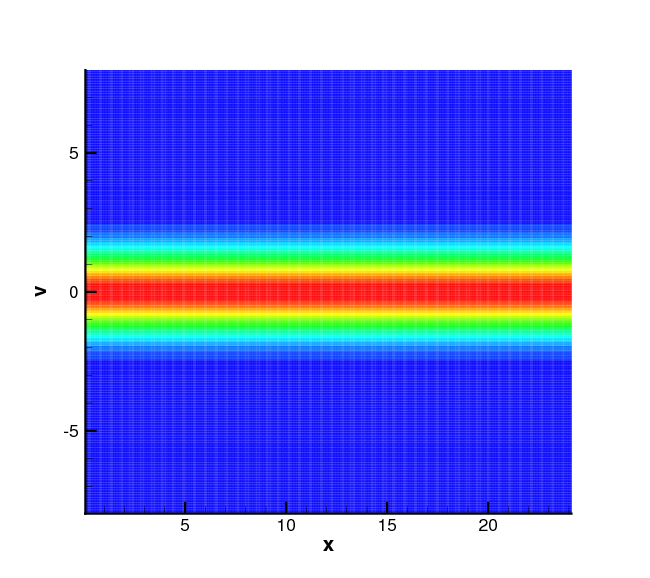}}
                    \subfigure[ $t=15$]{\includegraphics[width=2in,angle=0]{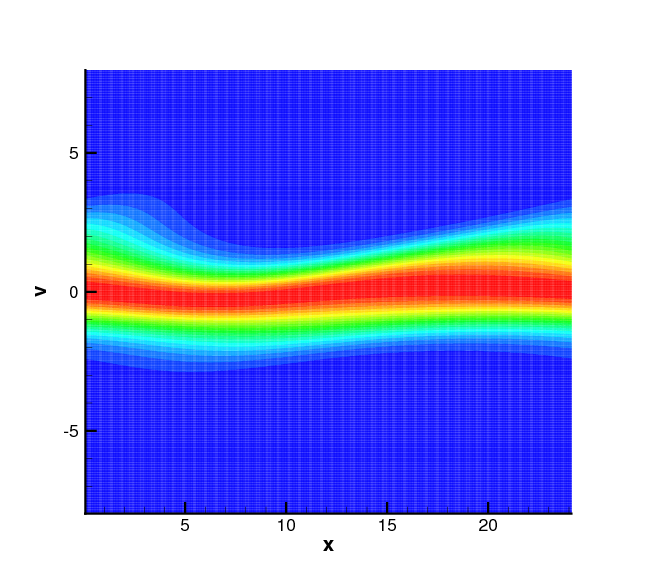}}
            \subfigure[ $t=30$]{\includegraphics[width=2in,angle=0]{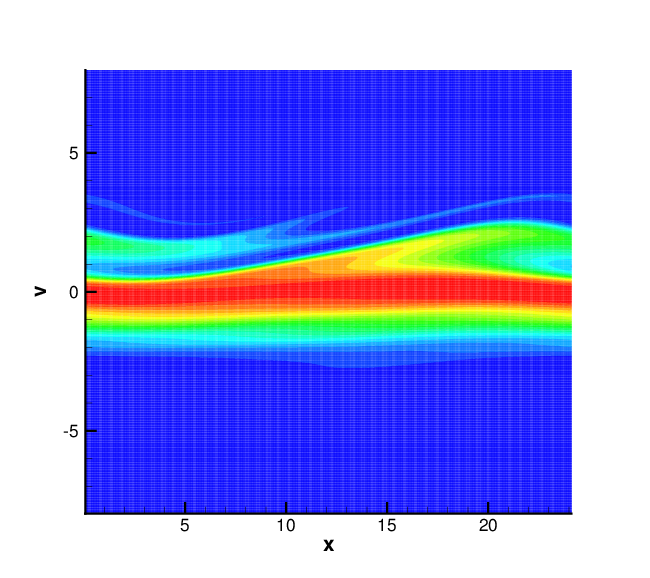}}
    \subfigure[ $t=45$]{\includegraphics[width=2in,angle=0]{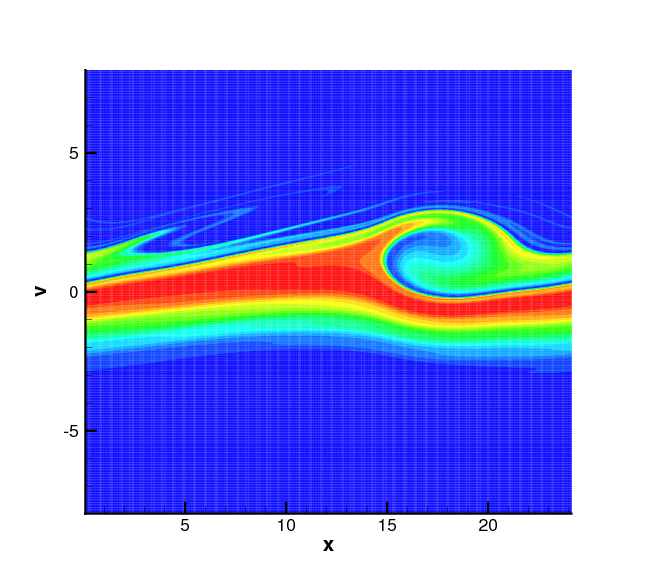}}
            \subfigure[ $t=60$]{\includegraphics[width=2in,angle=0]{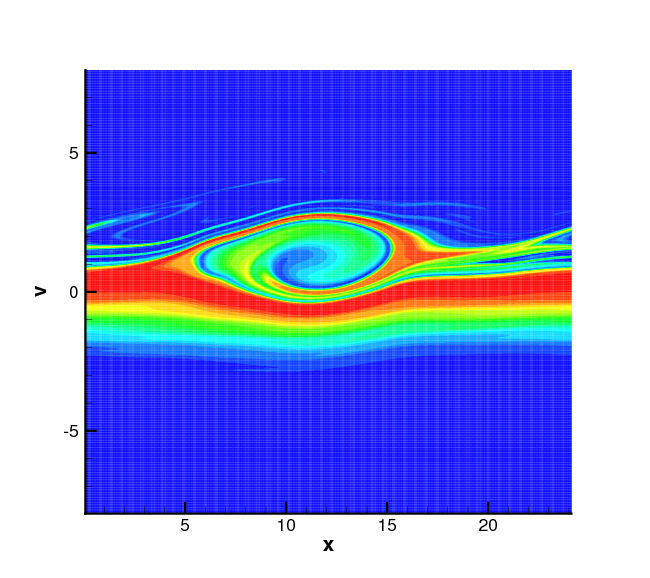}}
            \subfigure[ $t=90$]{\includegraphics[width=2in,angle=0]{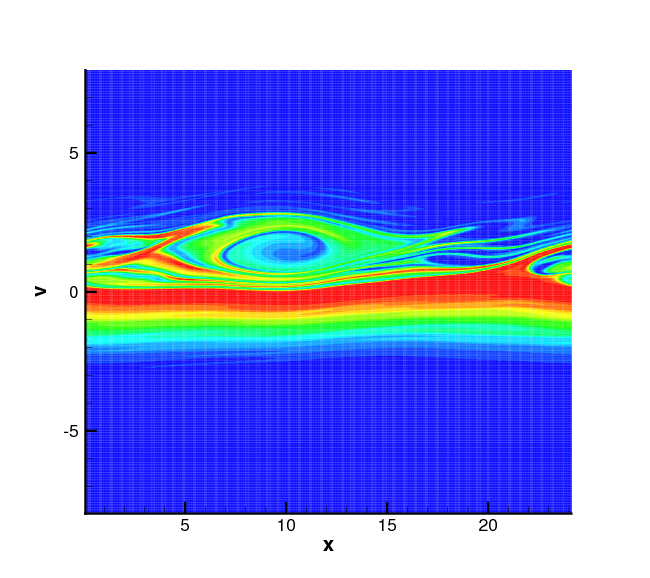}}
    \subfigure[ $t=105$]{\includegraphics[width=2in,angle=0]{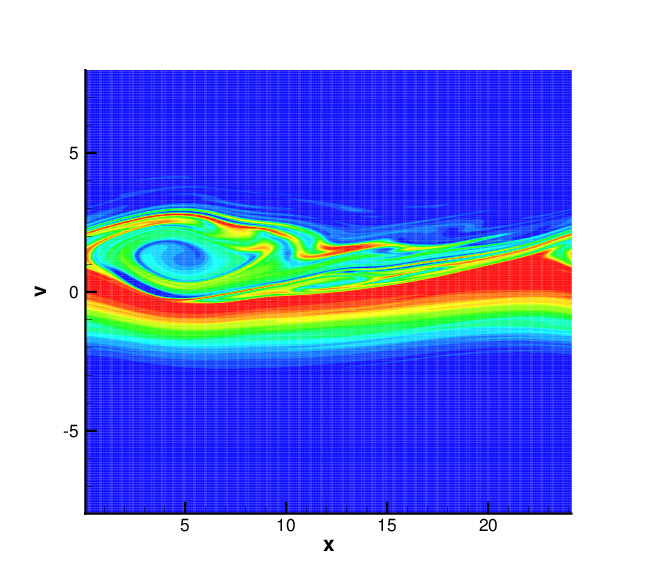}}
                \subfigure[ $t=120$]{\includegraphics[width=2in,angle=0]{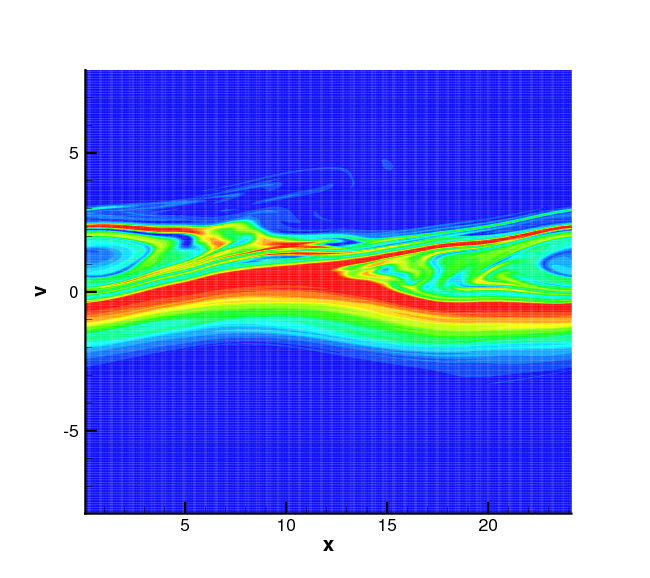}}
    \subfigure[ $t=135$]{\includegraphics[width=2in,angle=0]{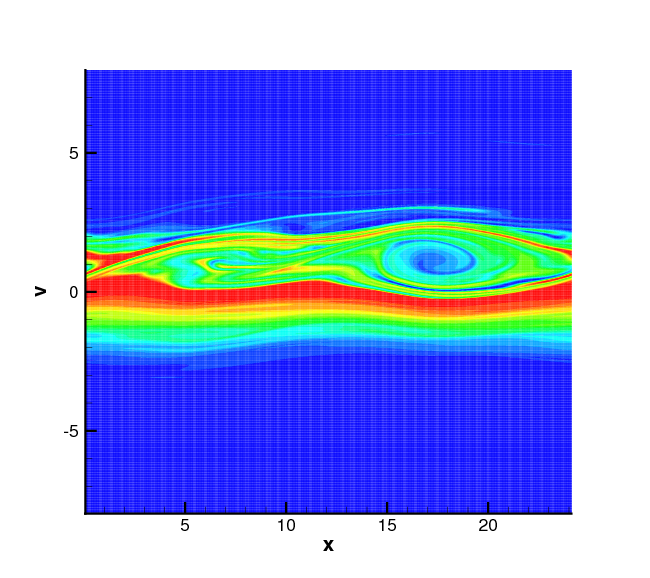}}
        \subfigure[ $t=160$]{\includegraphics[width=2in,angle=0]{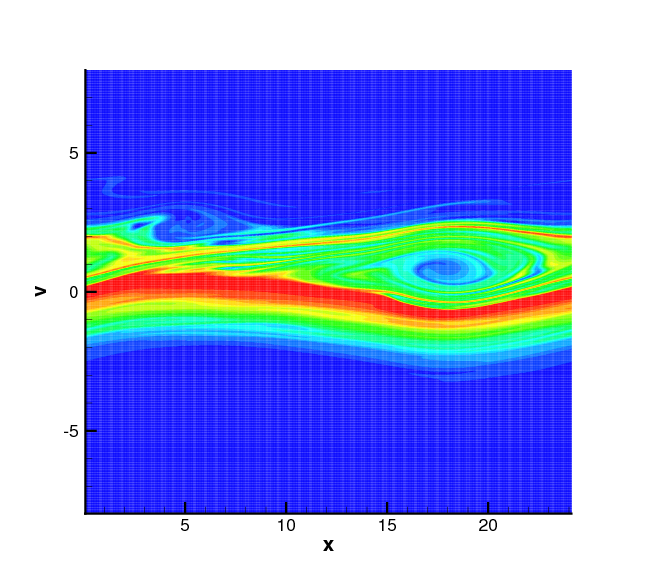}}
                        \subfigure[ $t=225$]{\includegraphics[width=2in,angle=0]{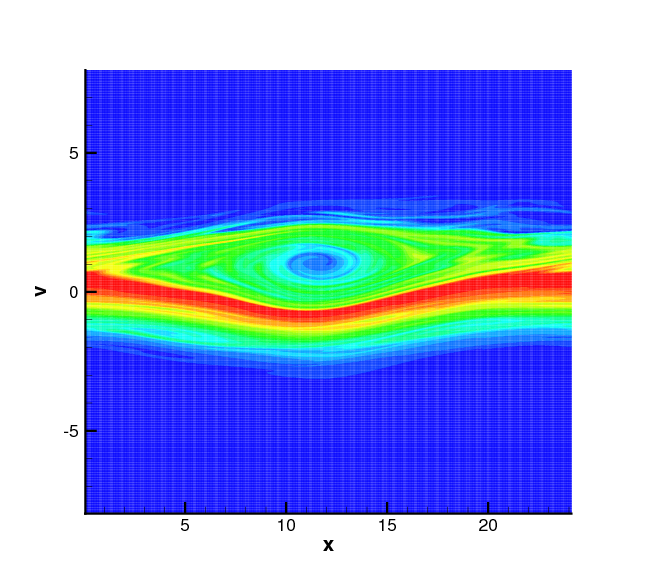}}
                        \subfigure[ $t=300$]{\includegraphics[width=2in,angle=0]{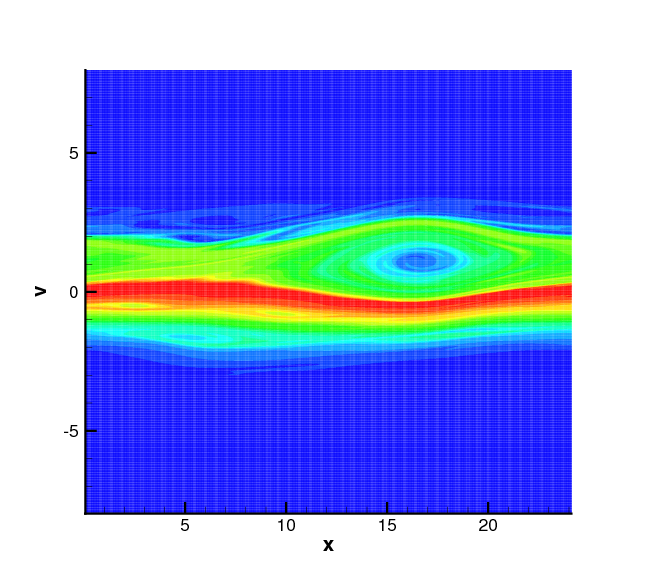}}
  \end{center}
  \caption{Phase space contour plots for the KEEN wave at the times indicated.  A large amplitude drive of  $A_m=0.4$ was used, along with the $\mathbb{P}^2$ basis and  a positivity-preserving limiter on a $200 \times 400$ mesh. }
\label{keencontour}
\end{figure}

Next, we increase the drive to compare with results of \cite{Afeyan_03}.  With the stronger drive of $A_d^A$ with $A_m=0.4$, the system does not approach a uniformly translating state, but approaches a structure with more complicated time dependence as seen in  the phase contour plots of Fig.~\ref{keencontour}.   These figures are in good agreement with those of \cite{Afeyan_03}.

\begin{figure}[htb]
 \hspace{-.90 in}
\subfigure{\includegraphics[width=8in,angle=0]{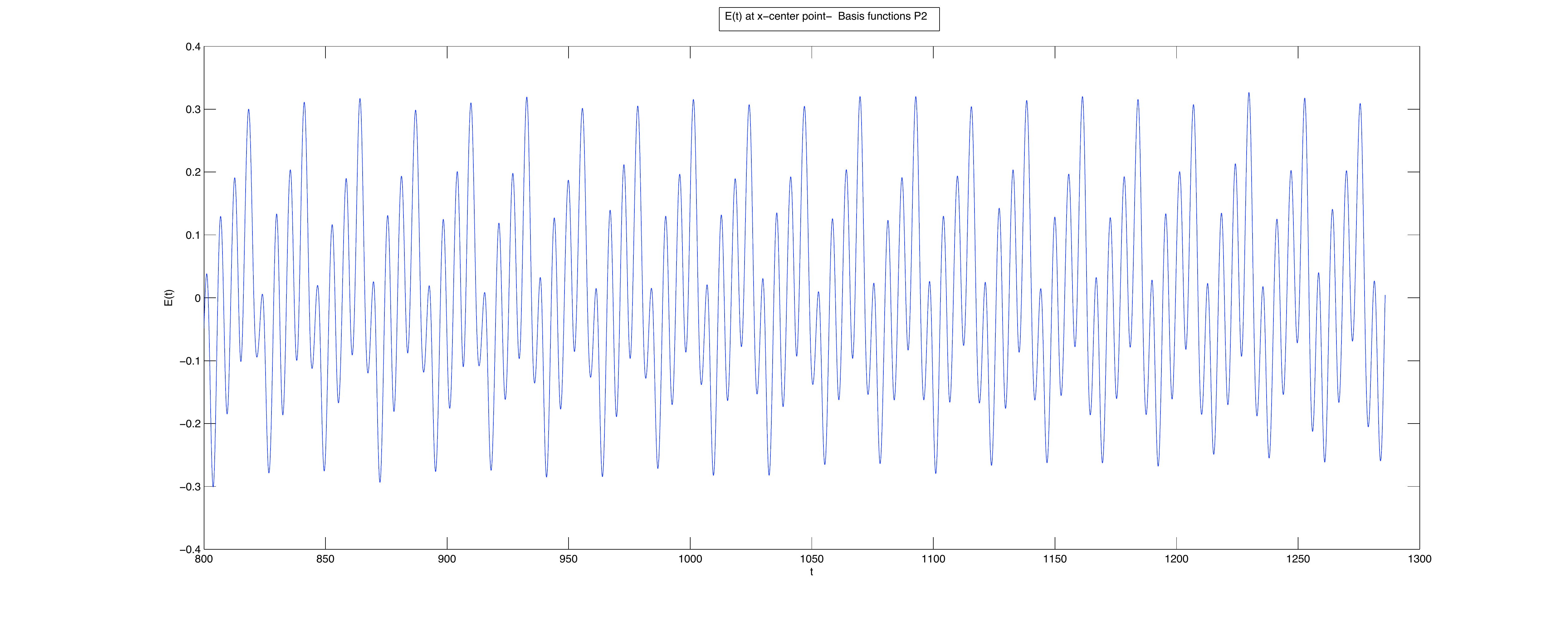}}
       \begin{center}
        \subfigure{\includegraphics[width=2.5in,angle=0]{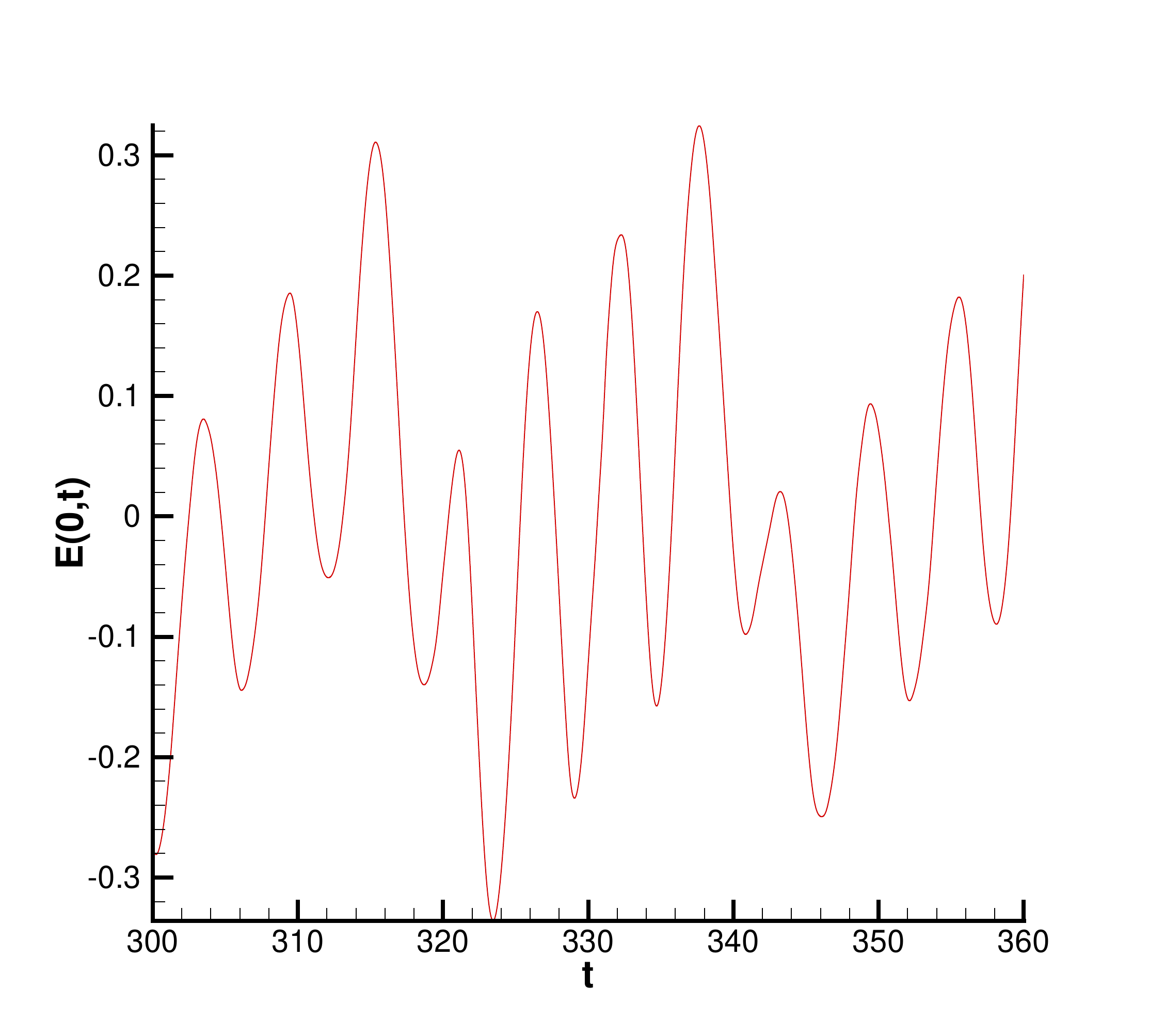}}
  \end{center}
  \caption{{\it(top)} The electric field E(0,t) at the center of the spacial domain at later times for the drive $A_d^{A}$.  The periodic structure is due to multiple interacting BGK-like states. {\it(bottom)} Blow up indicating a period-4 modulation of a central hole such as that of Fig.~\ref{JE}). The simulation was done with $\mathbb{P}^2$ elements with a limiter on a $200 \times 400$ mesh.}
\label{AE}
\end{figure}

The electric field in the middle of the spatial domain, $E(0,t)$, is plotted in Fig.~\ref{AE}, which shows more complicated behavior, which surprisingly heretofore has not  been plotted.   In the top part of this figure we see that there is regular periodic behavior at long times and from the bottom part of the figure we see that there is period-4 modulation of a basic  periodic structure similar to that of Fig.~\ref{JE}.  Closer examination of phase space plots shows that this modulation is cause by the existence of additional smaller BGK-like structures.  We note, that the existence of multiple BGK-like states is not new; for example, they were seen in the simulations of \cite{demeio}.   Thus, we propose that KEEN waves can be interpreted as the interaction of multiple BGK-states, which can also be interpreted as an infinite-dimensional version of Lyapunov-Moser-Weinstein periodic orbits in Hamiltonian systems (see, e.g.\ \cite{moser76}).    This will be the subject of a future publication, so we do not pursue it further here.

\begin{figure}[htb]
  \begin{center}
       \includegraphics[width=3in,angle=0]{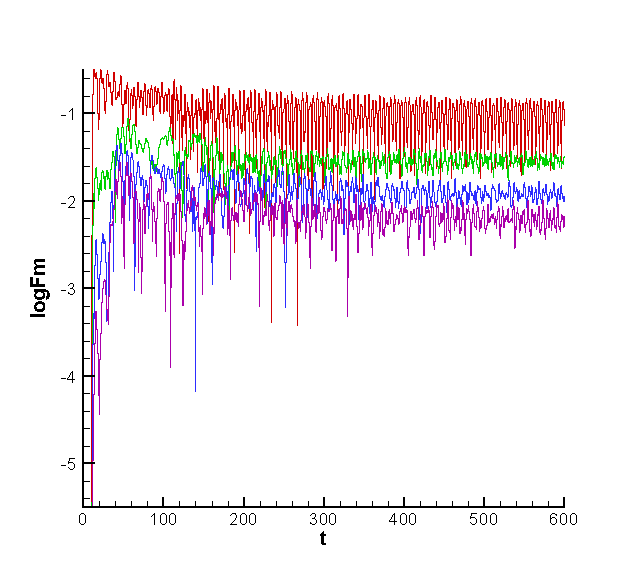}
  \end{center}
  \caption{Evolution of  the first four Log Fourier modes as a function of time for the drive of Eq.~(\ref{rampA}). The simulation used  $\mathbb{P}^2$ elements with a limiter on a $200 \times 400$ mesh.}
\label{keenlf}
\end{figure}

Finally, in Fig.~\ref{keenlf} we see from the evolution of log Fourier modes.  Prior to  $t=10$ the solution  remains roughly at  Maxwellian equilibrium.  However, at around $t=45$ we can observe the formation of the KEEN wave,  which continues to execute the behavior of Fig.~\ref{AE} well after the external field has been turned off at t=60.    We see from this figure the effects of mesh refinement and the use of different polynomial bases, as indicated in the figure.


%


\section{Conclusion}
\label{conclu}

In this paper, we considered the RKDG method for the VP system. We focused on two common solution spaces, viz., those with $\mathbb{P}^l$ and $\mathbb{Q}^l$ elements.  Ignoring  boundary contributions, the scheme can preserve the charge and  momentum, and maintain the total energy up to approximation errors when the polynomial order $l$ is taken big enough. However, when the positivity-preserving limiter was  used, some examples gave relatively large errors for the total energy.  A rigorous study of    numerical recurrence was performed  for the $\mathbb{Q}^l$ elements, and the eigenvalues of the amplification matrix were explicitly obtained. DG schemes of higher order were shown numerically  to give  a recurrence time that is close to the classical calculation $T_R=\frac{2 \pi}{k \triangle v}$. The qualitative behaviors of  the $\mathbb{P}^l$ and $\mathbb{Q}^l$ spaces were similar for most computational examples, except  the linear energy $H_L$  was much better conserved using  the $\mathbb{Q}^l$ space. The  schemes were used to compute the test cases of  Landau damping, the two-stream instability and the KEEN wave, and results  comparable to  those in the literature were obtained.

\section*{Acknowledgments}
 YC was supported by grant NSF DMS-1016001, IMG was supported by grants NSF DMS-0807712 and  DMS-0757450, and  PJM  was  supported by U.S. Dept.\ of Energy Contract \# DE-FG05-80ET-53088.   Also, support from the Department of Mathematics at Michigan State University and the Institute of Computational Engineering and Sciences at the University of Texas Austin are gratefully acknowledged.

\bibliographystyle{abbrv}

\begin{thebibliography}{10}

\bibitem{Afeyan_03}
B.~Afeyan, K.~Won, V.~Savchenko, T.~Johnston, A.~Ghizzo, and P.~Bertrand.
\newblock {Kinetic Electrostatic Electron Nonlinear (KEEN) Waves and their
  Interactions Driven by the Ponderomotive Force of Crossing Laser Beams}.
\newblock {\em Proc. IFSA 2003}, 213, 2003.

\bibitem{Ayuso2010}
B.~Ayuso, J.~A. Carrillo, and C.-W. Shu.
\newblock Discontinuous {Galerkin} methods for the multi-dimensional
  {Vlasov-Poisson} problems.
\newblock {\em Mathematical Models and Methods in Applied Sciences}.
\newblock to appear.

\bibitem{Ayuso2009}
B.~Ayuso, J.~A. Carrillo, and C.-W. Shu.
\newblock Discontinuous {Galerkin} methods for the one-dimensional
  {Vlasov-Poisson} system.
\newblock {\em Kinetic and Related Models}, 4:955--989, 2011.

\bibitem{barnes_86}
J.~Barnes and P.~Hut.
\newblock A hierarchical o(n log n) force-calculation algorithm.
\newblock {\em Nature}, 324:446--449, 1986.

\bibitem{BGK}
I.~Bernstein, J.~M. Greene, and M.~D. Kruskal.
\newblock Exact nonlinear plasma oscillations.
\newblock {\em Phys. Rev.}, 108:546--550, 1957.

\bibitem{Birdsall_book1991}
C.~K. Birdsall and A.~B. Langdon.
\newblock {\em Plasma physics via computer simulation}.
\newblock Institute of Physics Publishing, 1991.

\bibitem{Boris_1976}
J.~Boris and D.~Book.
\newblock Solution of continuity equations by the method of flux-corrected
  transport.
\newblock {\em J. Comp. Phys.}, 20:397--431, 1976.

\bibitem{chengknorr_76}
C.~Z. Cheng and G.~Knorr.
\newblock {The integration of the Vlasov equation in configuration space}.
\newblock {\em Journal of Computational Physics}, 22(3):330--351, 1976.

\bibitem{Cheng_jeans}
Y.~Cheng and I.~M. Gamba.
\newblock {Numerical study of {V}lasov-{P}oisson equations for infinite
  homogeneous stellar systems}.
\newblock {\em Comm. Nonlin. Sci. Num. Sim.}, 17, 2012.

\bibitem{cheng_mathcomp_10}
Y.~Cheng, I.~M. Gamba, and J.~Proft.
\newblock {Positivity-preserving discontinuous Galerkin schemes for linear
  Vlasov-Boltzmann transport equations}.
\newblock {\em Math. Comp.}, 2010.
\newblock to appear.

\bibitem{Cockburn_1990_MC_RK_DG}
B.~Cockburn, S.~Hou, and C.-W. Shu.
\newblock The {Runge-Kutta} local projection discontinuous {Galerkin} finite
  element method for conservation laws {IV}: the multidimensional case.
\newblock {\em Math. Comput.}, 54:545--581, 1990.

\bibitem{Cockburn_1989_JCP_RK_localDG}
B.~Cockburn, S.~Y. Lin, and C.-W. Shu.
\newblock {TVB Runge-Kutta} local projection discontinuous {Galerkin} finite
  element method for conservation laws {III}: one dimensional systems.
\newblock {\em J. Comput. Phys.}, 84:90--113, 1989.

\bibitem{Cockburn_1989_MC_RK_DG}
B.~Cockburn and C.-W. Shu.
\newblock {TVB Runge-Kutta} local projection discontinuous {Galerkin} finite
  element method for conservation laws {II}: general framework.
\newblock {\em Math. Comput.}, 52:411--435, 1989.

\bibitem{Cockburn_1991_MMNA_RK}
B.~Cockburn and C.-W. Shu.
\newblock The {Runge-Kutta} local projection p1-discontinuous {Galerkin} finite
  element method for scalar conservation laws.
\newblock {\em Math. Model. Num. Anal.}, 25:337--361, 1991.

\bibitem{Cockburn_1998_JCP_RK}
B.~Cockburn and C.-W. Shu.
\newblock The {Runge-Kutta} discontinuous {Galerkin} method for conservation
  laws {V}: multidimensional systems.
\newblock {\em J. Comput. Phys.}, 141:199--224, 1998.

\bibitem{demeio}
L.~Demeio and P.~F. Zweifel.
\newblock Numerical simulations of perturbed {V}lasov equilibria.
\newblock {\em Phys. Fluids B}, 2:1252--1254, 1990.

\bibitem{shadwick}
E.~G. Evstatiev and B.~A. Shadwick.
\newblock {\em J. Comp. Phys.}, preprint 2012.
\newblock to appear.

\bibitem{Fijalkow_1999}
E.~Fijalkow.
\newblock A numerical solution to the {Vlasov} equation.
\newblock {\em Comput. Phys. Comm.}, 116:319--328, 1999.

\bibitem{Filbet_PFC_2001}
F.~Filbet, E.~Sonnendr{\" u}cker, and P.~Bertrand.
\newblock Conservative numerical schemes for the {Vlasov} equation.
\newblock {\em J. Comp. Phys.}, 172:166--187, 2001.

\bibitem{fried}
B.~D. Fried and S.~D. Conte.
\newblock {\em The plasma dispersion function}.
\newblock Academic Press, London, 1961.

\bibitem{glassey}
R.~T. Glassey.
\newblock {\em The Cauchy problem in kinetic theory}.
\newblock Society for Industrial and Applied Mathematics (SIAM), Philadelphia,
  PA, 1996.

\bibitem{Heath_thesis}
R.~E. Heath.
\newblock {Numerical analysis of the discontinuous Galerkin method applied to
  plasma physics}.
\newblock 2007.
\newblock Ph. D. dissertation, the University of Texas at Austin.

\bibitem{Heath}
R.~E. Heath, I.~M. Gamba, P.~J. Morrison, and C.~Michler.
\newblock A discontinuous {Galerkin} method for the {Vlasov-Poisson} system.
\newblock {\em J. Comp. Phys.}, 231:1140--1174, 2012.

\bibitem{Hockney_book1981}
R.~W. Hockney and J.~W. Eastwood.
\newblock {\em Computer simulation using particles}.
\newblock McGraw-Hill, New York, 1981.

\bibitem{Johnston_09}
T.~W. Johnston, Y.~Tyshetskiy, A.~Ghizzo, and P.~Bertrand.
\newblock {Persistent subplasma-frequency kinetic electrostatic electron
  nonlinear waves}.
\newblock {\em Phys. Plasmas}, 16:042105, 2009.

\bibitem{jung}
S.~Jung, P.~J. Morrison, and H.~L. Swinney.
\newblock On the statistical mechanics of two-dimensional turbulence.
\newblock {\em J. Fluid Mech.}, 554:433--456, 2006.

\bibitem{Klimas_1987}
A.~J. Klimas.
\newblock A method for overcoming the velocity space filamentation problem in
  collisionless plasma model solutions.
\newblock {\em J. Comp. Phys.}, 68:202--226, 1987.

\bibitem{Klimas_1994}
A.~J. Klimas and W.~M. Farrell.
\newblock A splitting algorithm for {Vlasov} simulation with filamentation
  filtration.
\newblock {\em J. Comp. Phys.}, 110:150--163, 1994.

\bibitem{kraichnan80}
R.~H. Kraichnan and D.~Montgomery.
\newblock Two-dimensional turbulence.
\newblock {\em Rep. Prog. Phys.}, 43:548--618, 1980.

\bibitem{KO58}
M.~D. Kruskal and C.~Oberman.
\newblock On the stability of plasma in static equilibrium.
\newblock {\em Phys. Fluids}, 1:275--280, 1958.

\bibitem{lee}
T.~D. Lee.
\newblock {On some statistical properties of hydrodynamical and
  magneto-hydrodynamical fields}.
\newblock {\em Q. Appl. Math.}, 10:69--74, 1952.

\bibitem{Lesaint_r_74}
P.~Lesaint and P.-A. Raviart.
\newblock On a finite element method for solving the neutron transport
  equation.
\newblock In {\em Mathematical aspects of finite elements in partial
  differential equations (Proc. Sympos., Math. Res. Center, Univ. Wisconsin,
  Madison, Wis., 1974)}, pages 89--123. Math. Res. Center, Univ. of
  Wisconsin-Madison, Academic Press, New York, 1974.

\bibitem{Mont}
S.~Montgomery, J.~A. Cobble, J.~C. Fern‡ndez, R.~J. Focia, R.~P. Johnson,
  N.~Renard-LeGalloudec, H.~A. Rose, and D.~A. Russell.
\newblock Recent {T}rident single hot spot experiments: Evidence for kinetic
  effects, and observation of {L}angmuir decay instability cascade.
\newblock {\em Phys. Plasmas}, 9:2311--2320, 2002.

\bibitem{morrison_98}
P.~J. Morrison.
\newblock {Hamiltonian description of the ideal fluid}.
\newblock {\em Rev. Mod. Phys.}, 70:467--521, 1998.

\bibitem{Morrison_00}
P.~J. Morrison.
\newblock {Hamiltonian description of Vlasov dynamics: action-angle variables
  for the continuous spectrum}.
\newblock {\em Transport Theory and Statistical Physics}, 29:397--414, 2000.

\bibitem{Morrison_89}
P.~J. Morrison and D.~Pfirsch.
\newblock {Free Energy Expressions for {V}lasov-{M}axwell Equilibria}.
\newblock {\em Phys. Rev.}, 40A:3898--3910, 1989.

\bibitem{Morrison_90}
P.~J. Morrison and D.~Pfirsch.
\newblock {The free energy of {M}axwell-{V}lasov equilibria}.
\newblock {\em Phys. Fluids}, 2B:1105--1113, 1990.

\bibitem{Morrison_92}
P.~J. Morrison and D.~Pfirsch.
\newblock {Dielectric energy versus plasma energy, and Hamiltonian action-angle
  variables for the Vlasov equation}.
\newblock {\em Phys. Fluids}, 4B:3038--3057, 1992.

\bibitem{moser76}
J.~Moser.
\newblock Periodic orbits near an equilibrium and a theorem by {A}lan
  {W}einstein.
\newblock {\em Comm. Pure App. Math.}, 29:727--747, 1976.

\bibitem{burgers}
F.~T.~M. Nieuwstadt and J.~A. Steketee.
\newblock {\em Selected papers of J. M. Burgers}.
\newblock Kluwer Academic Publishers, Dodrecht, 1995.

\bibitem{qiu_ppdg_11}
J.-M. Qiu and C.-W. Shu.
\newblock {Positivity preserving semi-Lagrangian discontinuous Galerkin
  formulation: theoretical analysis and application to the Vlasov-Poisson
  system}.
\newblock 2011.
\newblock submitted to J. Comp. Phys.

\bibitem{Reed_hill_73}
W.~Reed and T.~Hill.
\newblock Tiangular mesh methods for the neutron transport equation.
\newblock Technical report, Los Alamos National Laboratory, Los Alamos, NM,
  1973.

\bibitem{rossmanith_11}
J.~Rossmanith and D.~Seal.
\newblock {A positivity-preserving high-order semi-Lagrangian discontinuous
  Galerkin scheme for the Vlasov-Poisson equations}.
\newblock 2011.
\newblock submitted to J. Comp. Phys.

\bibitem{Shu_1988_JCP_NonOscill}
C.-W. Shu and S.~Osher.
\newblock Efficient implementation of essentially non-oscillatory
  shock-capturing schemes.
\newblock {\em J. Comput. Phys.}, 77:439--471, 1988.

\bibitem{Sonnendrucker_1999}
E.~Sonnendr{\" u}cker, J.~Roche, P.~Bertrand, and A.~Ghizzo.
\newblock The {semi-Lagrangian} method for the numerical resolution of the
  {Vlasov} equation.
\newblock {\em J. Comp. Phys.}, 149(2):201--220, 1999.

\bibitem{Valentini}
F.~Valentini, T.~M. O'Neil, and D.~H.~E. Dubin.
\newblock Excitation of nonlinear electron acoustic waves.
\newblock {\em Phys. Plasmas}, 13:052303, 2006.

\bibitem{Valentini2}
F.~Valentini, D.~Perrone, F.~Califano, F.~Pegoraro, P.~Veltri, P.~J. Morrison,
  and T.~M. O'Neil.
\newblock Undamped electrostatic plasma waves.
\newblock {\em Phys. Plasmas}, 19:092103, 2012.

\bibitem{Zhang_Water_2010}
Y.~Xing, X.~Zhang, and C.-W. Shu.
\newblock {Positivity preserving high order well balanced discontinuous
  Galerkin methods for the shallow water equations}.
\newblock {\em Advances in Water Resources}, 33:1476--1493, 2010.

\bibitem{Zaki_1988_1}
S.~Zaki, L.~Gardner, and T.~Boyd.
\newblock A finite element code for the simulation of one-dimensional {Vlasov}
  plasmas. i. theory.
\newblock {\em J. Comp. Phys.}, 79:184--199, 1988.

\bibitem{Zaki_1988_2}
S.~Zaki, L.~Gardner, and T.~Boyd.
\newblock A finite element code for the simulation of one-dimensional {Vlasov}
  plasmas. ii. applications.
\newblock {\em J. Comp. Phys.}, 79:200--208, 1988.

\bibitem{Zhang_2005_CF_DG}
M.~Zhang and C.-W. Shu.
\newblock An analysis of and a comparison between the discontinuous {Galerkin}
  and the spectral finite volume methods.
\newblock {\em Computers and Fluids}, 34:581--592, 2005.

\bibitem{Zhang_2010_Max}
X.~Zhang and C.-W. Shu.
\newblock On maximum-principle-satisfying high order schemes for scalar
  conservation laws.
\newblock {\em J. Comput. Phys.}, 229:3091--3120, 2010.

\bibitem{Zhang_JCP_2010_Euler}
X.~Zhang and C.-W. Shu.
\newblock {On positivity preserving high order discontinuous Galerkin schemes
  for compressible Euler equations on rectangular meshes}.
\newblock {\em J. Comput. Phys.}, 229:8918--8934, 2010.

\bibitem{Zhang_mpp_review}
X.~Zhang and C.-W. Shu.
\newblock {Maximum-principle-satisfying and positivity-preserving high order
  schemes for conservation laws: Survey and new developments}.
\newblock {\em Proceedings of the Royal Society A}, 2011.
\newblock to appear.

\bibitem{Zhang_JCP_2011_Euler_source}
X.~Zhang and C.-W. Shu.
\newblock {Positivity-preserving high order discontinuous Galerkin schemes for
  compressible Euler equations with source terms}.
\newblock {\em J. Comput. Phys.}, 230:1238--1248, 2011.

\bibitem{Zhang_JSC_triangle}
X.~Zhang, Y.~Xia, and C.-W. Shu.
\newblock {Maximum-principle-satisfying and positivity-preserving high order
  discontinuous Galerkin schemes for conservation laws on triangular meshes}.
\newblock {\em J. Sci. Comp.}
\newblock to appear.

\bibitem{Zhong_11}
X.~Zhong and C.-W. Shu.
\newblock Numerical resolution of discontinuous galerkin methods for time
  dependent wave equations.
\newblock 200:2814--2827, 2011.

\bibitem{Guo_landau_2001}
T.~Zhou, Y.~Guo, and C.-W. Shu.
\newblock {Numerical study on Landau damping}.
\newblock {\em Physica D: Nonlinear Phenomena}, 157(4):322--333, 2001.

\end{thebibliography}


\end{document}